\documentclass[a4paper,11pt,final]{article}

\usepackage{graphicx,psfrag}
\usepackage{amssymb,amsmath,amsthm}
\usepackage{enumerate}
\usepackage{a4wide}
\usepackage{hyperref}

\numberwithin{equation}{section}
\allowdisplaybreaks[4]

\newtheorem{proposition}{Proposition}[section]

\newtheorem{lemma}[proposition]{Lemma}

\newtheorem{definition}[proposition]{Definition}
\theoremstyle{definition}
\newtheorem{remark}[proposition]{Remark}

\newenvironment{proofof}[1]{\smallskip\noindent{\textbf{Proof~of~#1.}}%
  \hspace{1pt}}{\hspace{-5pt}{\nobreak\quad\nobreak\hfill\nobreak%
    $\square$\vspace{2pt}\par}\smallskip\goodbreak}

\newcommand{\C}[1]{\mathbf{C}^{#1}}

\newcommand{\caratt}[1]{\chi_{\strut #1}}

\renewcommand{\L}[1]{{\mathbf{L}^#1}}
\newcommand{\Lloc}[1]{{\mathbf{L}_{loc}^{#1}}}

\newcommand{\W}[2]{{\mathbf{W}^{#1,#2}}}

\newcommand{\modulo}[1]{{\left|#1\right|}}
\newcommand{\norma}[1]{{\left\|#1\right\|}}
\newcommand{\reali}{{\mathbb{R}}}
\newcommand{\naturali}{{\mathbb{N}}}
\newcommand{\interi}{{\mathbb{Z}}}
\renewcommand{\epsilon}{\varepsilon}
\renewcommand{\phi}{\varphi}
\renewcommand{\theta}{\vartheta}

\newcommand{\tv}{\mathinner{\rm TV}}
\newcommand{\spt}{\mathop{\rm spt}}
\newcommand{\sgn}{\mathop{\rm sgn}}

\renewcommand{\d}[1]{\mathinner{\mathrm{d}{#1}}}

\newcommand \be           {\begin{equation}}
\newcommand \ee            {\end{equation}}

\begin{document}

\title{A Numerical Approach to \\ Scalar Nonlocal Conservation Laws}

\author{Paulo Amorim$^1$ \and Rinaldo M.~Colombo$^2$ \and Andreia
  Teixeira$^1$}

\footnotetext[1]{Centro de Matem\'atica e Aplica\c c\~oes
  Fundamentais, Departamento de Matem\'atica, Universidade de Lisboa,
  Av. Prof. Gama Pinto 2, 1649-003 Lisboa, Portugal.}

\footnotetext[2]{Unit\`a INdAM, Universit\`a di Brescia, Via Branze
  38, 25123 Brescia, Italy.}

\maketitle

\begin{abstract}

  \noindent We address the study of a class of 1D nonlocal
  conservation laws from a numerical point of view. First, we present
  an algorithm to numerically integrate them and prove its
  convergence. Then, we use this algorithm to investigate various
  analytical properties, obtaining evidence that usual properties of
  standard conservation laws fail in the nonlocal setting. Moreover,
  on the basis of our numerical integrations, we are lead to
  conjecture the convergence of the nonlocal equation to the local
  ones, although no analytical results are, to our knowledge,
  available in this context.

  \medskip

  \noindent\textit{2000~Mathematics Subject Classification:} 35L65

  \medskip

\end{abstract}

\section{Introduction}
\label{sec:I}

Conservation laws with nonlocal fluxes have appeared recently in the
literature, arising naturally in many fields of application, such as
in crowd dynamics (see~\cite{ColomboGaravelloLecureux2012,
  ColomboHertyMercier, ColomboMercierXDafermos, PiccoliTosin} and the
references therein), or in models inspired from biology,
see~\cite{CarrilloColomboGwiazdaUlikowska, Eftimie2012, Friedman2012,
  GwiazdaLorenzEtAl}.

In this paper, we initiate the study of these
equations from a numerical point of view. First, we prove the
convergence of a finite volume algorithm to numerically integrate a class of
one-dimensional conservation laws with a nonlocal flow. Then, we use this algorithm to
show peculiar properties of these nonlocal equations and, in
particular, how they differ from the usual local ones.

Consider the scalar equation
\begin{equation}
  \label{eq:1}
  \left\{
    \begin{array}{l}
      \partial_t \rho
      +
      \partial_x \left(f(t,x,\rho) \, v (\rho * \eta)\right) =0
      \\
      \rho (0,x) = \rho^o (x)
    \end{array}
  \right.
  \qquad (t,x) \in \reali^+ \times \reali
\end{equation}
which slightly extends, in the 1D case, the class of equations
considered in~\cite{ColomboHertyMercier, ColomboMercierXDafermos}. We
present below a numerical scheme to integrate~\eqref{eq:1} and prove
its convergence. As a byproduct, we also establish an existence result
for~\eqref{eq:1}, thus slightly
extending~\cite[Theorem~2.2]{ColomboHertyMercier} in the 1D case.

\smallskip

This numerical algorithm is then implemented and used to investigate
various properties of~\eqref{eq:1}. First, we provide evidence that
the usual \emph{Maximum Principle} for scalar conservation laws fails
in the case of~\eqref{eq:1}. Another integration shows that the total
variation of the solution to~\eqref{eq:1} may well sharply increase,
contrary to what happens in the standard local situation. Remark that
both these examples are in agreement with the estimates we rigorously
obtain on the approximate solutions.

Of particular interest is the limit $\eta \to \delta$, $\delta$ being
the Dirac measure centered at the origin. Numerical integrations show
that the solutions to~\eqref{eq:1} converge to that of
\begin{equation}
  \label{eq:1loc}
  \left\{
    \begin{array}{l}
      \partial_t \rho
      +
      \partial_x \left(f(t,x,\rho) \, v (\rho)\right) =0
      \\
      \rho (0,x) = \rho^o (x)
    \end{array}
  \right.
  \qquad (t,x) \in \reali^+ \times \reali
\end{equation}
although no rigorous proof of this convergence is, to our knowledge,
known. Remark that in the nonlocal case, well posedness results are
available also in the case of systems in several space dimensions,
see~\cite{ColomboMercierXDafermos, CrippaLecureux}. Hence, the ability
of passing to the limit $\eta \to \delta$ might help in the search for
analytical results about systems of conservation laws in several space
dimensions.

\smallskip
Let us make the following remarks. First, the scheme presented below has
an associated CFL condition.
The CFL condition is often interpreted through a comparison between
the numerical propagation speed and the analytical one, see for
instance~\cite[\S~4.4, p.~68]{LeVequeBook2002}. In the present
nonlocal case~\eqref{eq:1}, information propagates at an infinite
speed, due to the presence of the term $\eta * \rho$. Nevertheless,
also in the nonlocal case~\eqref{eq:1} a suitable CFL condition plays
a key role, see~\eqref{eq:CFL}.

Second, the scheme presented below is not \emph{monotone} in the sense of the
usual definition~\cite[Formula~(12.42)]{LeVequeBook2002}, as follows
from the integration in Section~\ref{subs:TV}. There, both constant
initial data $\bar \rho = 0$ and $\bar \rho = 1$ yield constant
solutions, but the initial datum~\eqref{eq:idtv}, although it
attains values in $[0,1]$, yields a solution exceeding $1$.
Nevertheless, the scheme~\eqref{LFa} enjoys several properties of
monotone schemes, proved in the lemmas in Section~\ref{sec:MR}.

\begin{remark}
  Throughout this work, we follow the usual habit of referring
  to~\eqref{eq:1} as to a \emph{nonlocal} equation and, hence, to the
  standard case~\eqref{eq:1loc} as to the \emph{local} case. However,
  whenever the support of $\eta$ is bounded, it might seem more
  appropriate to call~\eqref{eq:1} a \emph{local} equation
  and~\eqref{eq:1loc} the \emph{punctual} case.
\end{remark}

The next section deals with the definition of the algorithm and with
the statement of the estimates which
ensure its convergence, as well as the entropicity of the
limit solution. Section~\ref{sec:NI} deals with various numerical
integrations of~\eqref{eq:1}. All proofs are deferred to the last
Section~\ref{sec:TD}.

\section{Main Results}
\label{sec:MR}

Throughout, we set $\reali^+ = \left[ 0, +\infty\right[$.

As a starting point, we state what we mean by \emph{solution}
to~\eqref{eq:1}, see also~\cite[Definition~2.1]{ColomboHertyMercier}.

\begin{definition}
  \label{def:sol1}
  Let $T > 0$. Fix $\rho^o \in \L\infty (\reali; \reali)$. A weak
  entropy solution to~\eqref{eq:1} on $[0, T]$ is a bounded measurable
  Kru\v zkov solution $\rho \in \C0 \left( [0, T]; \Lloc1(\reali;
    \reali) \right)$ to
  \begin{displaymath}
    \left\{
      \begin{array}{@{\,}l@{}}
        \partial_t \rho + \partial_x \left( f(t,x,\rho) \, V(t,x) \right) = 0
        \\
        \rho(0,x) = \rho^o(x)
      \end{array}
    \right.
    \quad \mbox{ where} \quad
    V(t,x) =
    v \! \left( (\rho (t) * \eta) (x) \right)\,.
  \end{displaymath}
\end{definition}

\noindent For the definition of Kru\v zkov solution, see for
instance~\cite[Paragraph~6.2]{DafermosBook}
or~\cite[Definition~1]{Kruzkov}.  Here, as usual,
\begin{displaymath}
  \left(\rho (t) * \eta\right) (x)
  =
  \int_{\reali} \rho (t,\xi) \; \eta (x-\xi) \, \d\xi\,.
\end{displaymath}
Remark that the assumptions
\begin{eqnarray}
  \label{eq:hyp1}
  f \in \C2(\reali^+ \times \reali \times \reali; \reali)
  &\mbox{ and }&
  \left\{
    \begin{array}{rcl}
      \displaystyle
      \sup_{t,x,\rho} \modulo{\partial_\rho f (t,x,\rho)}
      & < &
      +\infty
      \\
      \displaystyle
      \sup_{t,x} \modulo{\partial_x f (t,x,\rho)}
      & < &
      C \,\modulo{\rho}
      \\
      \displaystyle
      \sup_{t,x} \modulo{\partial^2_{xx} f (t,x,\rho)}
      & < &
      C \,\modulo{\rho}
      \\
      \displaystyle
      \forall\,t, x \quad f (t,x,0)
      & = &
      0
    \end{array}
  \right.
  \\
  \label{eq:hyp2}
  v \in (\C2 \cap \W1\infty)(\reali; \reali)
  &\mbox{ and } &\qquad
  \eta \in  (\C2 \cap \W2\infty) (\reali; \reali)
\end{eqnarray}
ensure that the transport equation in Definition~\ref{def:sol1} fits
in Kru\v zkov framework, see~\cite{DafermosBook, Kruzkov}.  From the
modeling point of view, it is natural to require that the kernel
$\eta$ attains only positive (or non-negative) values. However, this
requirement is not necessary for the analytical results below.

Below, Remark~\ref{rem:MP} and Lemma~\ref{lem:Linfty} provide uniform
$\L\infty$ bounds on the solution to~\eqref{eq:1} under
conditions~\eqref{eq:hyp1}--\eqref{eq:hyp2} on the equations and for
data in $\L\infty$. Therefore, the apparently strong requirement
$\norma{\partial_\rho f}_{\L\infty} < +\infty$ can be easily relaxed
to
\begin{displaymath}
  \sup_{t \in \reali^+,\, x \in \reali, \, \rho \in [-M, M]}
  \modulo{\partial_\rho f (t, x, \rho)}
  <
  +\infty
\end{displaymath}
for a suitable positive $M$.  Moreover, the usual sublinearity
condition $\sup_{t,x} \modulo{\partial_x f (t,x,\rho)} < C( 1 +
\modulo{\rho})$ takes the form $\sup_{t,x} \modulo{\partial_x f
  (t,x,\rho)} < C \modulo{\rho}$ in~\eqref{eq:hyp1} due to the
assumption $f (t,x,0) = 0$ for all $t$ and $x$.

Introduce a uniform mesh with size $h$ along the $x$ axis and size
$\tau$ along the $t$ axis. Throughout, we assume that
\begin{equation}
  \label{eq:lambda}
  h < 1/C
\end{equation}
with $C$ as in~\eqref{eq:hyp1} and that the following \emph{CFL
  condition} is satisfied:
\begin{equation}
  \label{eq:CFL}
  \lambda \,
  \left(
    1
    +
    2 \, \norma{\partial_\rho f}_{\L\infty}
  \right) \norma{v}_{\L\infty}
  \leq
  \frac{1}{6}
\end{equation}
where, as usual, $\lambda = \tau / h$.

Consider the following Lax--Friedrichs type scheme:
\begin{equation}
  \label{LF}
  \left\{
    \begin{array}{lcl}
      \rho^{n+1}_j
      & = &
      \rho^n_j
      -
      \lambda
      \left(
        \mathbf{f}^n_{j+1/2}(\rho^n_j, \rho^n_{j+1})
        -
        \mathbf{f}^n_{j-1/2}(\rho^n_{j-1}, \rho^n_{j})
      \right)
      \\[6pt]
      \rho^o_j
      & = &
      \displaystyle
      \frac{1}{h} \int_{x_{j-1/2}}^{x_{j+1/2}} \rho^o (x) \d x
    \end{array}
  \right.
\end{equation}
where the numerical flux $\mathbf{f}^n_{j+1/2}$ in~\eqref{LF} is given
by
\begin{equation}
  \label{LFa}
  \mathbf{f}^n_{j+1/2} (\rho_1,\rho_2)
  :=
  \frac{f(t^n, x_{j+1/2},\rho_1)
    +
    f(t^n, x_{j+1/2},\rho_2)}{2}
  \,
  v (c^n_{j+1/2})
  -
  \frac{1}{6\,\lambda} (\rho_2 - \rho_1) \,.
\end{equation}
Here, the convolution is computed through a standard quadrature
formula using the same space mesh, as follows
\begin{equation}
  \label{eq:cnj}
  c^n_{j+1/2}
  =
  \sum_{k\in\interi} h \, \rho^n_{k+1/2} \, \eta_{j+1/2-k}
\end{equation}
where $\rho^n_{k+1/2}$ is a suitable convex combination of $\rho^n_k$,
$\rho^n_{k+1}$ and $\eta_{j+1/2} = \frac{1}{h} \int_{x_j}^{x_{j+1}}
\eta (x) \d x$, for instance.

\medskip

The next three lemmas provide the basic properties of the
algorithm~(\ref{LF}), namely positivity, $\L1$ and $\L\infty$
bounds. All proofs are deferred to Section~\ref{sec:TD}.

\begin{lemma}[Positivity]
  \label{lem:pos}
  Let conditions~\eqref{eq:hyp1}--\eqref{eq:hyp2} hold. Assume that
  $h$ and $\tau$ satisfy~\eqref{eq:lambda} and the CFL
  condition~\eqref{eq:CFL}.  If $\rho^o_j \geq 0$ for all $j$, then
  the approximate solution constructed by the algorithm~\eqref{LF} is
  such that $\rho^n_j \geq 0$ for all $j$ and $n$.
\end{lemma}

\begin{remark}
  \label{rem:MP}
  The proof of the above lemma clearly shows that if we assume $\rho^o
  \leq 0$, then $\rho^n \leq 0$ for all $n$. Moreover, under the same
  assumptions~\eqref{eq:hyp1}--\eqref{eq:hyp2}--\eqref{eq:lambda}--\eqref{eq:CFL},
  a straightforward modification of the proof of Lemma~\ref{lem:pos}
  ensures that if there exists a $\bar \rho \in \reali^+$ such that $f
  (t,x,\bar\rho) =0$, then the inequality $\rho^o \geq \bar \rho$,
  respectively $\rho^o \leq \bar \rho$, implies that $\rho^n \geq
  \bar\rho$, respectively $\rho^n \leq \bar \rho$, for all $n$.
\end{remark}

\begin{lemma}[$\L1$ bound]
  \label{lem:L1}
  Let conditions~\eqref{eq:hyp1}--\eqref{eq:hyp2} hold. Assume that
  $h$ and $\tau$ satisfy~\eqref{eq:lambda} and the CFL
  condition~\eqref{eq:CFL}.  If $\rho^o_j \geq 0$ for all $j$, then
  the approximate solution constructed by the algorithm~\eqref{LF}
  satisfies
  \begin{displaymath}
    \norma{\rho^n}_{\L1} \leq \norma{\rho^o}_{ \L1} \,.
  \end{displaymath}
\end{lemma}

\begin{lemma}[$\L\infty$ bound]
  \label{lem:Linfty}
  Let conditions~\eqref{eq:hyp1}--\eqref{eq:hyp2} hold. Assume that
  $h$ and $\tau$ satisfy~\eqref{eq:lambda} and the CFL
  condition~\eqref{eq:CFL}. If $\rho^o_j \geq 0$ for all $j$, then the
  solution constructed by the algorithm~\eqref{LF} satisfies
  \begin{displaymath}
    \norma{\rho^{n}}_{\L\infty}
    \leq
    \norma{\rho^o}_{\L\infty} \, e^{\mathcal{L} t} \,,
  \end{displaymath}
  where $\mathcal{L}$ depends on C in~(\ref{eq:hyp1}), on various
  norms of $f, v, \eta$ and on the $\L1$ norm of the initial datum,
  see~\eqref{eq:L}.
\end{lemma}

The next result concerns the bound on the total variation of the
approximate solution constructed in~(\ref{LF}). In the standard
Kru\v{z}kov case, when the flow is independent from $t$ and $x$, the
total variation of the solution is well know to be a non-increasing
function of time, see~\cite[Theorem~6.1]{BressanLectureNotes}. Here,
on the contrary, the total variation and the $\L\infty$ norm of the
solution to~(\ref{eq:1})
may well sharply increase due to the nonlocal terms, even when
the flow is independent
from $t$ and $x$, see Section~\ref{subs:TV}.

\begin{proposition}[Total variation bound]
  \label{prop:TV}
  Let conditions~\eqref{eq:hyp1}--\eqref{eq:hyp2} hold. Assume that
  $h$ and $\tau$ satisfy~\eqref{eq:lambda} and the CFL
  condition~\eqref{eq:CFL}. If $\rho^o_j \geq 0$ for all $j$, then the
  approximate solution constructed by the algorithm~\eqref{LF}
  satisfies the following total variation estimate, for all $n \geq
  0$:
  \begin{equation}
    \label{TV1}
    \sum_{j\in\interi} \modulo{\rho^n_{j+1} - \rho^n_j}
    \leq
    \left(
      \mathcal{K}_2 \, t
      +
      \sum_{j\in\interi} \modulo{\rho^o_{j+1} - \rho^o_j}
    \right)
    e^{\mathcal{K}_1 t} \,,
  \end{equation}
  where the constants $\mathcal{K}_1$ and $\mathcal{K}_2$ depend on
  $C$ in~(\ref{eq:hyp1}), on various norms of $f, v, \eta$ and of the
  initial datum, see~\eqref{eq:K1K2}.
\end{proposition}

A first consequence of the bound on the total variation is the
$\L1$--Lipschitz continuity in time of the approximate solution,
proved in the following lemma.

\begin{lemma}[$\L1$-Lipschitz continuity in time]
  \label{lem:DepOnT}
  Fix a positive $T$.  Let conditions~\eqref{eq:hyp1}--\eqref{eq:hyp2}
  hold. Assume that $h$ and $\tau$ satisfy~\eqref{eq:lambda} and the
  CFL condition~\eqref{eq:CFL}. If $\rho^o_j \geq 0$ for all $j$, then
  the approximate solution constructed by the algorithm~\eqref{LF} is
  an $\L1$-Lipschitz continuous function of time, in the sense that
  for any $n,m \in \naturali$ such that $n\, \tau \leq T$ and $m\,\tau
  \leq T$,
  \begin{displaymath}
    \norma{\rho^n - \rho^m}_{\L1}
    \leq
    {\cal C}(T) \, \modulo{n-m} \, \tau
  \end{displaymath}
  where the quantity $\mathcal{C} (T)$ grows exponentially in time and
  depends on $C$ in~(\ref{eq:hyp1}), on various norms of $f, v, \eta$
  and of the initial datum, see~\eqref{eq:C}.
\end{lemma}

The $\L\infty$ bound proved in Lemma~\ref{lem:Linfty}, the total
variation bound proved in Proposition~\ref{prop:TV} and the uniform
continuity in time that follows from Lemma~\ref{lem:DepOnT} allow to
apply Helly Theorem, for instance in the form
of~\cite[Theorem~2.6]{BressanLectureNotes}, to the sequence of
approximate solutions constructed through~(\ref{LF}).  A
straightforward limiting procedure, see for
instance~\cite[Section~6.2]{BressanLectureNotes}, thus ensures the
existence of weak solutions to the Cauchy problem for~(\ref{eq:1}).

To obtain uniqueness, we prove that the approximate
solutions~(\ref{LF}) also satisfy a discrete entropy condition. To
this end, define for each $k\in\reali$ the Kru\v zkov numerical
entropy flux as
\begin{equation}
  \label{eq:Kruzhkov}
  F^k_{j+1/2} (\rho_1, \rho_2)
  =
  \mathbf{f}^n_{j+1/2}(\rho_1 \vee k, \rho_2 \vee k)
  -
  \mathbf{f}^n_{j+1/2}(\rho_1 \wedge k, \rho_2 \wedge k) \,,
\end{equation}
where $a\vee b = \max(a,b)$ and $a \wedge b = \min(a,b)$.

\begin{proposition}[Discrete entropy condition]
  \label{prop:Entropy}
  Let conditions~\eqref{eq:hyp1}--\eqref{eq:hyp2} hold. Assume that
  $h$ and $\tau$ satisfy~\eqref{eq:lambda} and the CFL
  condition~\eqref{eq:CFL}. If $\rho^o_j \geq 0$ for all $j$, then the
  approximate solution constructed by the algorithm~\eqref{LF}
  verifies the discrete entropy inequality
  \begin{equation}
    \label{eq:Entropy}
    \begin{array}{rcl}
      \modulo{\rho^{n+1}_{j} - k}
      -
      \modulo{\rho^{n}_{j} - k}
      +
      \lambda
      \left(
        F^k_{j+1/2} (\rho^{n}_{j}, \rho^{n}_{j+1})
        -
        F^k_{j-1/2}(\rho^{n}_{j-1}, \rho^{n}_{j})
      \right)
      \\[6pt]
      +
      \lambda \,
      \sgn ( \rho^{n+1}_{j} - k)
      \left( f(t^n, x_{j+1/2}, k) - f(t^n, x_{j-1/2}, k) \right)
      & \leq & 0 \,
    \end{array}
  \end{equation}
  for all $k\in\reali$.
\end{proposition}

\section{Numerical Integrations}
\label{sec:NI}

\subsection{A Nonlocal Traffic Model}
\label{subs:LWR}

The classical Lighthill--Whitham~\cite{LighthillWhitham} and
Richards~\cite{Richards} (LWR) model for vehicular traffic consists of
the continuity equation $\partial_t \rho + \partial_x (\rho \, V) =0$
supplied with a suitable \emph{speed law} $V = V (\rho)$. Here, as
usual, $t$ is time, $x$ an abscissa along a rectilinear road with
neither entries nor exits and $\rho \in [0,1]$ is the (average)
vehicular density.

Equation~(\ref{eq:1}) with
\begin{equation}
  \label{eq:Traffic}
  f (\rho) = \rho \, (1-\rho)
  \quad
  v (r) = V_{\max} \, (1-r)
  \quad \mbox{and} \quad
  \eta (x) = \alpha \left((x-a) (b-x)\right)^{5/2} \, \caratt{[a,b]} (x) \,,
\end{equation}
where $V_{\max} >0$, can be used as an LWR-type macroscopic model for
vehicular traffic, where drivers adjust their speed according to the
local traffic density, so that the speed law takes the functional form
\begin{displaymath}
  V (\rho) = V_{\max} \, (1-\rho) \, (1-\rho*\eta) \,.
\end{displaymath}
The coefficient $\alpha$ in~\eqref{eq:Traffic} is chosen so that
$\int_{\reali} \eta = 1$.  The parameters $a$ and $b$ are the
\emph{horizon} of each driver, in the sense that a driver situated at
$x$ adjusts his speed according to the average vehicular density he
sees on the interval $[x-b, x-a]$. To emphasize their roles, we select
below the two situations
\begin{equation}
  \label{eq:6}
  \begin{array}{rcl}
    a & = & -1/4 \\ b & = & 0
  \end{array}
  \quad \mbox{ and } \quad
  \begin{array}{rcl}
    a & = & 0 \\ b & = & 1/4 \,.
  \end{array}
\end{equation}
In the former case, drivers look forward, while in the latter they
look backward. We consider the initial datum
\begin{equation}
  \label{eq:5}
  \rho^o (x) =
  \frac{1}{2} \, \caratt{[-2.8, \, -1.8]} (x)
  +
  \frac{3}{4} \, \caratt{[-1.2, \, -0.2]} (x)
  +
  \frac{3}{4} \, \caratt{[0.6, \, 1.0]} (x)
  +
  \caratt{[1.5, +\infty[} (x)
\end{equation}
representing three groups of vehicles lining up in a queue.

The results in Section~\ref{sec:MR} ensure that for any $\rho^o \in
\L1 (\reali; [0,1])$, the Cauchy problem consisting
of~(\ref{eq:1})--(\ref{eq:Traffic}) with initial datum $\rho^o$ admits
a unique solution $\rho = \rho (t,x)$ attaining values in $[0,1]$.
\begin{figure}[!htpb]
  \centering
  \includegraphics[width=0.2\textwidth, trim=20 0 20
  0]{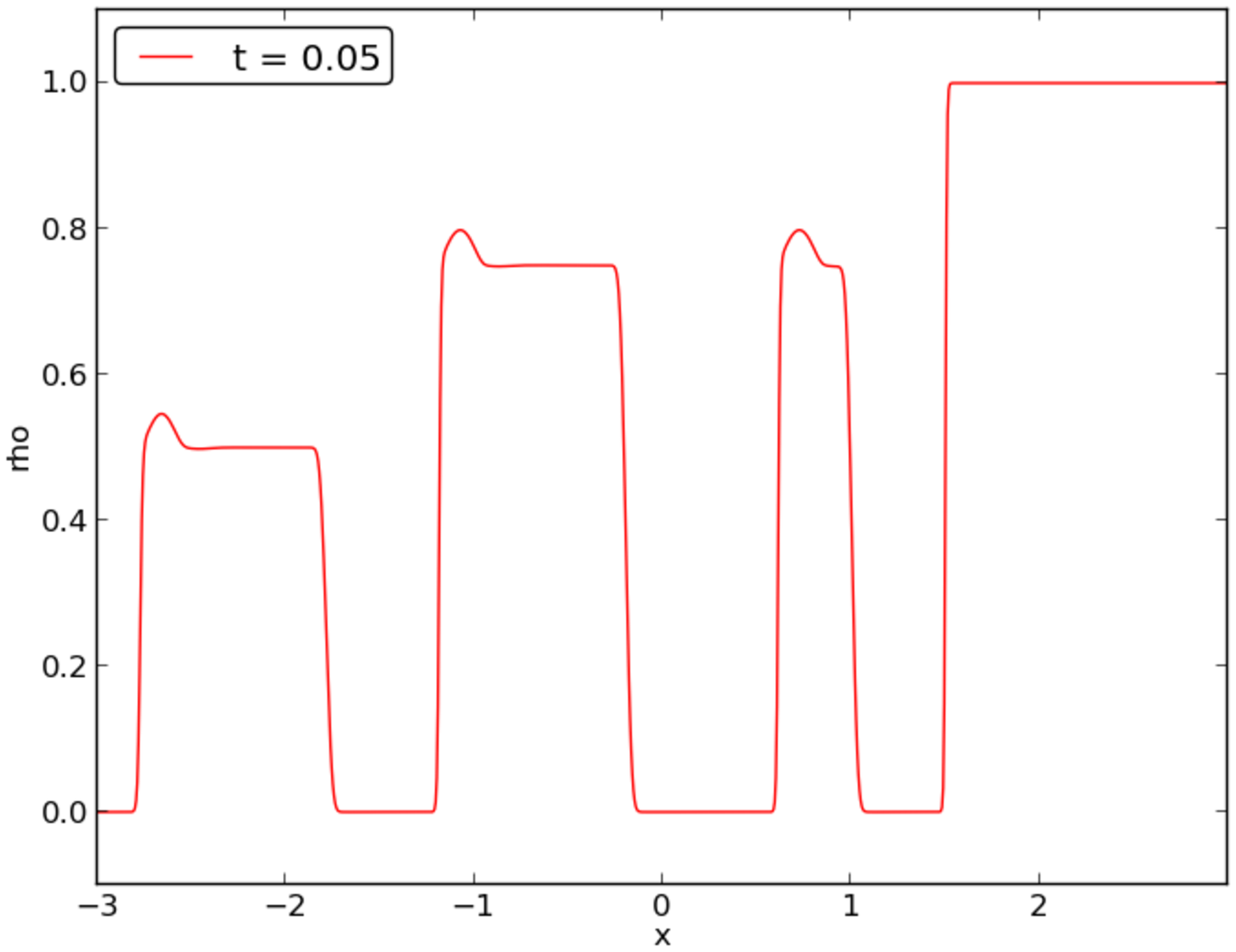}%
  \includegraphics[width=0.2\textwidth, trim=20 0 20
  0]{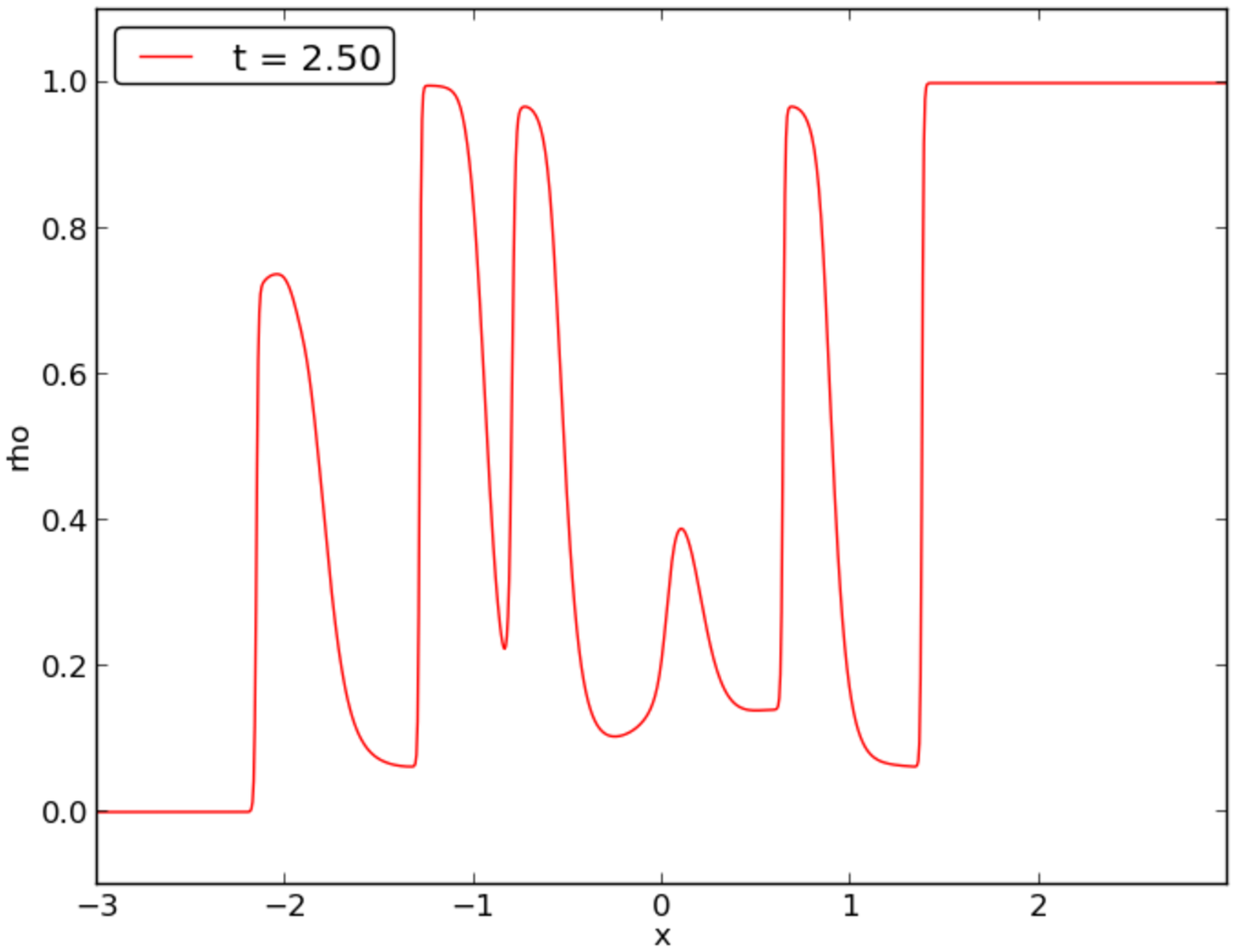}%
  \includegraphics[width=0.2\textwidth, trim=20 0 20
  0]{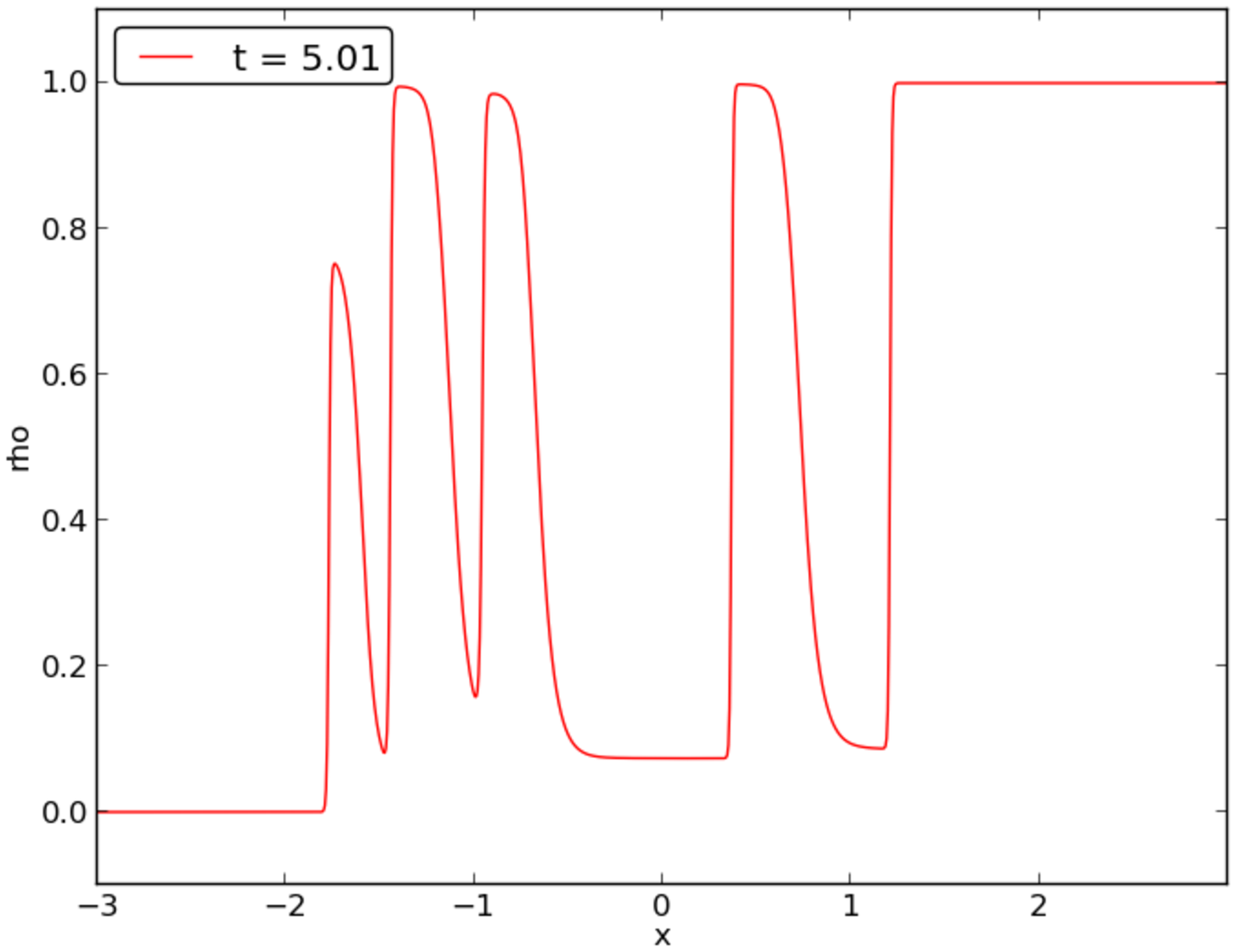}%
  \includegraphics[width=0.2\textwidth, trim=20 0 20
  0]{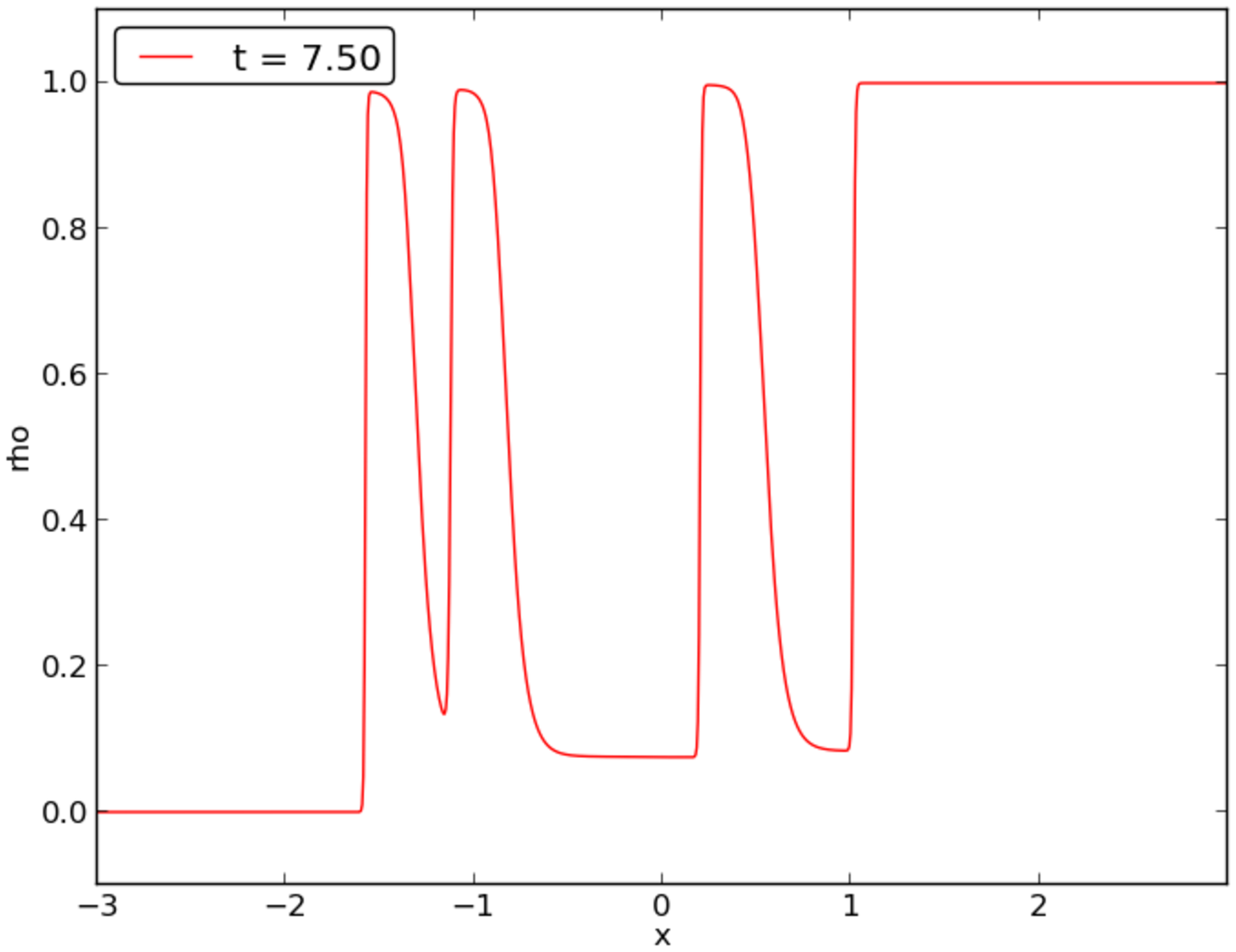}% %%%%%%%%%%%%600!!!!!!!!!!!!!!!!!!!!!!!!!!!!
  \includegraphics[width=0.2\textwidth, trim=20 0 20 0]{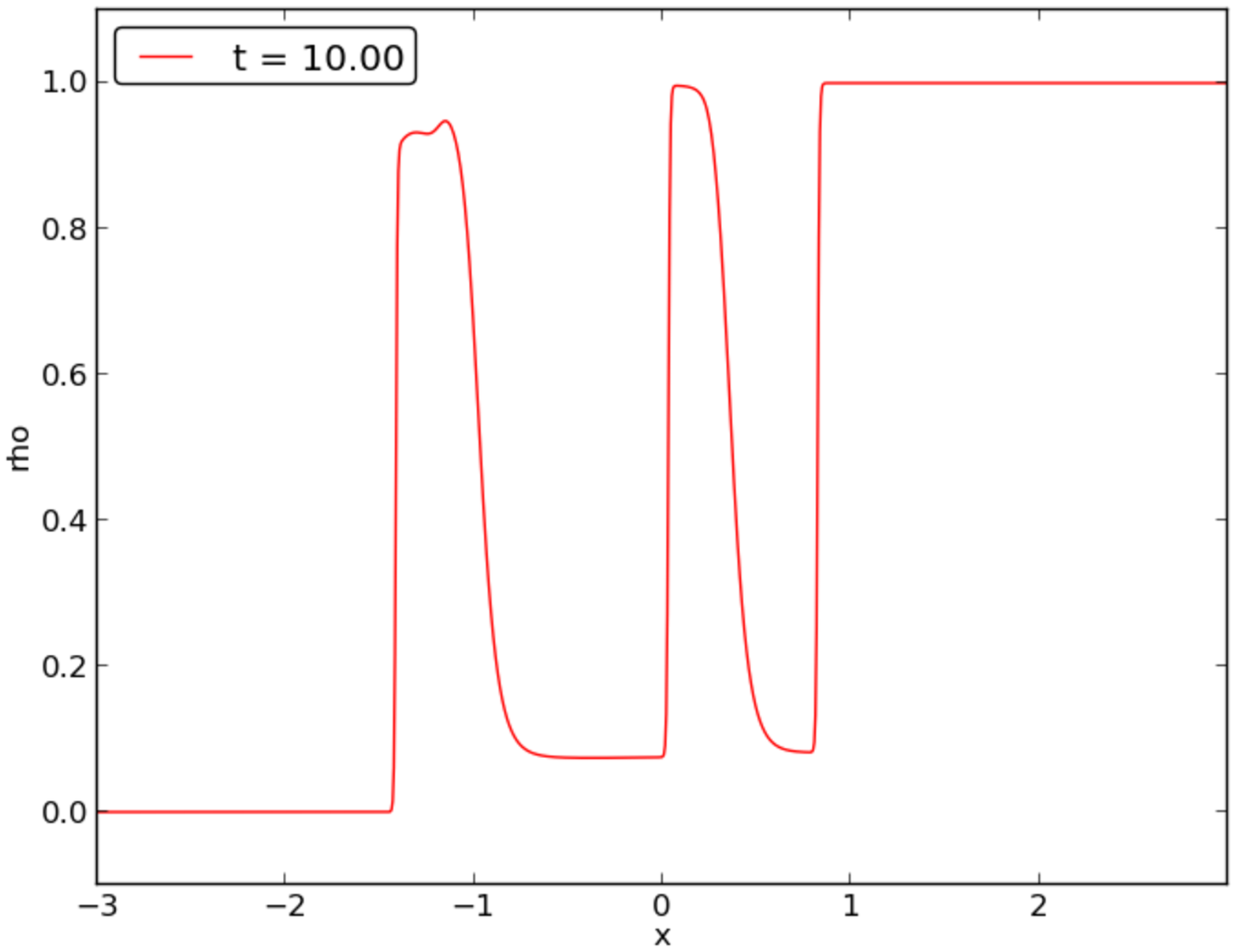}\\
  \includegraphics[width=0.2\textwidth, trim=20 0 20
  0]{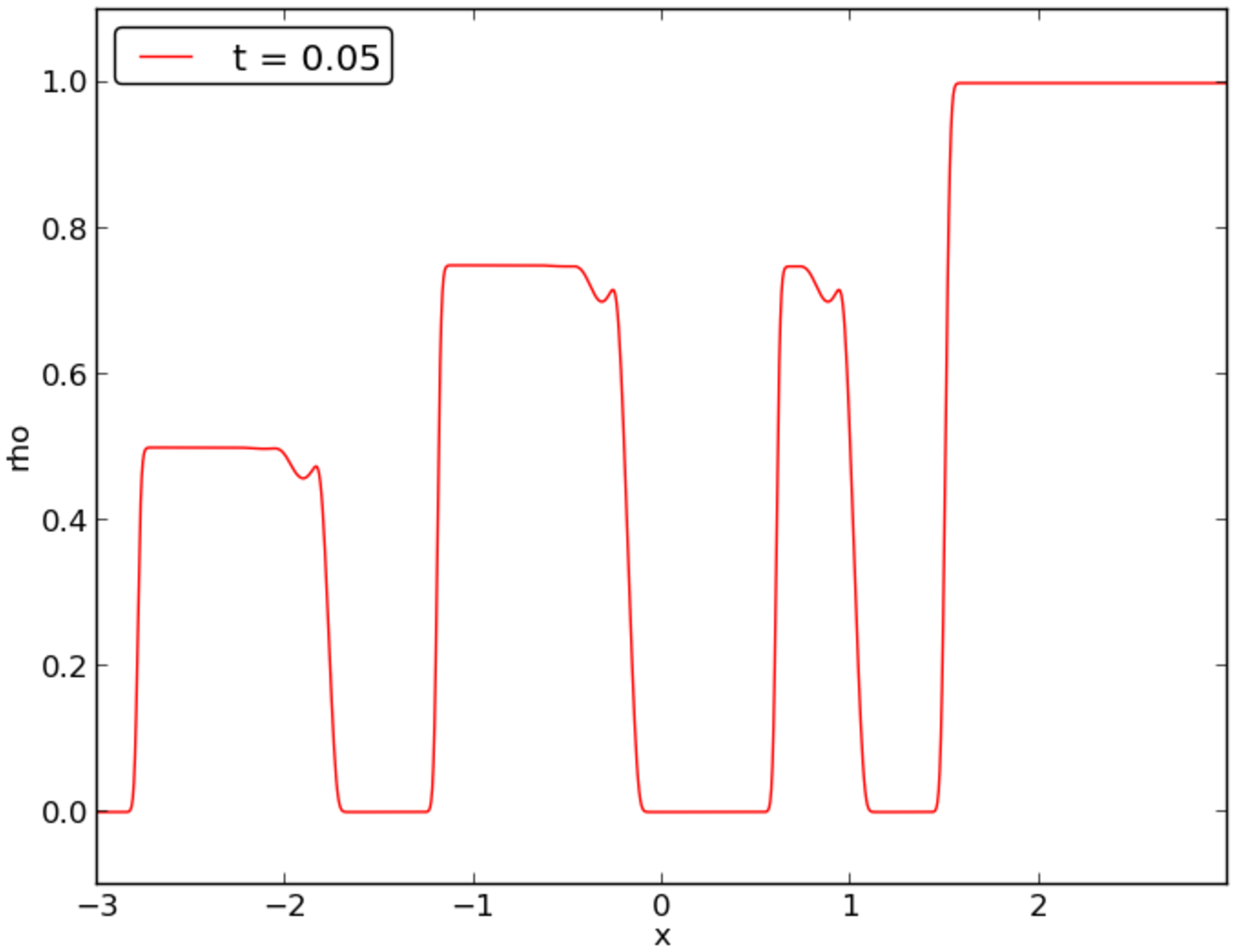}%
  \includegraphics[width=0.2\textwidth, trim=20 0 20
  0]{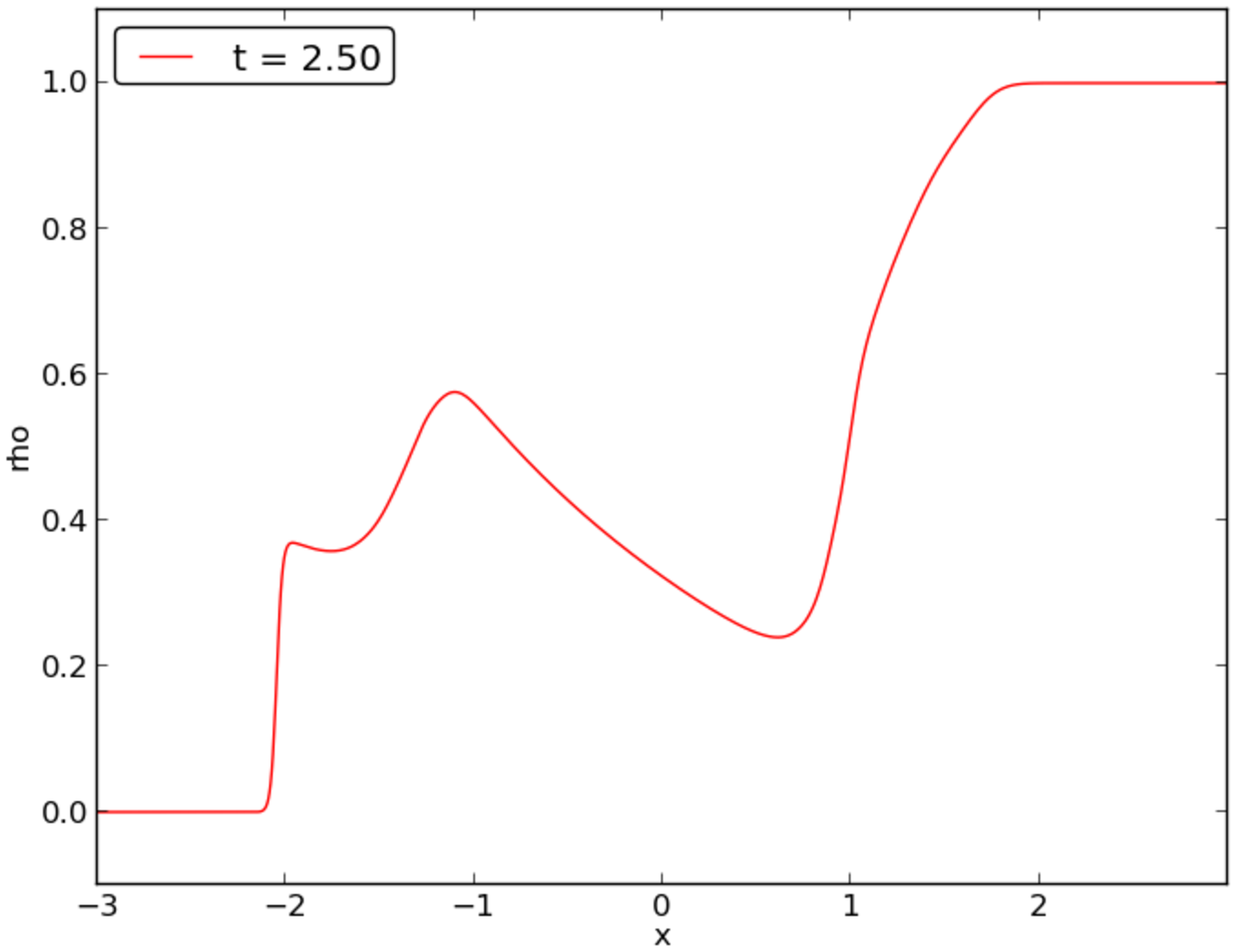}%
  \includegraphics[width=0.2\textwidth, trim=20 0 20
  0]{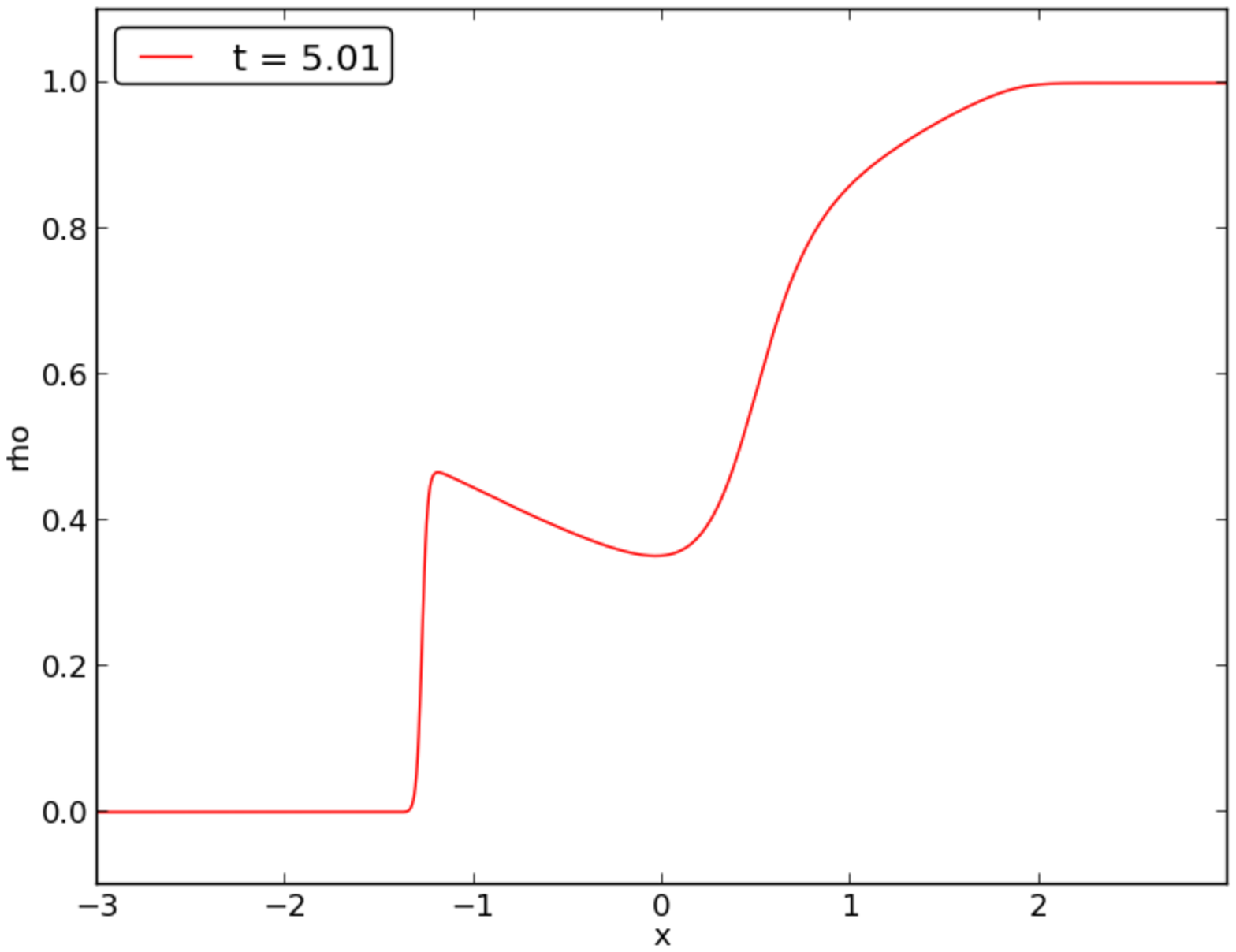}%
  \includegraphics[width=0.2\textwidth, trim=20 0 20
  0]{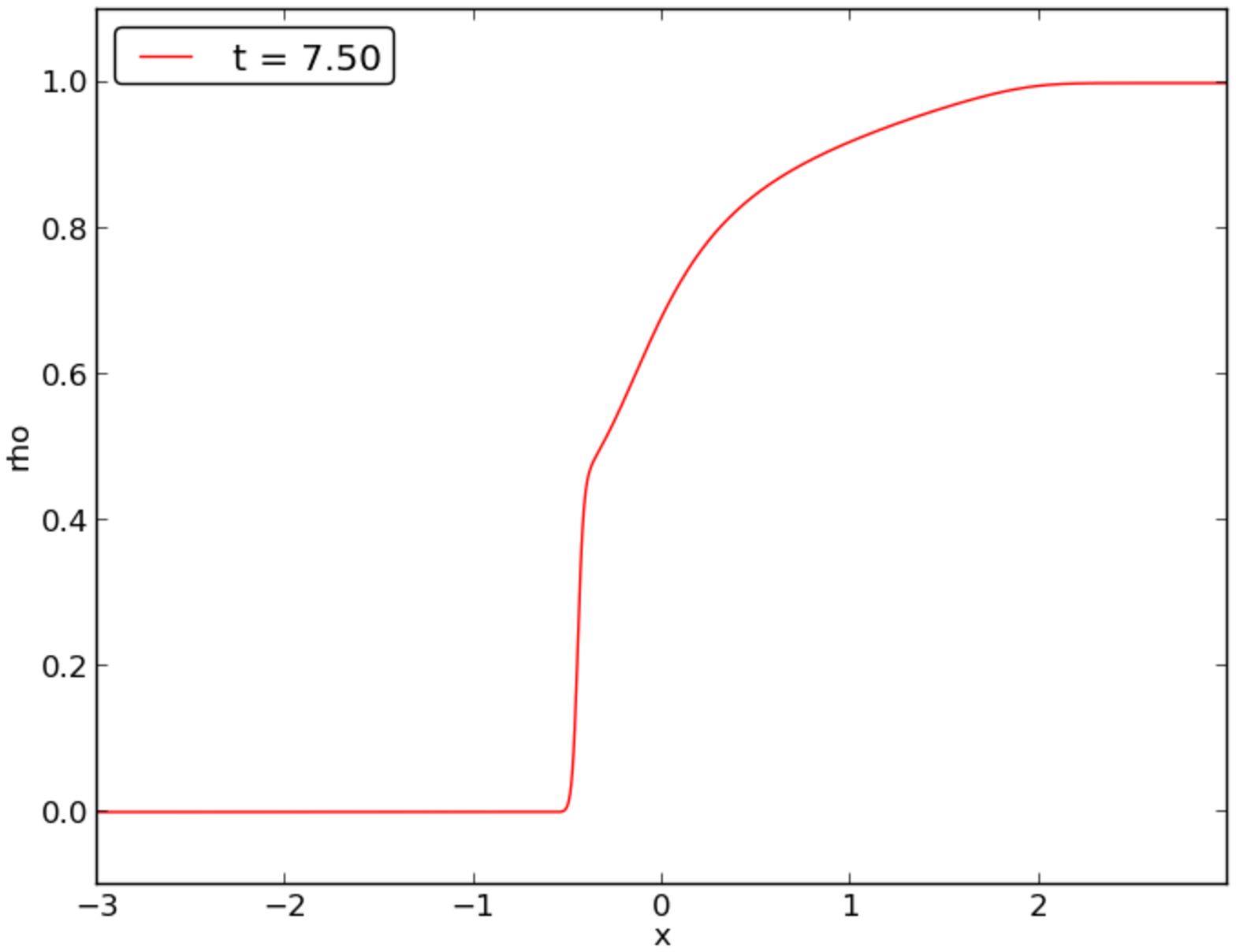}%%%%%%%%%%%%%600!!!!!!!!!!!!!!!!!!!!!!!!!!!!1
  \includegraphics[width=0.2\textwidth, trim=20 0 20 0]{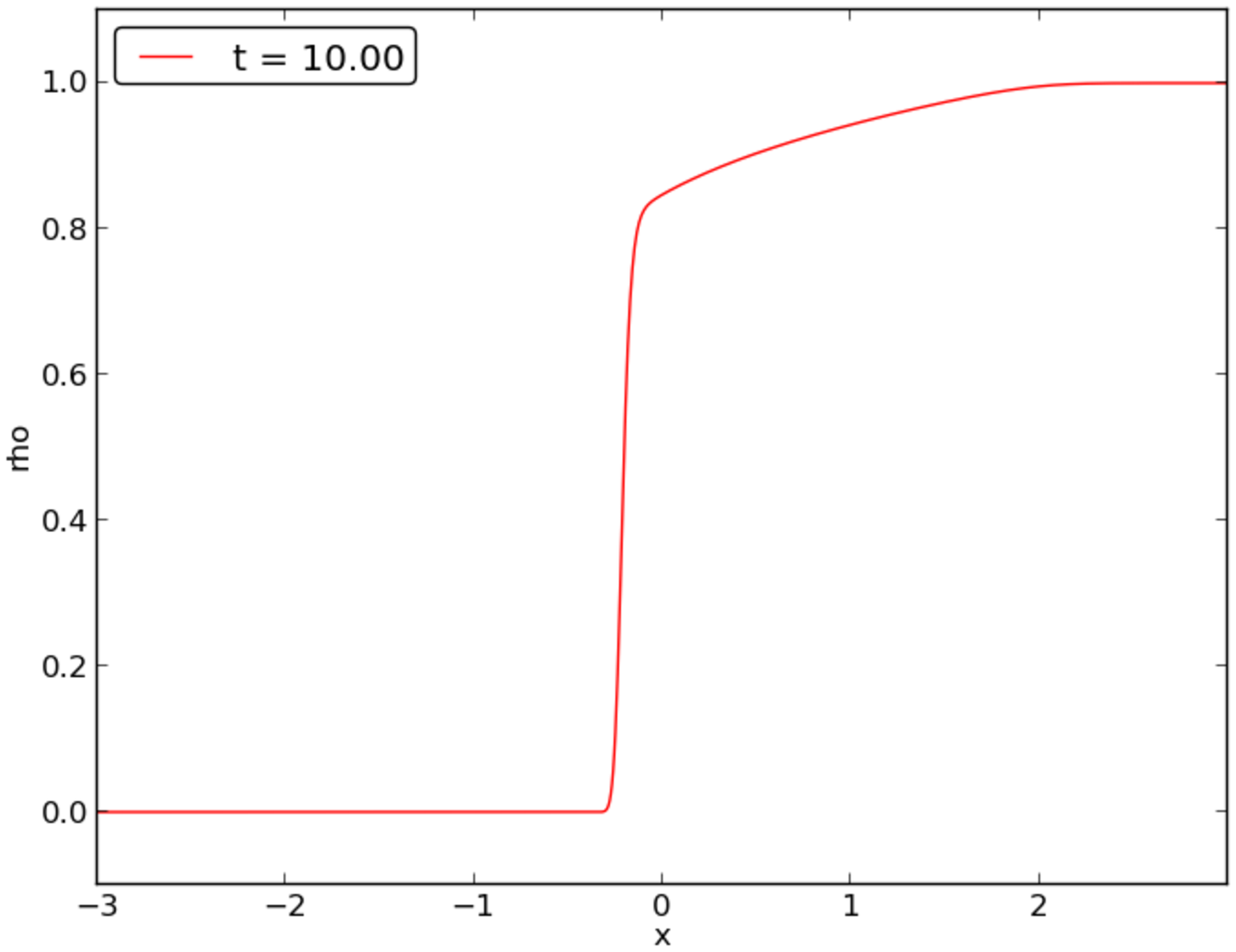}\\
  \caption{Integration
    of~\eqref{eq:1}--\eqref{eq:Traffic}--\eqref{eq:5} in the two
    cases~\eqref{eq:6} at times $t=0.05$, $2.50$, $5.01$, $7.50$,
    $10.00$. Above, drivers look backward while below they look
    forward. Note, already on the first column, the difference in the
    two evolutions, clearly due to the position of the support of
    $\eta$.}
  \label{fig:Traffic}
\end{figure}
However, the qualitative behaviors of the solutions are rather
different in the two situations in~\eqref{eq:6}, see
Figure~\ref{fig:Traffic}. Clearly, the evolution in the case of
drivers looking forward (second line in Figure~\ref{fig:Traffic}) is
far more reasonable, while the backward looking case leads to big
oscillations in the vehicular density.

\subsection{Increase of the Total Variation and of the
  \texorpdfstring{$\L\infty$}{L infinity} Norm}
\label{subs:TV}

This paragraph is devoted to show that Lemma~\ref{lem:pos} and the
total variation bound~(\ref{TV1}) are, at least qualitatively,
optimal. Moreover, the example below shows that the nonlocal
equation~\eqref{eq:1} does not enjoy two standard properties typical
of 1D scalar conservation laws, namely the maximum
principle~\cite[(iv) in Theorem~6.3]{BressanLectureNotes}, see
also~\cite[Theorem~3]{Kruzkov}, and the diminishing of the total
variation~\cite[Theorem~6.1]{BressanLectureNotes}.

\begin{figure}[!htpb]
  \centering
  \includegraphics[width=0.2\textwidth, trim = 30 0 30
  0]{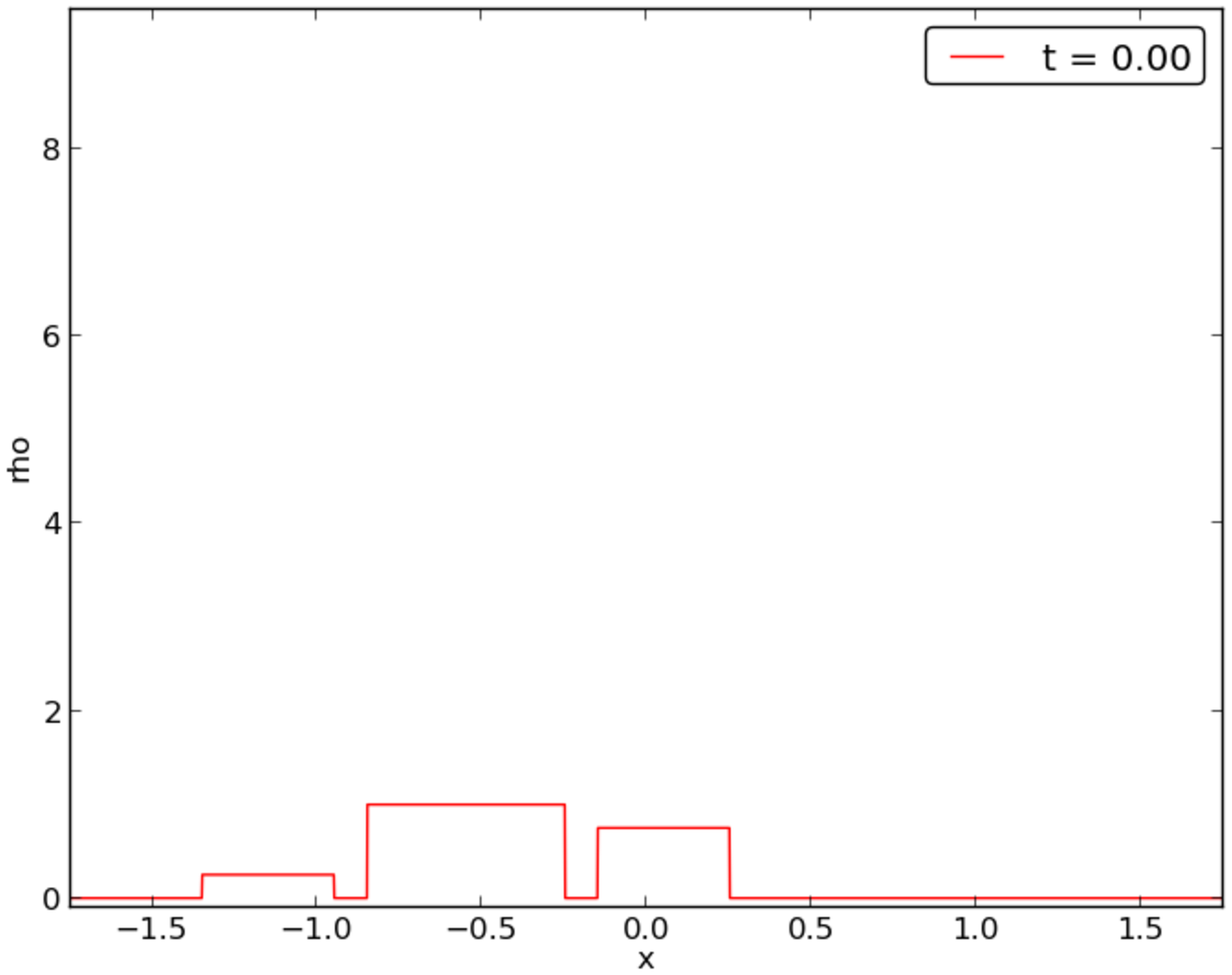}%
  \includegraphics[width=0.2\textwidth, trim = 30 0 30
  0]{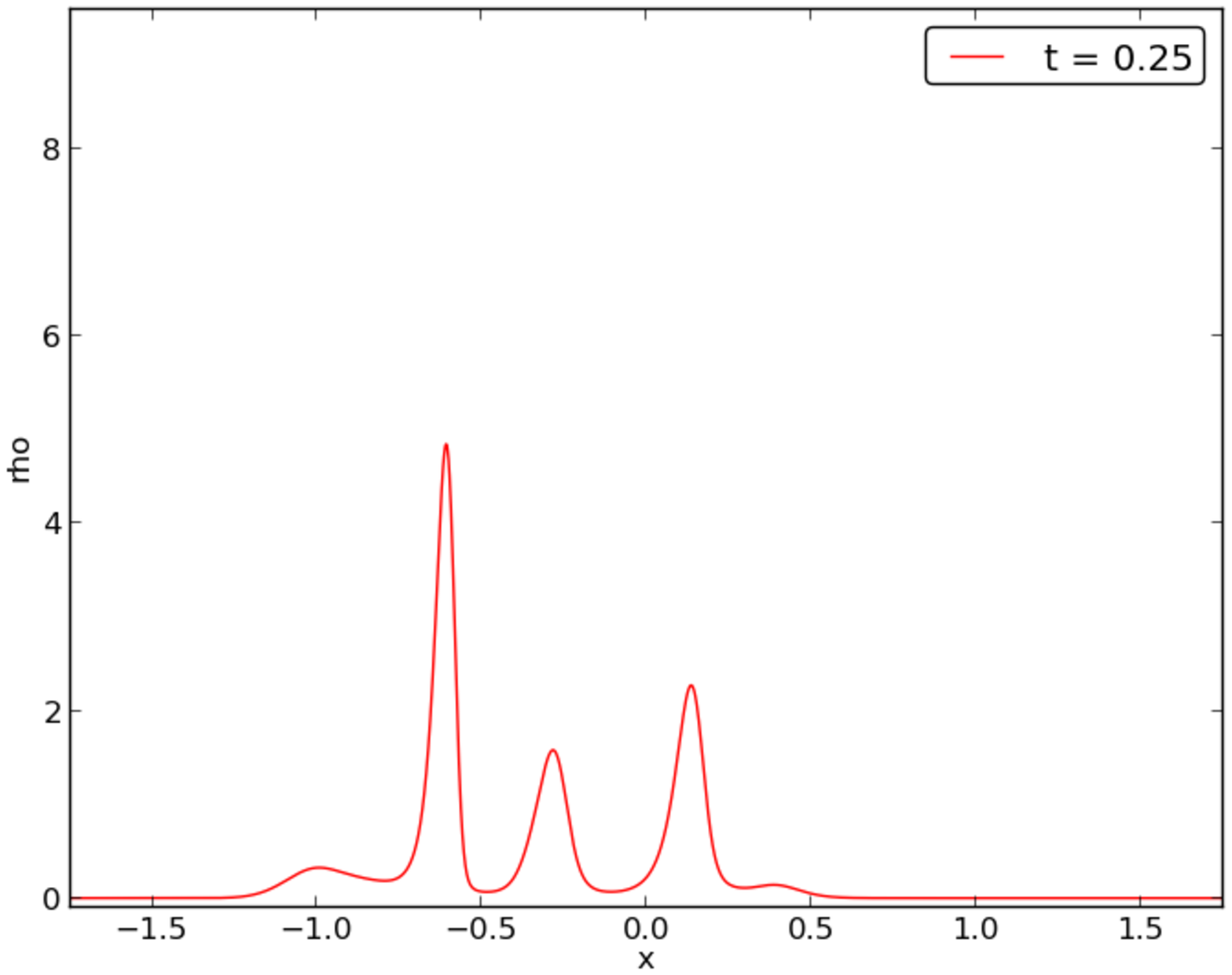}%
  \includegraphics[width=0.2\textwidth, trim = 30 0 30
  0]{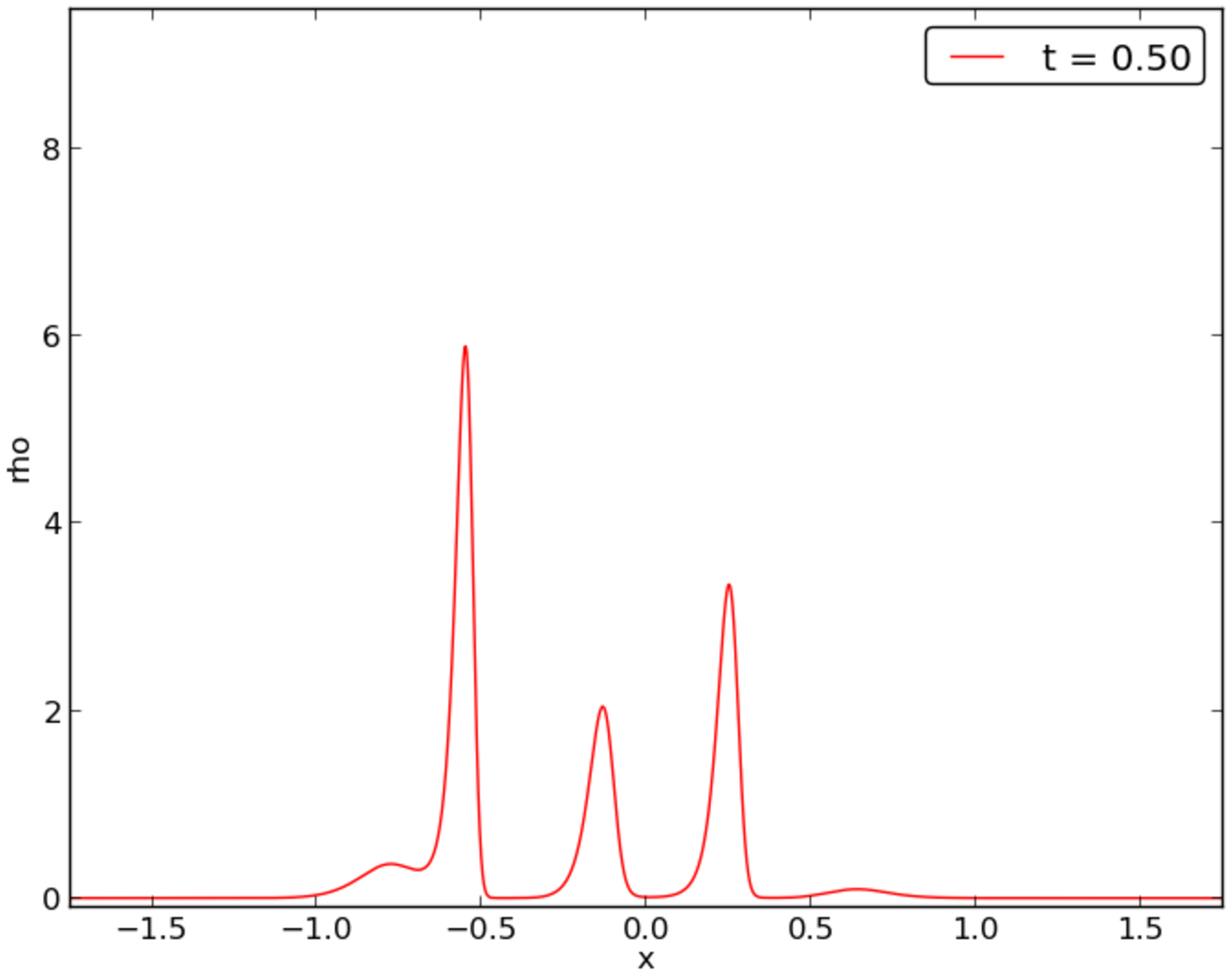}%
  \includegraphics[width=0.2\textwidth, trim = 30 0 30
  0]{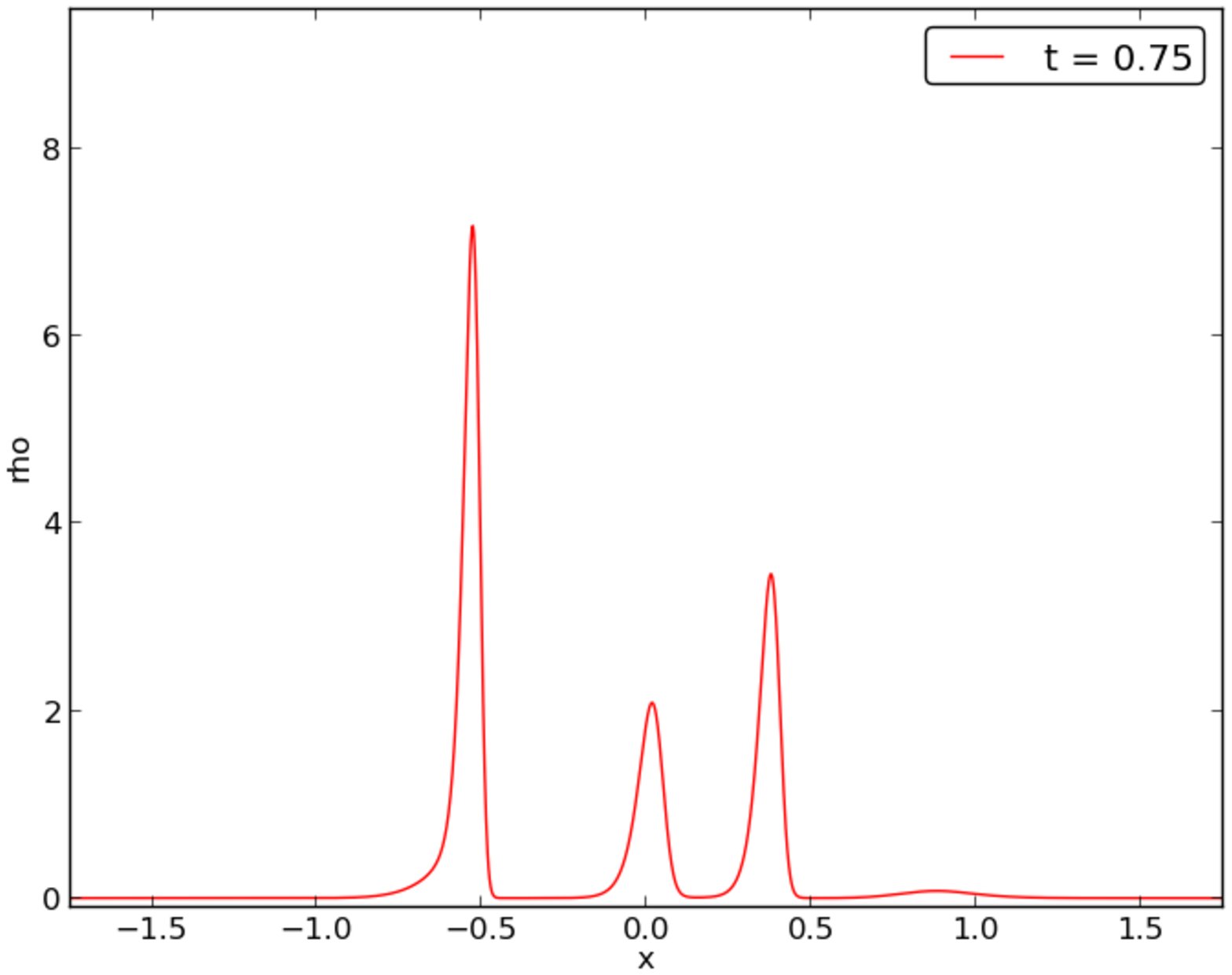}%
  \includegraphics[width=0.2\textwidth, trim = 30 0 30
  0]{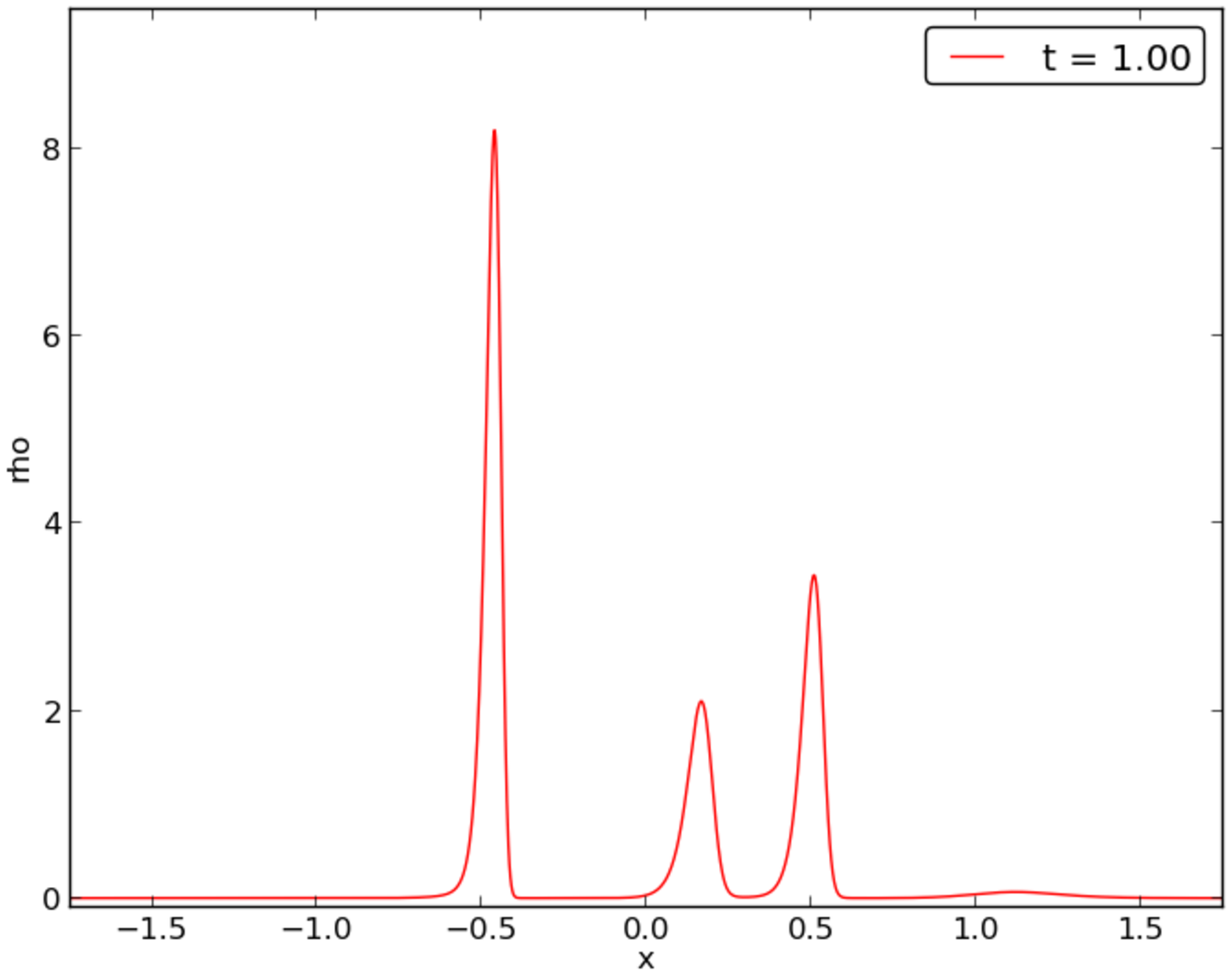}%
  \caption{Numerical integration of~\eqref{eq:1}--\eqref{eq:1bis} with
    $a=0$, $b=0.2$ and with the initial datum~\eqref{eq:idtv} at times
    $t=0.00,\, 0.25,\, 0.50,\, 0.75,\, 1.00$. Note the sharp increase
    in both the $\L\infty$--norm and in the total variation.}
  \label{fig:tv1}
\end{figure}

In Remark~\ref{rem:MP} the assumption that $f (\bar \rho)= 0$ can not
be replaced by $v (\bar \rho) = 0$ to ensure that the solution remains
bounded between $0$ and $\bar \rho$. Let $\bar \rho = 1$ and
consider~(\ref{eq:1}) in the case
\begin{equation}
  \label{eq:1bis}
  f (\rho) = \rho\,, \qquad  v (r) = 1-r
  \quad \mbox{ and } \quad
  \eta (x) = \alpha \left((x-a) (b-x)\right)^{5/2} \, \caratt{[a,b]} (x)
\end{equation}
with $\alpha$ chosen so that $\int_{\reali} \eta = 1$. Then, clearly,
the initial data $\bar\rho (x) \equiv 1$ and $\bar\rho (x) \equiv 0$
are stationary solutions to~(\ref{eq:1})--(\ref{eq:1bis}). However, as
the numerical integration below shows, the initial datum
\begin{equation}
  \label{eq:idtv}
  \rho^o (x) =
  0.25 \, \caratt{[-1.35,\,-0,95]} (x)
  +
  \caratt{[-0.85, \, -0.25]} (x)
  +
  0.75 \, \caratt{[-0.15, \, .25]} (x)
\end{equation}
which satisfies $\rho^o (x) \in [0,1]$ for all $x \in \reali$, yields
a solution $\rho = \rho (t,x)$ that exceeds $\bar \rho = 1$, showing
that~(\ref{eq:1})--(\ref{eq:1bis}) does not satisfy the Maximum
Principle, see Figure~\ref{fig:tv1}.

Remark that the choice~(\ref{eq:1bis}) leads to a flow in~(\ref{eq:1})
which is independent both from $t$ and $x$. In the standard case of
local scalar conservation laws,
\cite[Theorem~6.1]{BressanLectureNotes} ensures that the total
variation of the solution may not increase in time. The numerical
integration below shows that the total variation of the solution
to~(\ref{eq:1})--(\ref{eq:1bis}) may well sharply increase in a very
short time, coherently with~(\ref{TV1}).

\begin{figure}[!htpb]
  \centering
  \includegraphics[width=0.32\textwidth]{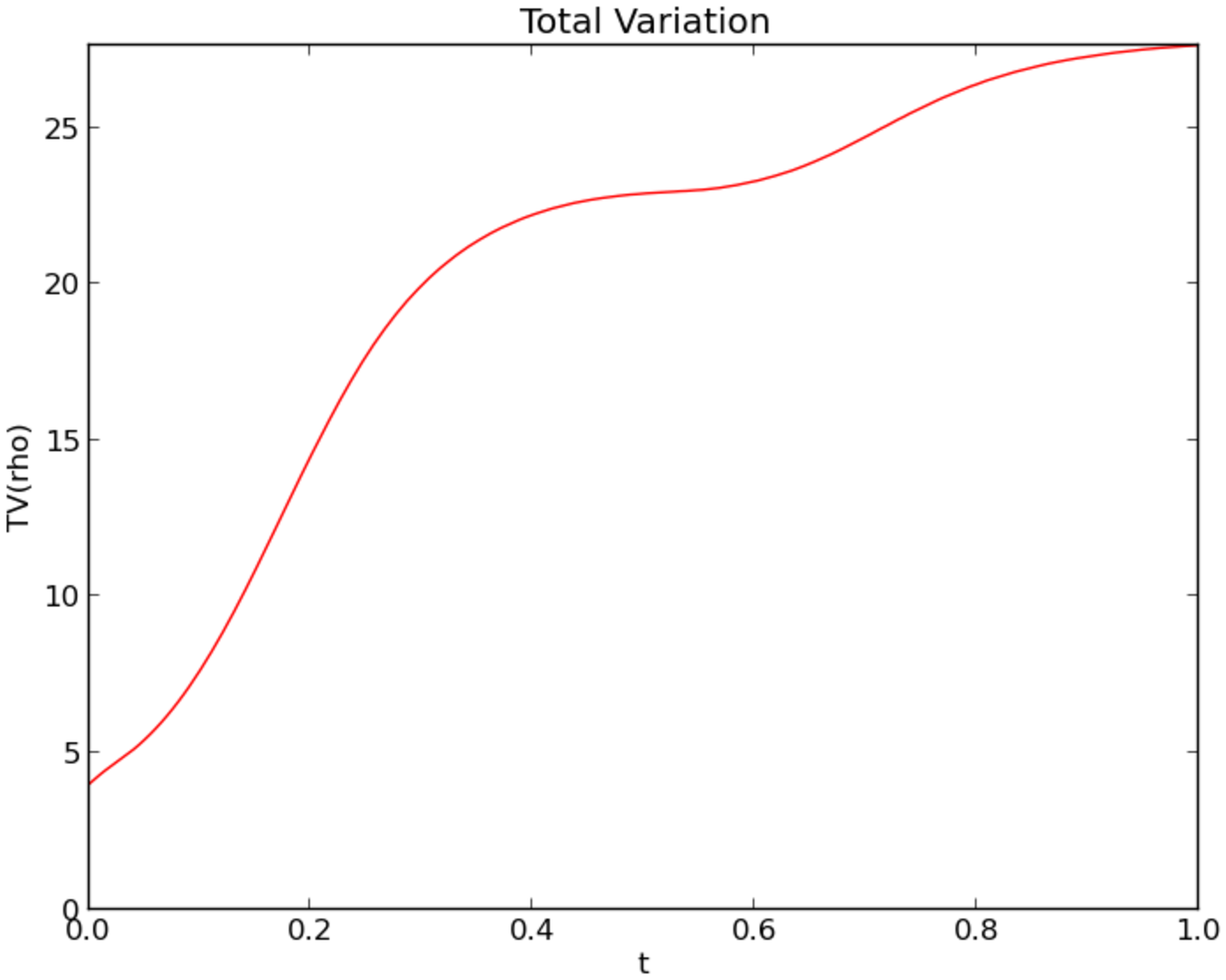}%
  \includegraphics[width=0.32\textwidth]{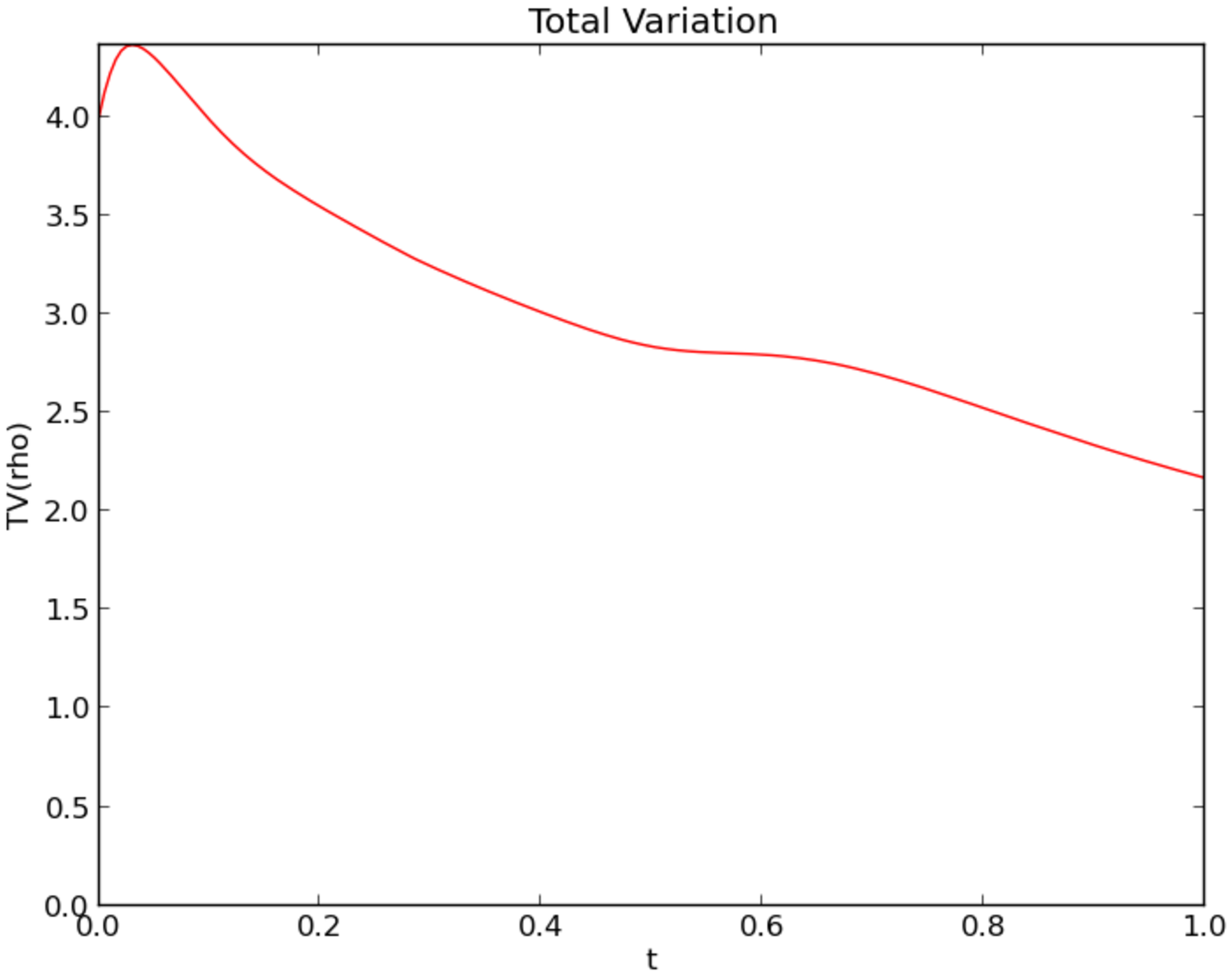}%
  \includegraphics[width=0.32\textwidth]{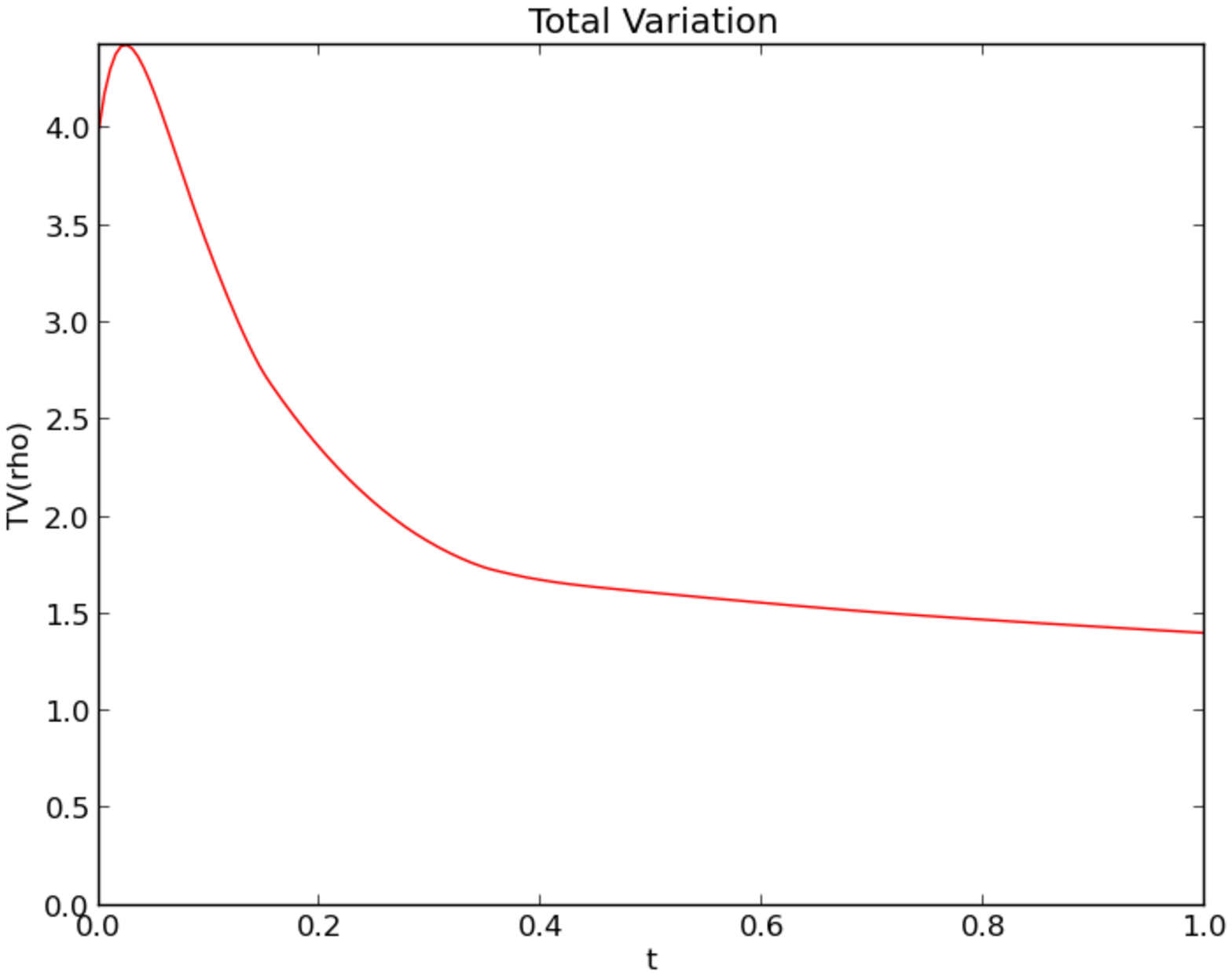}%
  \caption{Total variation of the solution
    to~\eqref{eq:1}--\eqref{eq:1bis}--\eqref{eq:idtv} versus time, in
    the three cases $a=0, \, b=0.2$; $a=-0.1, \, b=0.1$ and $a=-0.2,
    \, b=0$. Remark that the vertical scales in the leftmost diagram
    differs from that used in the middle and on the right. Indeed, the
    initial total variation is the same, $4$, in all cases.}
  \label{fig:tv2}
\end{figure}

It is of interest to note that this behavior depends from the geometry
of the support of $\eta$. Indeed, a translation of the convolution
kernel leads to very different solutions, see Figure~\ref{fig:tv2}.
When the support is contained in $\reali^+$, there is a sharp increase
in the total variation. In the other two cases, when $\spt\eta$ is
centered about the origin or contained in $\reali^-$, there is a small
increase in $\tv (\rho)$ for a small time interval, with a subsequent
decrease.

\subsection{The \emph{Nonlocal to Local} Limit}
\label{subs:Limit}

In this section we use the algorithm~(\ref{LF}) to investigate the
limit in which $\eta$ tends to a Dirac $\delta$, so that the nonlocal
equation~(\ref{eq:1}) tends, at least formally, to the local
conservation law~\eqref{eq:1loc}.  To our knowledge, no analytical
result is at present available on this limit.

\begin{figure}[!htpb]
  \centering
  \includegraphics[width=0.45\textwidth]{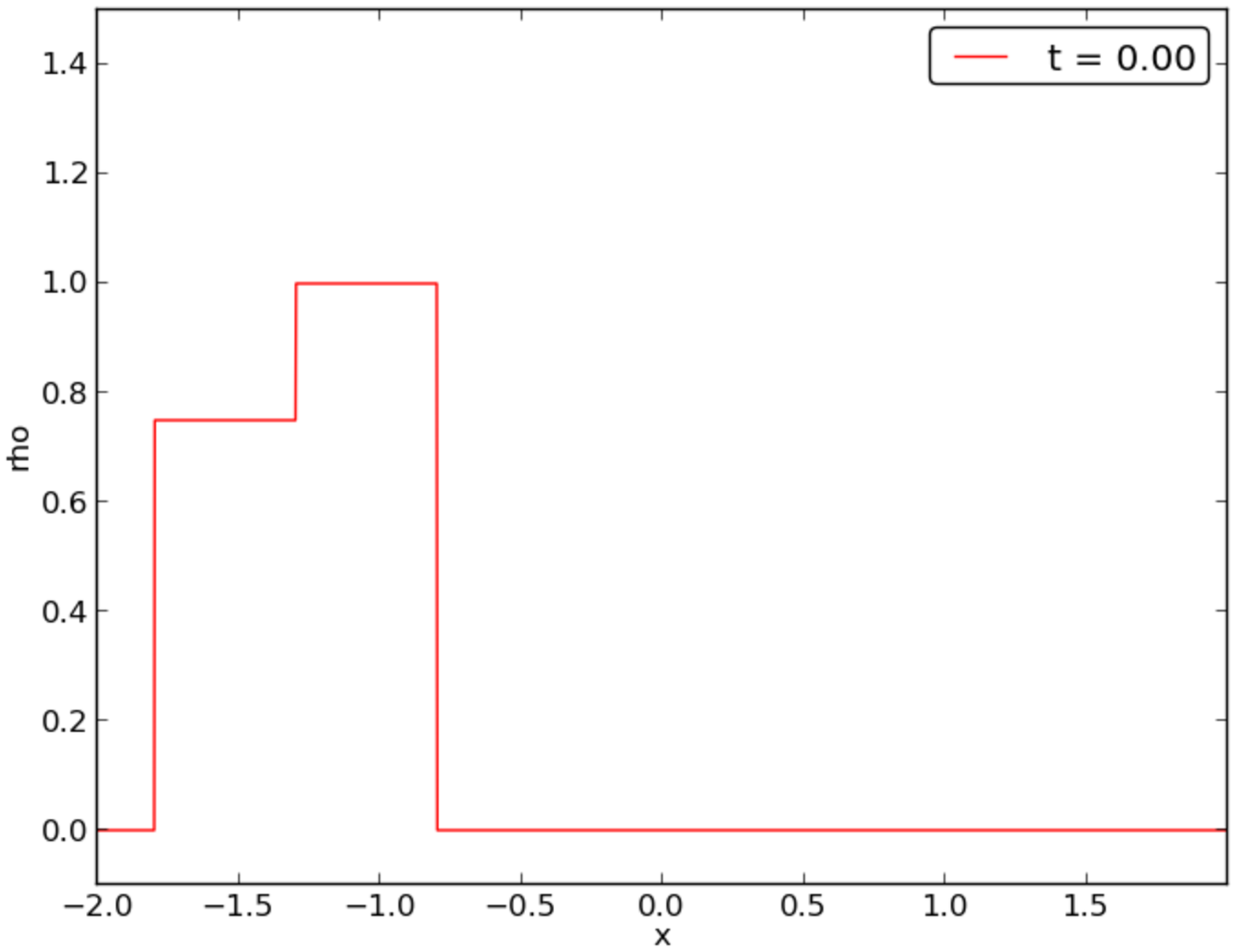}
  \includegraphics[width=0.45\textwidth]{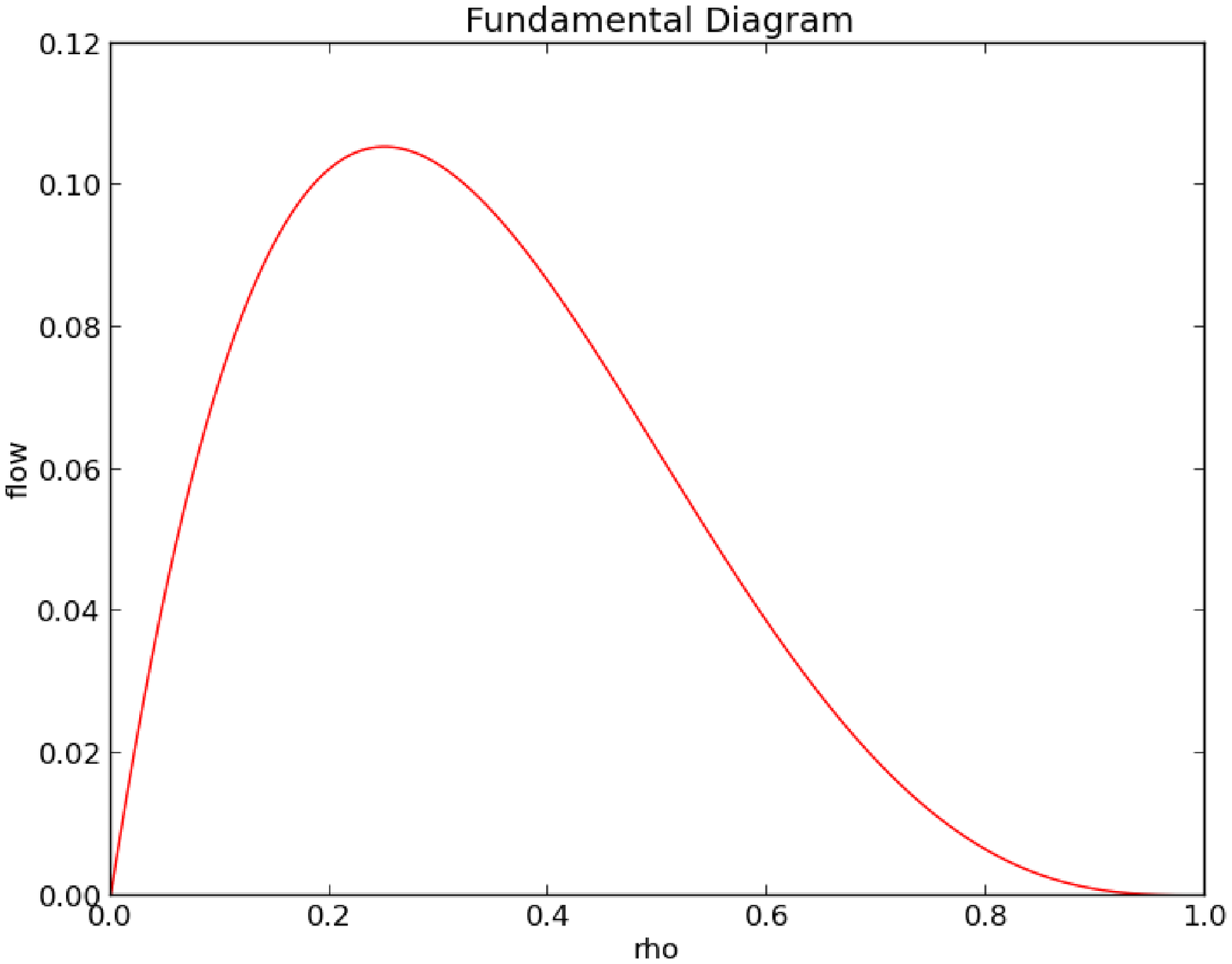}
  \caption{Left, the initial datum in~\eqref{eq:2} and, right, the
    flow~\eqref{eq:2} used in~\eqref{eq:1loc}. Note its change of
    convexity.}
  \label{fig:IDandFlow}
\end{figure}

Consider~(\ref{eq:1}) with the initial datum and the parameters below,
see also Figure~\ref{fig:IDandFlow}, left:
\begin{equation}
  \label{eq:2}
  \begin{array}{@{}r@{\;}c@{\;}l@{}@{\qquad}@{}r@{\;}c@{\;}l@{}}
    f (t,x,\rho)
    & = &
    \rho
    &
    v (r) & = & (1-r)^3 \, \caratt{]-\infty,1[} (r)
    \\[10pt]
    \eta_a (x)
    & = &
    \alpha_a \, (a^2-x^2)^{5/2} \, \caratt{[-a, a]} (x)
    &
    \rho_o (x)
    & = &
    \frac{3}{4} \, \caratt{[-1.8, -1.3[} (x)
    +
    \caratt{[-1.3, \, -0.8]}(x) \,,
  \end{array}
\end{equation}
and with the following choices for $a$:
\begin{displaymath}
  a = 0.25, \quad 0.1, \quad 0.05,
\end{displaymath}%
with $\alpha_a$ computed so that $\int_{\reali} \eta_a (x)
\d{x}=1$. As limit case, we consider the standard conservation
law~\eqref{eq:1loc} with $f$ and $v$ as in~\eqref{eq:2}, see also
Figure~\ref{fig:IDandFlow}, right.
\begin{figure}[!htpb]
  \centering
  \includegraphics[width=0.25\textwidth, trim=20 0 20 0]{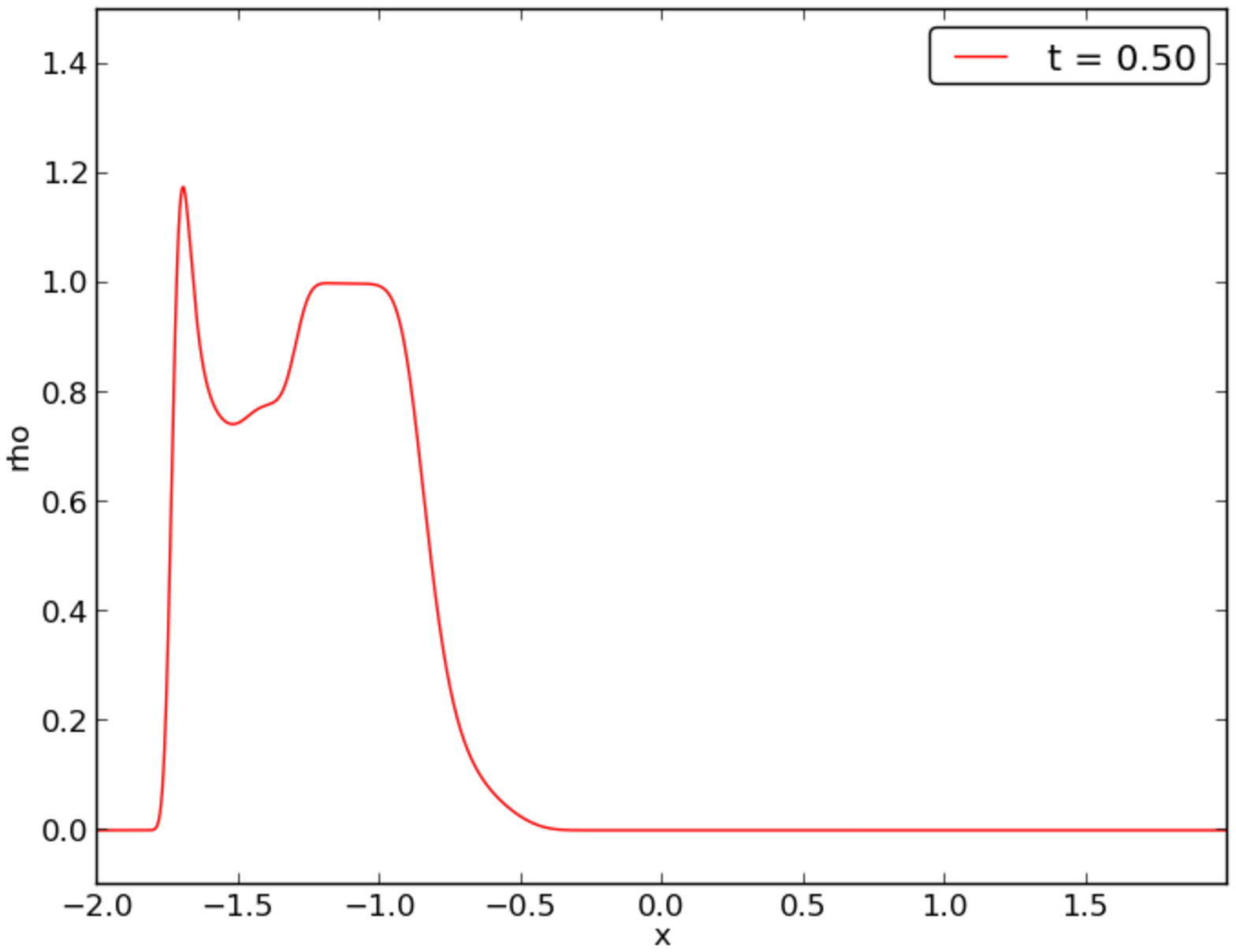}%
  \includegraphics[width=0.25\textwidth, trim=20 0 20 0]{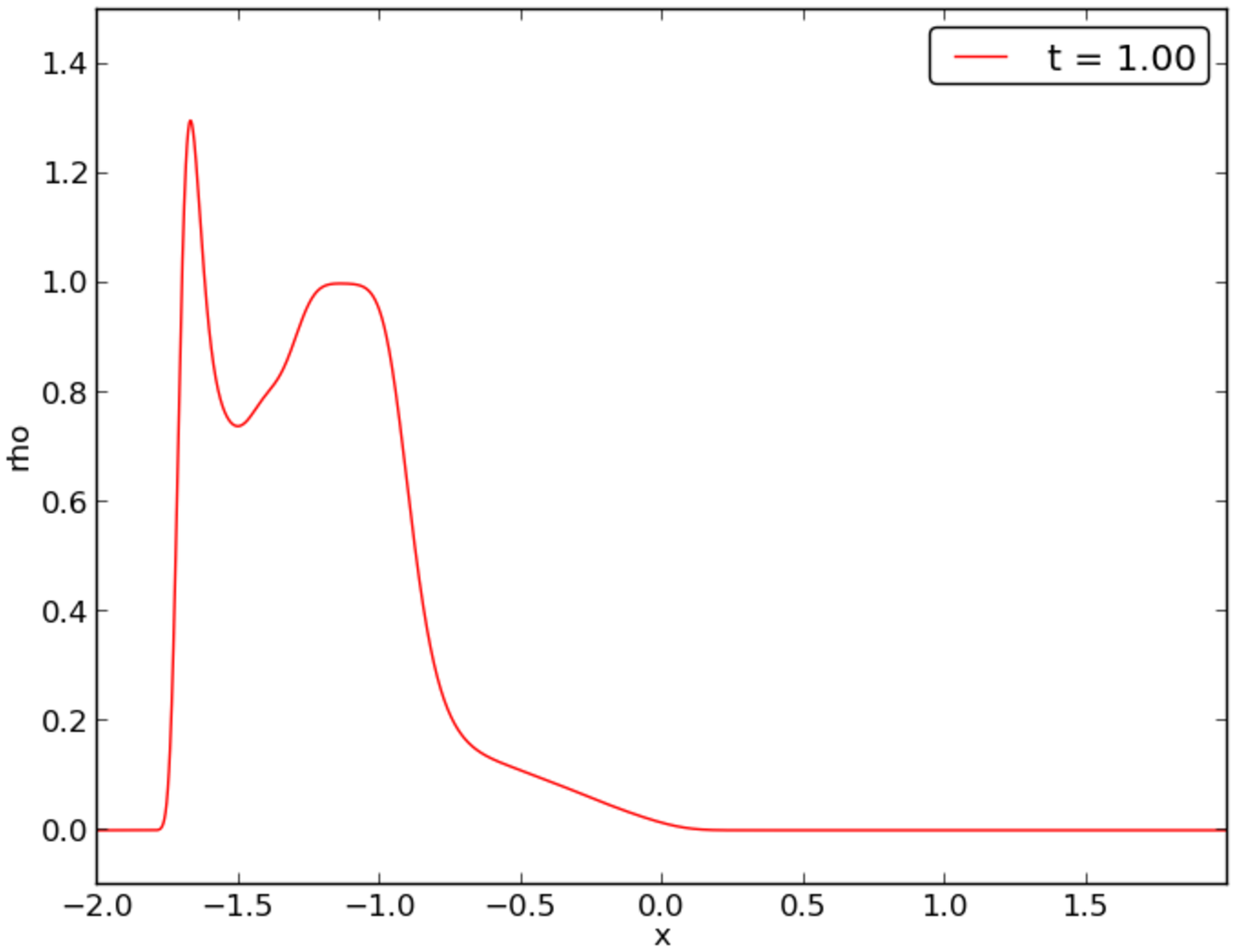}%
  \includegraphics[width=0.25\textwidth, trim=20 0 20 0]{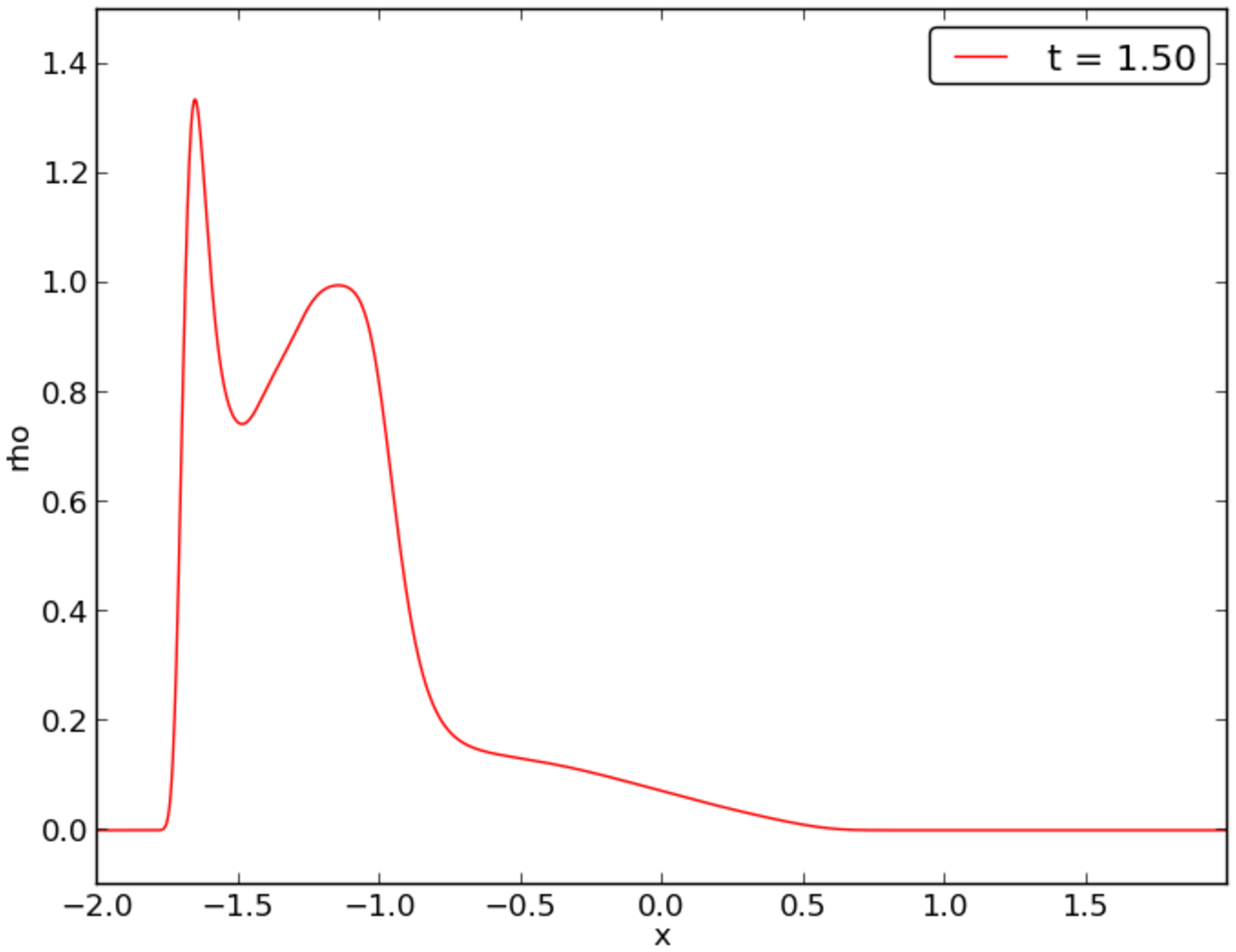}%
  \includegraphics[width=0.25\textwidth, trim=20 0 20 0]{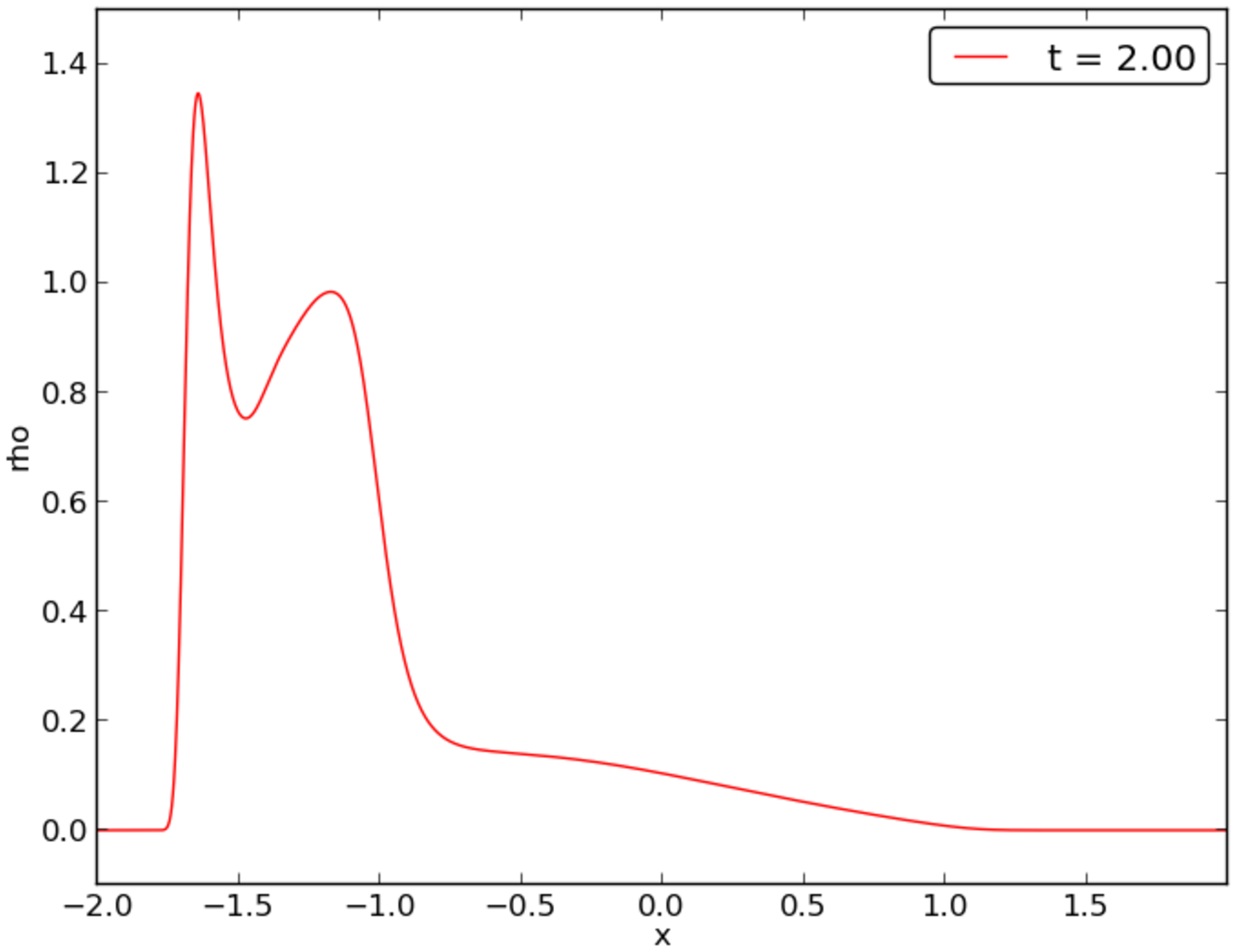}\\
  \includegraphics[width=0.25\textwidth, trim=20 0 20 0]{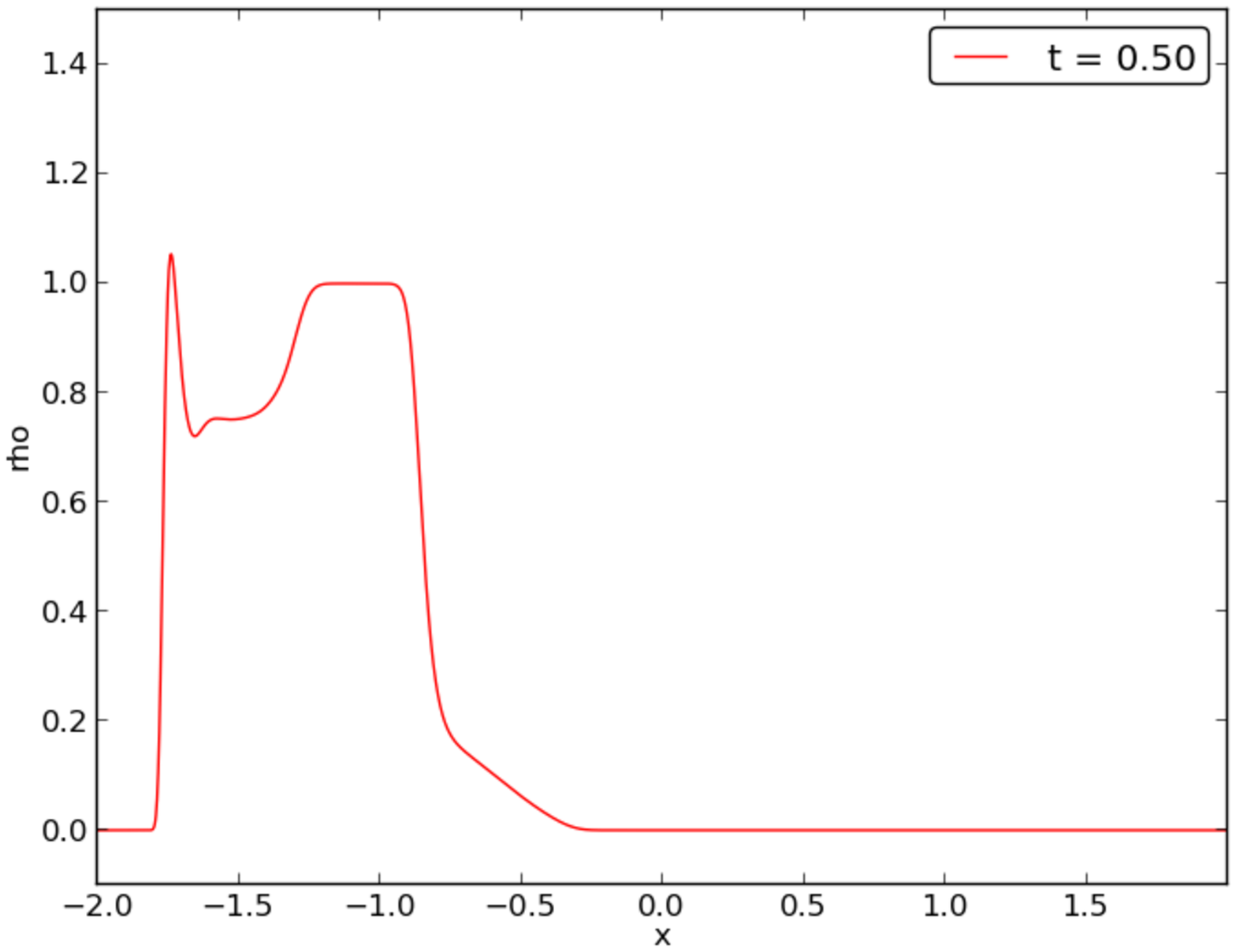}%
  \includegraphics[width=0.25\textwidth, trim=20 0 20 0]{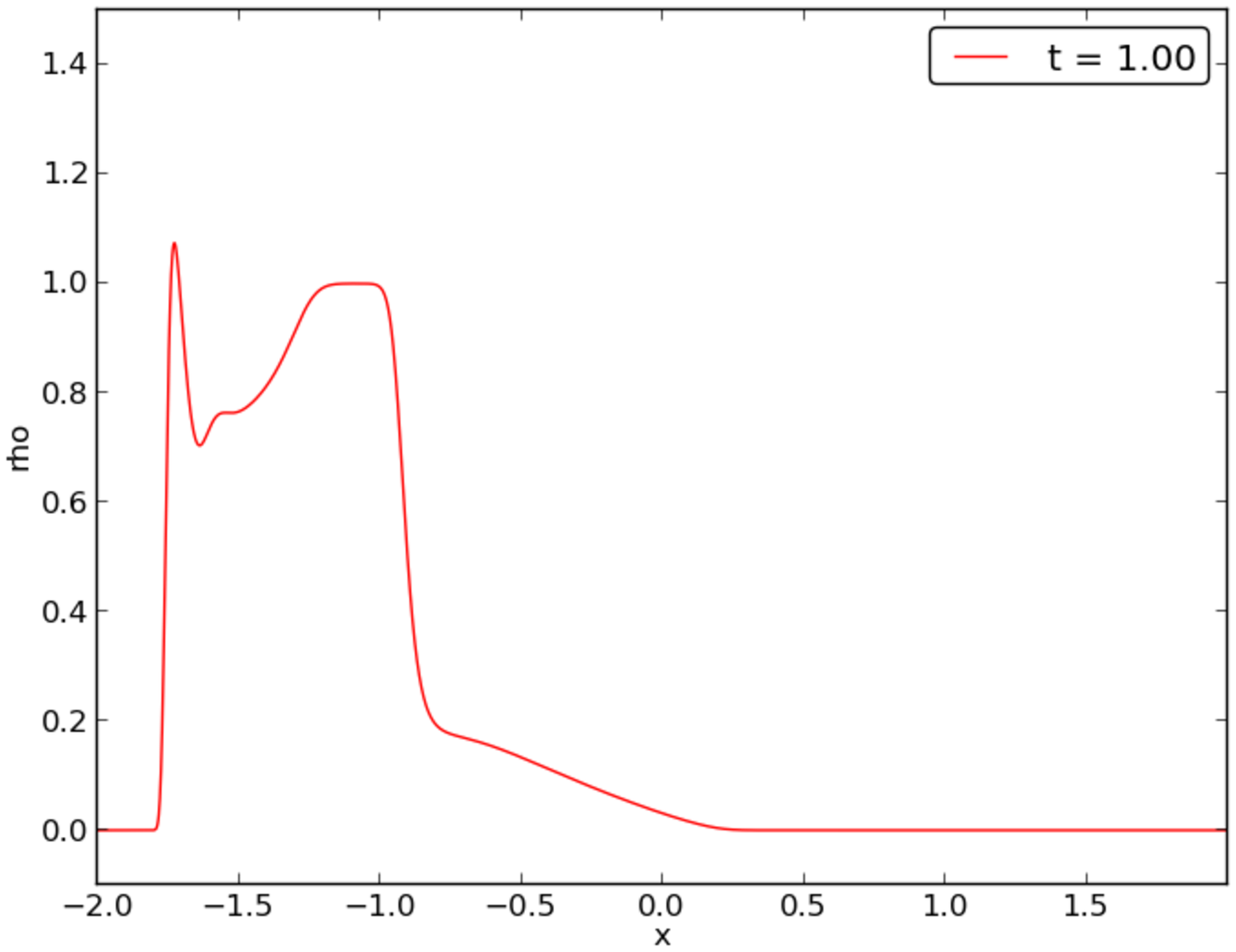}%
  \includegraphics[width=0.25\textwidth, trim=20 0 20 0]{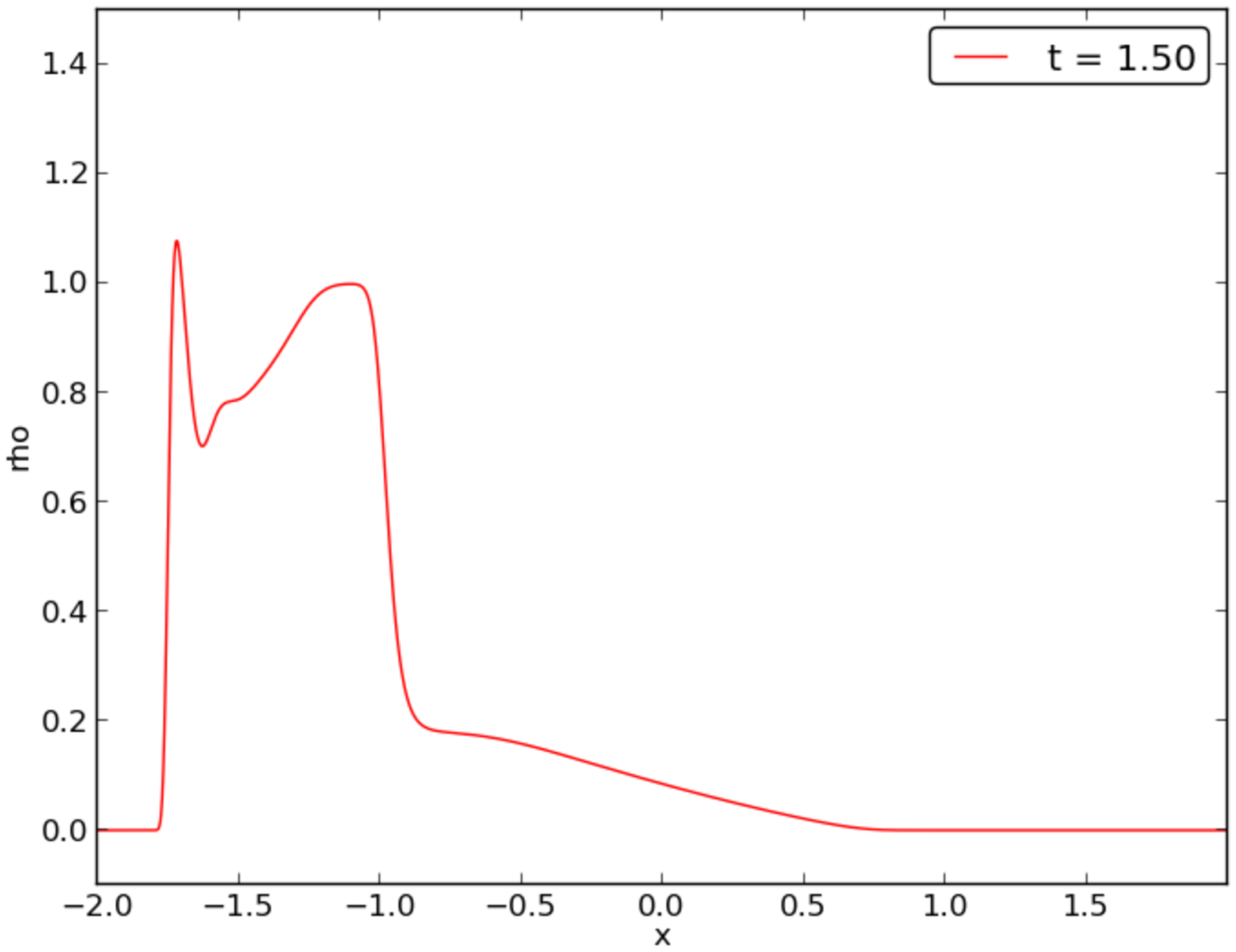}%
  \includegraphics[width=0.25\textwidth, trim=20 0 20 0]{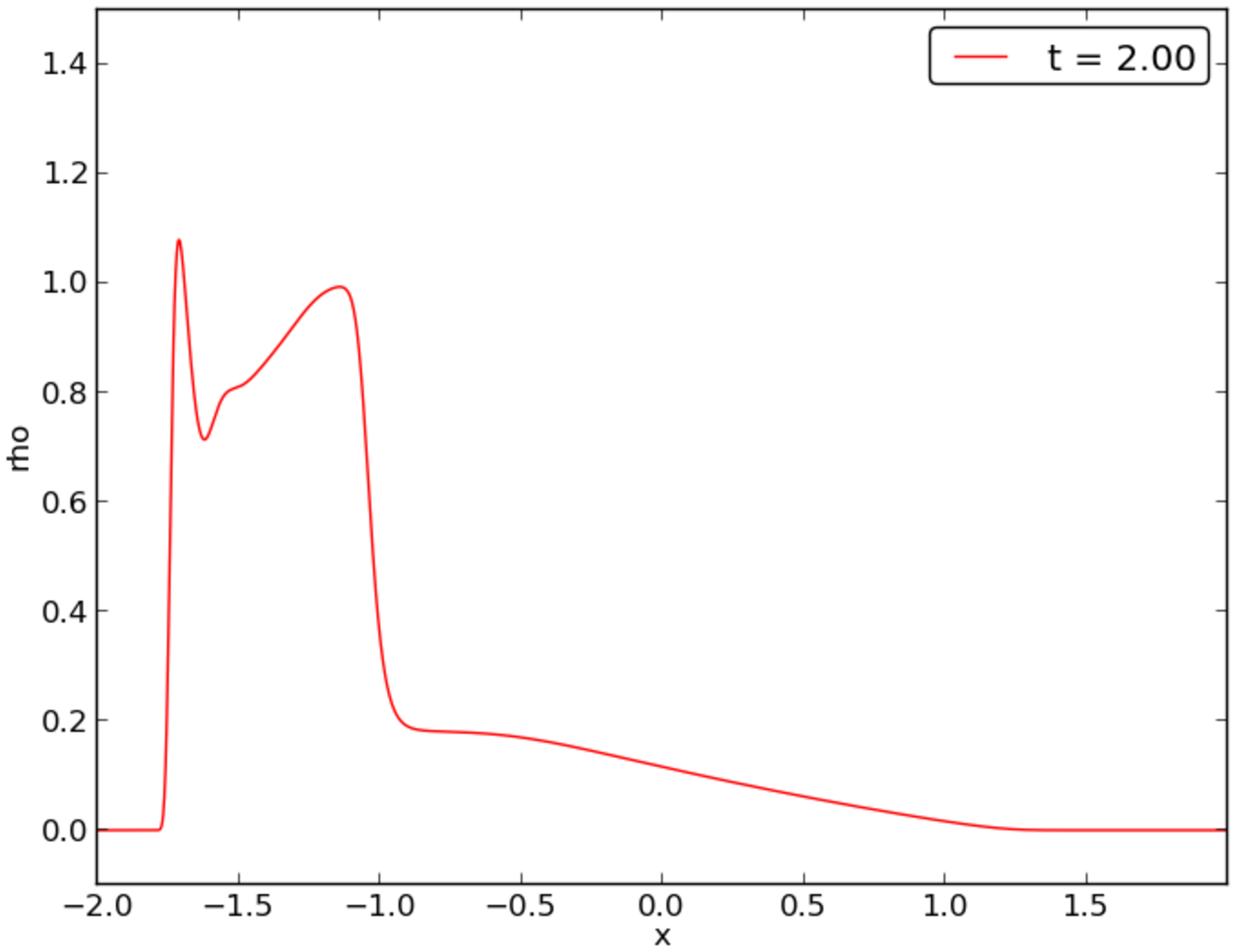}\\
  \includegraphics[width=0.25\textwidth, trim=20 0 20 0]{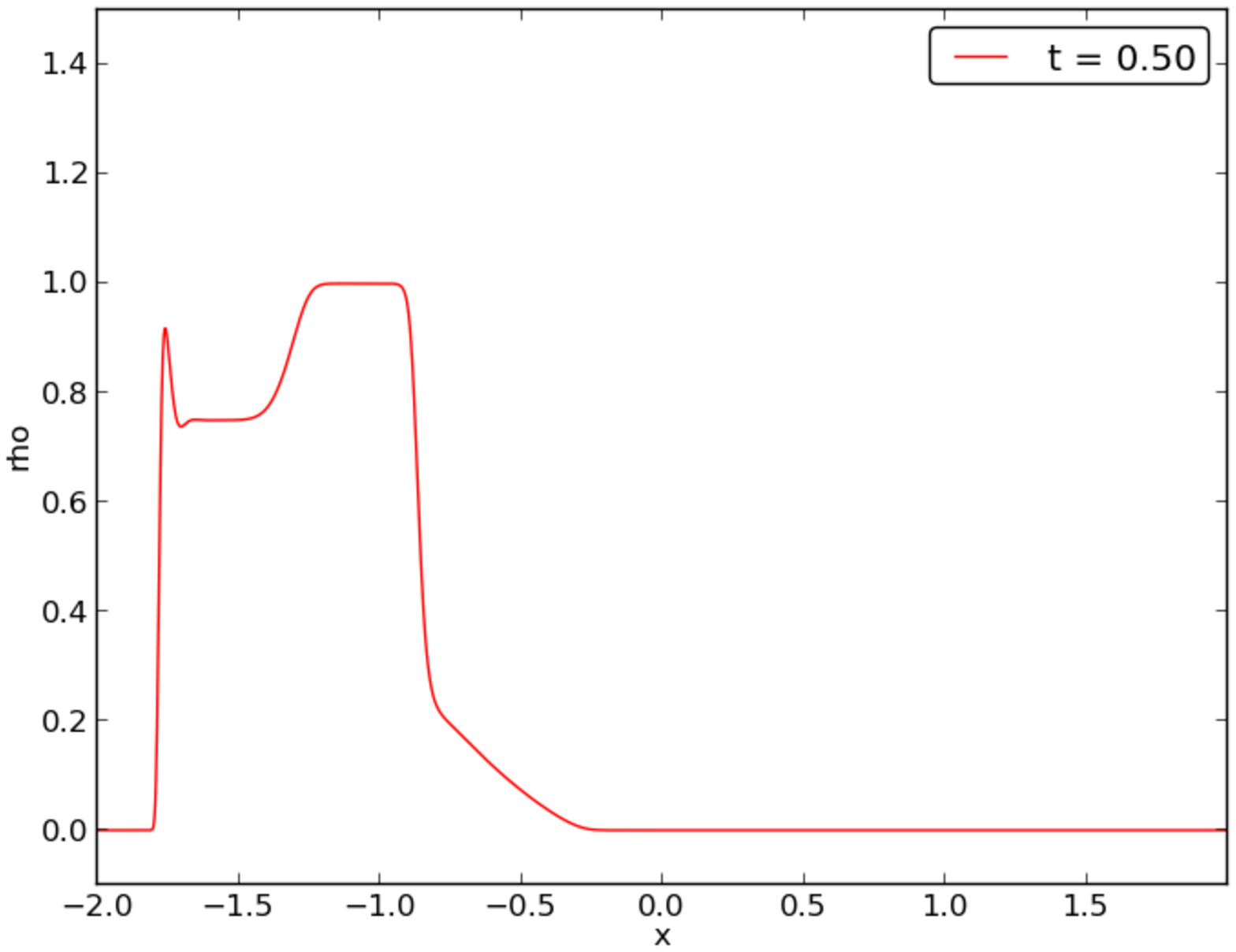}%
  \includegraphics[width=0.25\textwidth, trim=20 0 20 0]{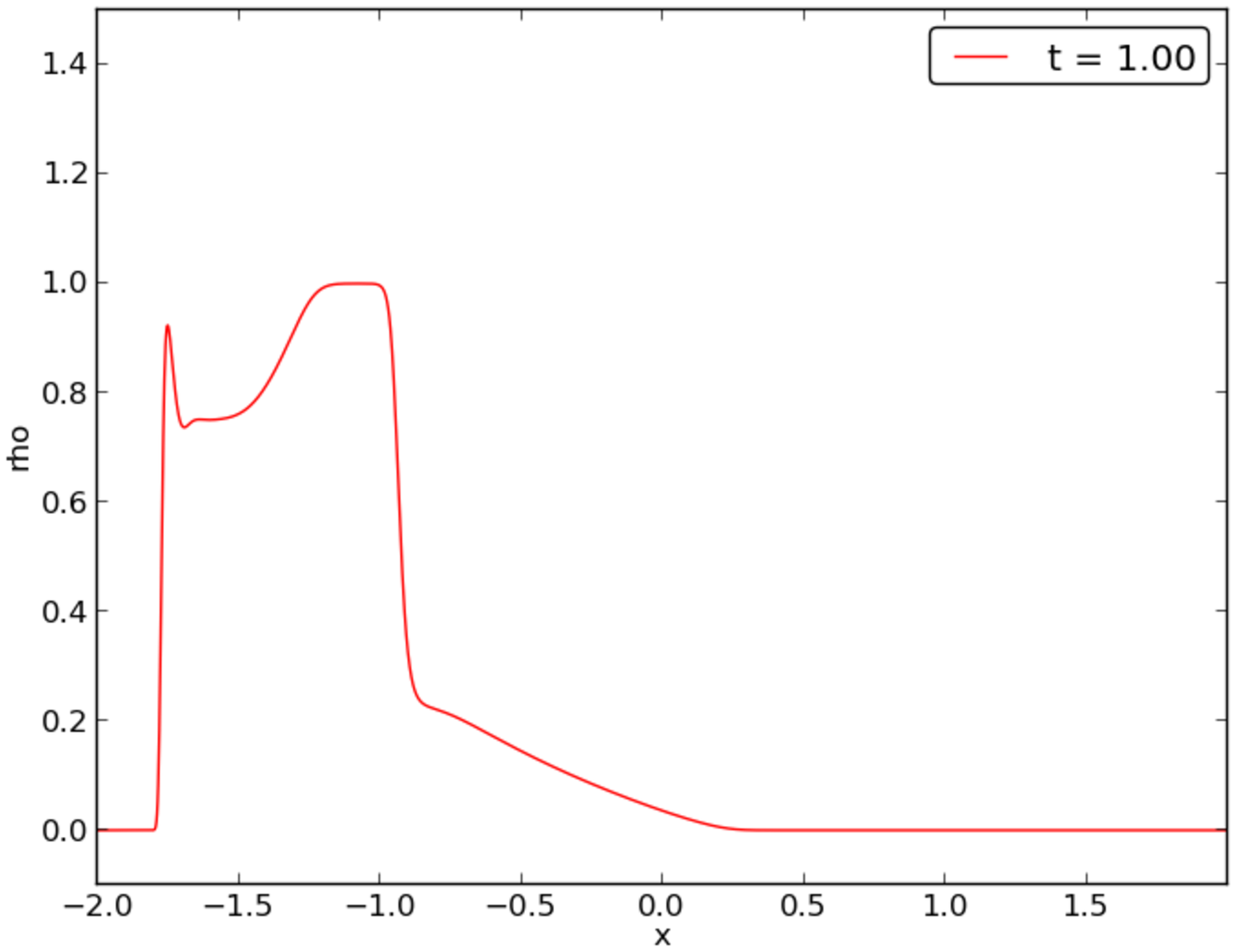}%
  \includegraphics[width=0.25\textwidth, trim=20 0 20 0]{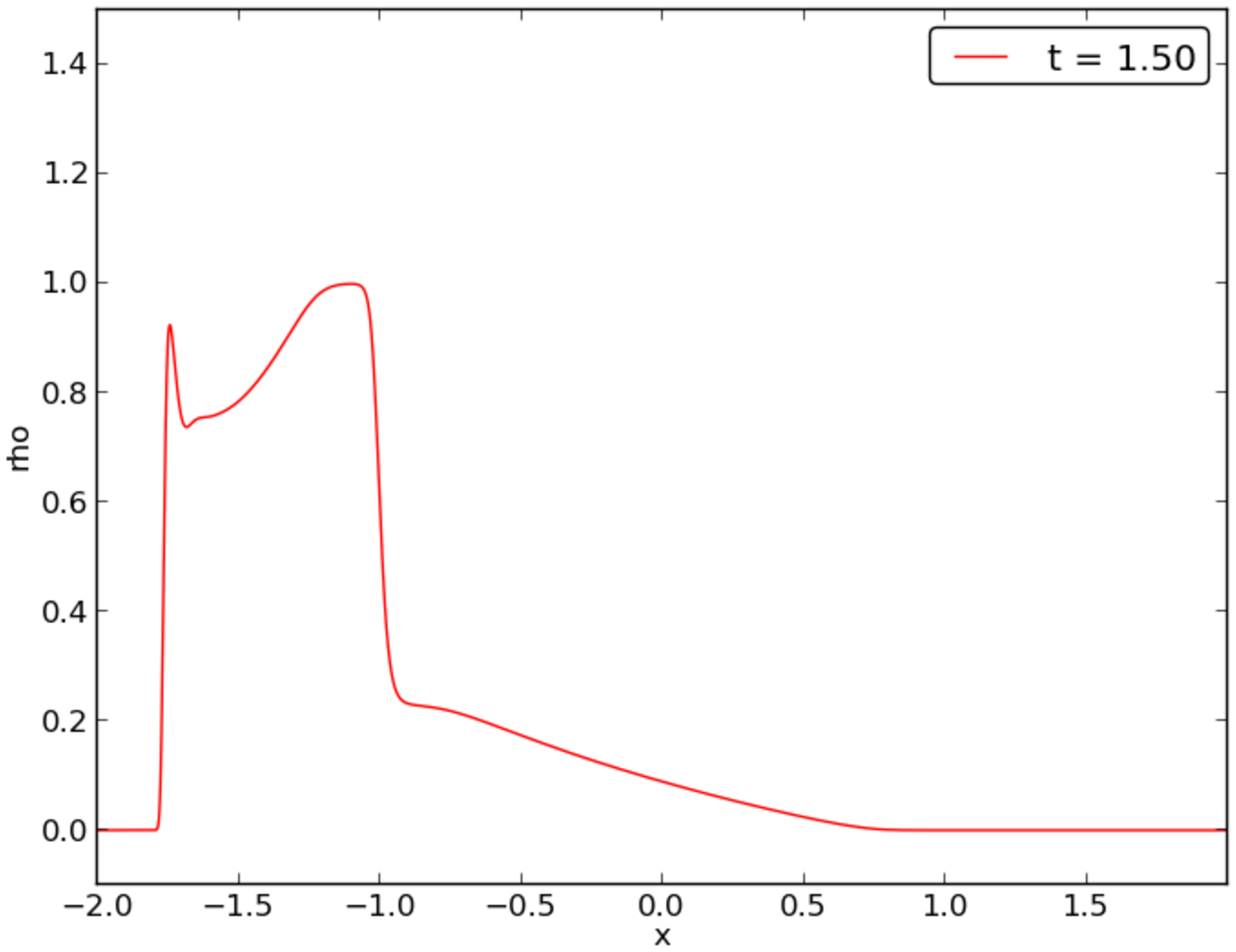}%
  \includegraphics[width=0.25\textwidth, trim=20 0 20 0]{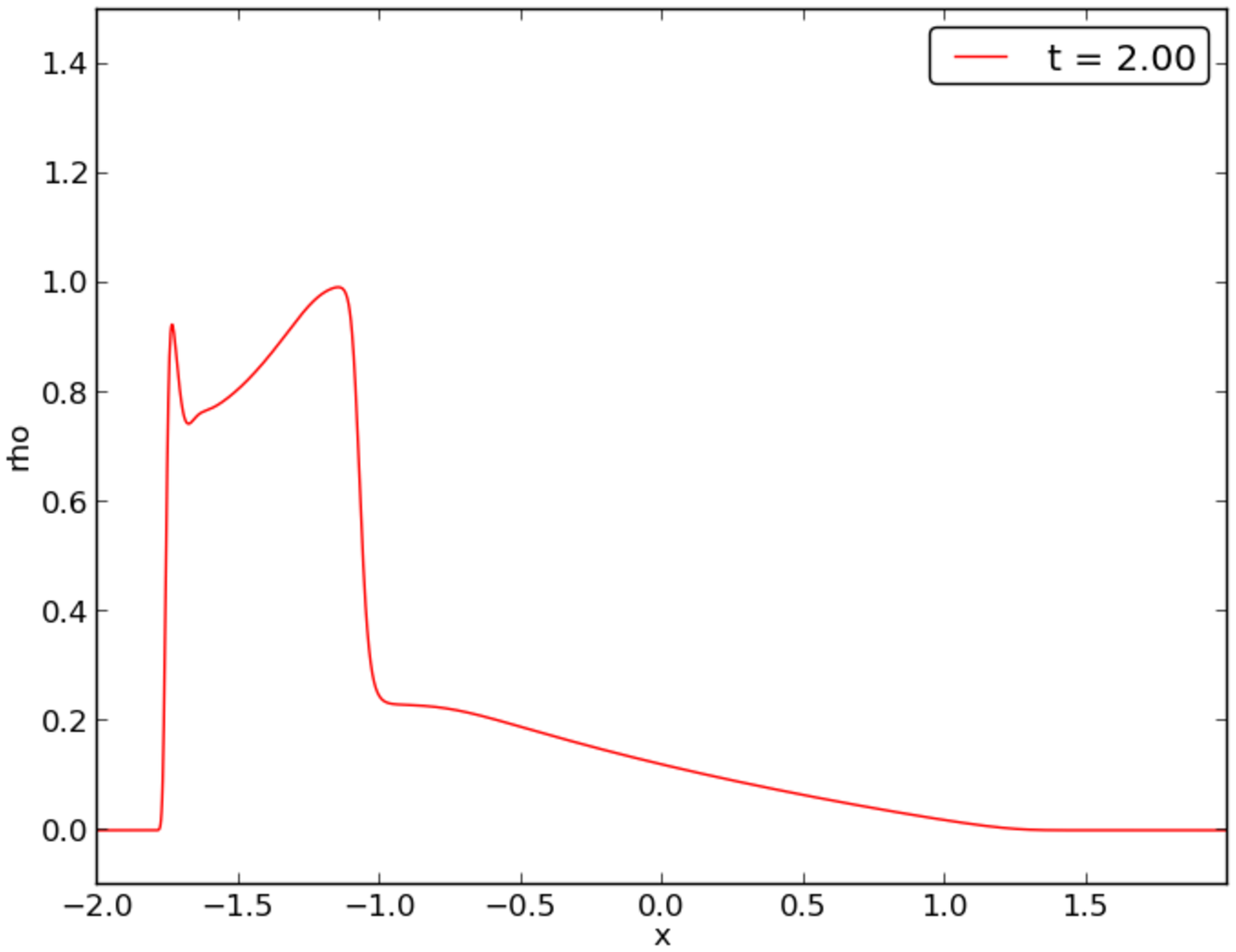}\\
  \includegraphics[width=0.25\textwidth, trim=20 0 20 0]{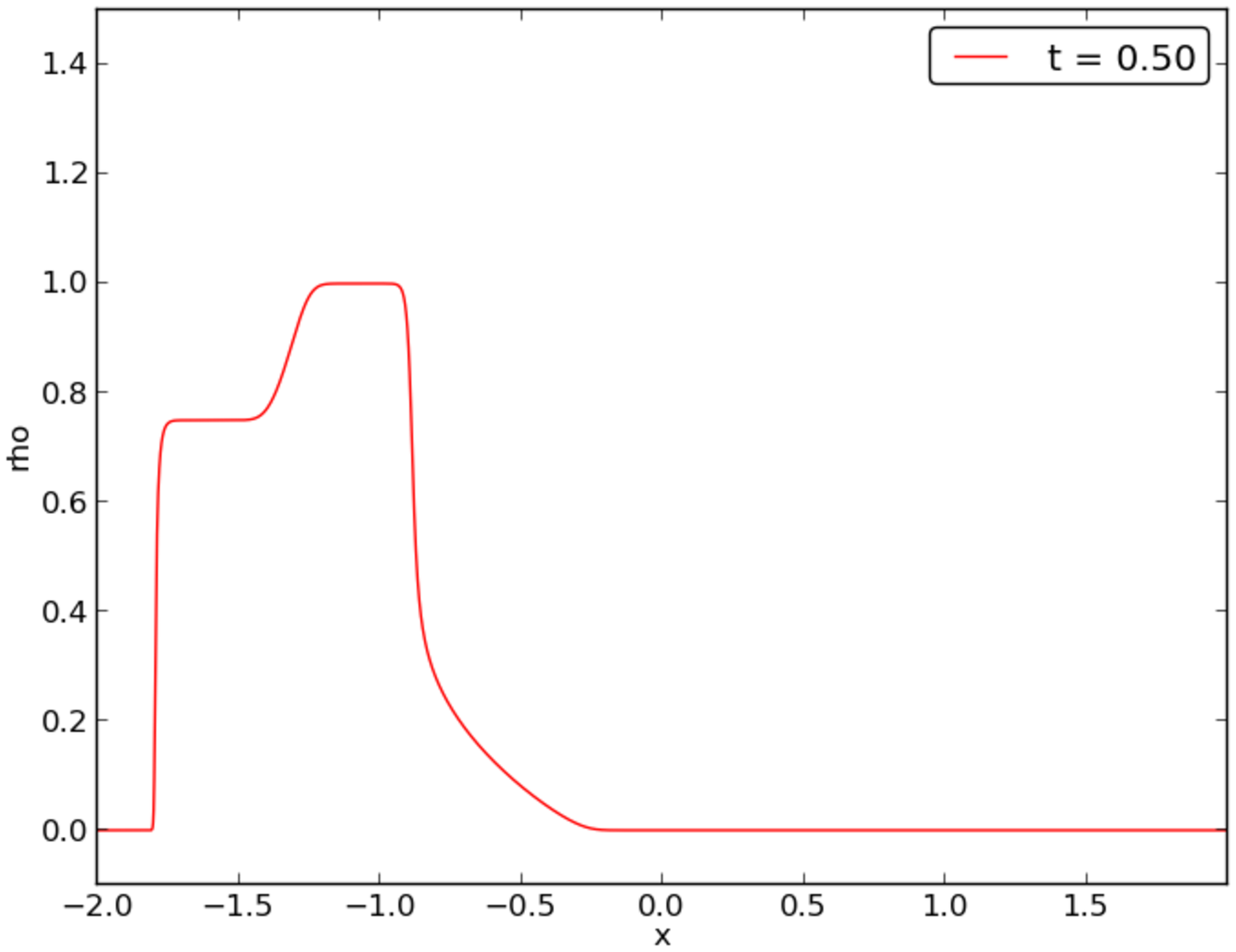}%
  \includegraphics[width=0.25\textwidth, trim=20 0 20 0]{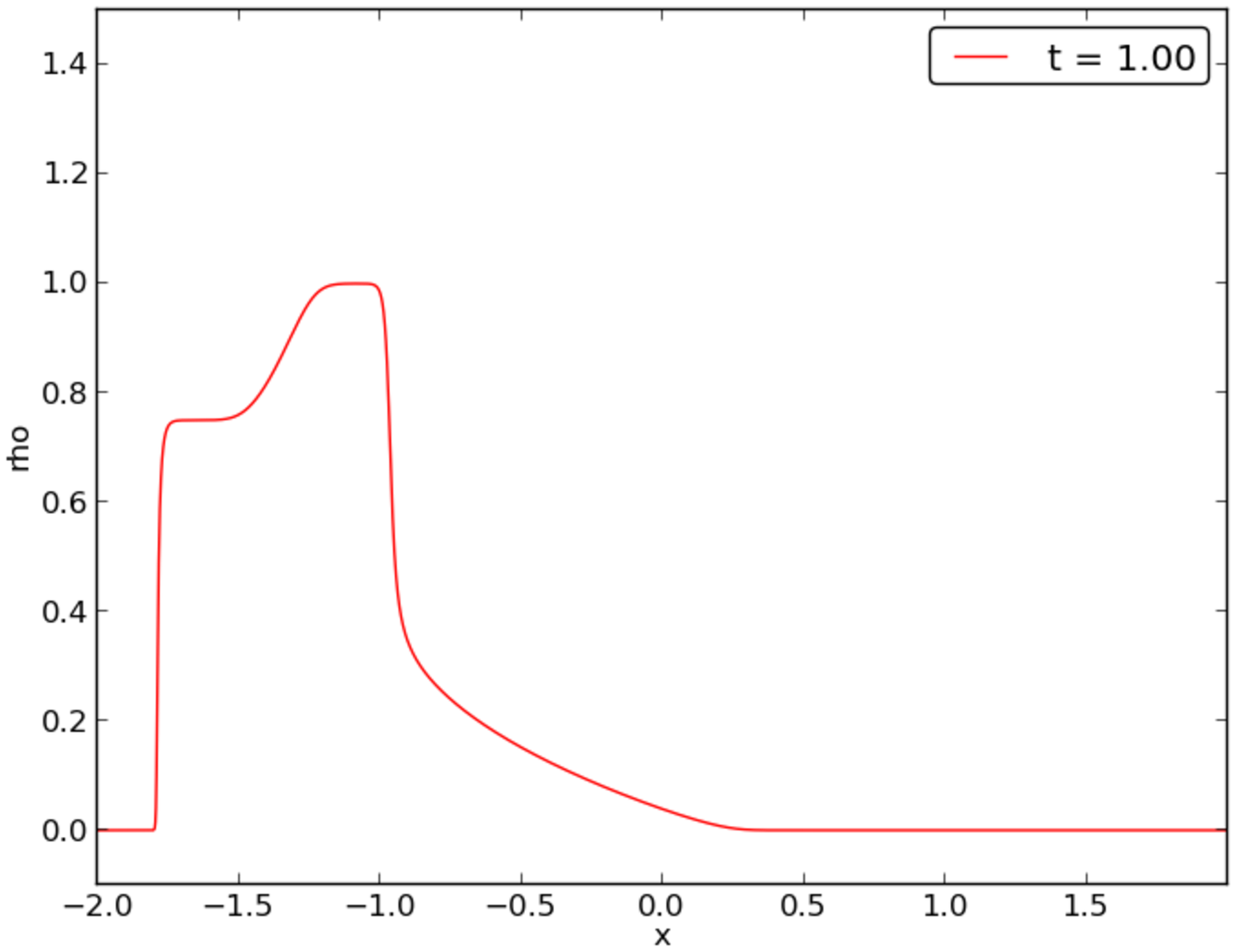}%
  \includegraphics[width=0.25\textwidth, trim=20 0 20 0]{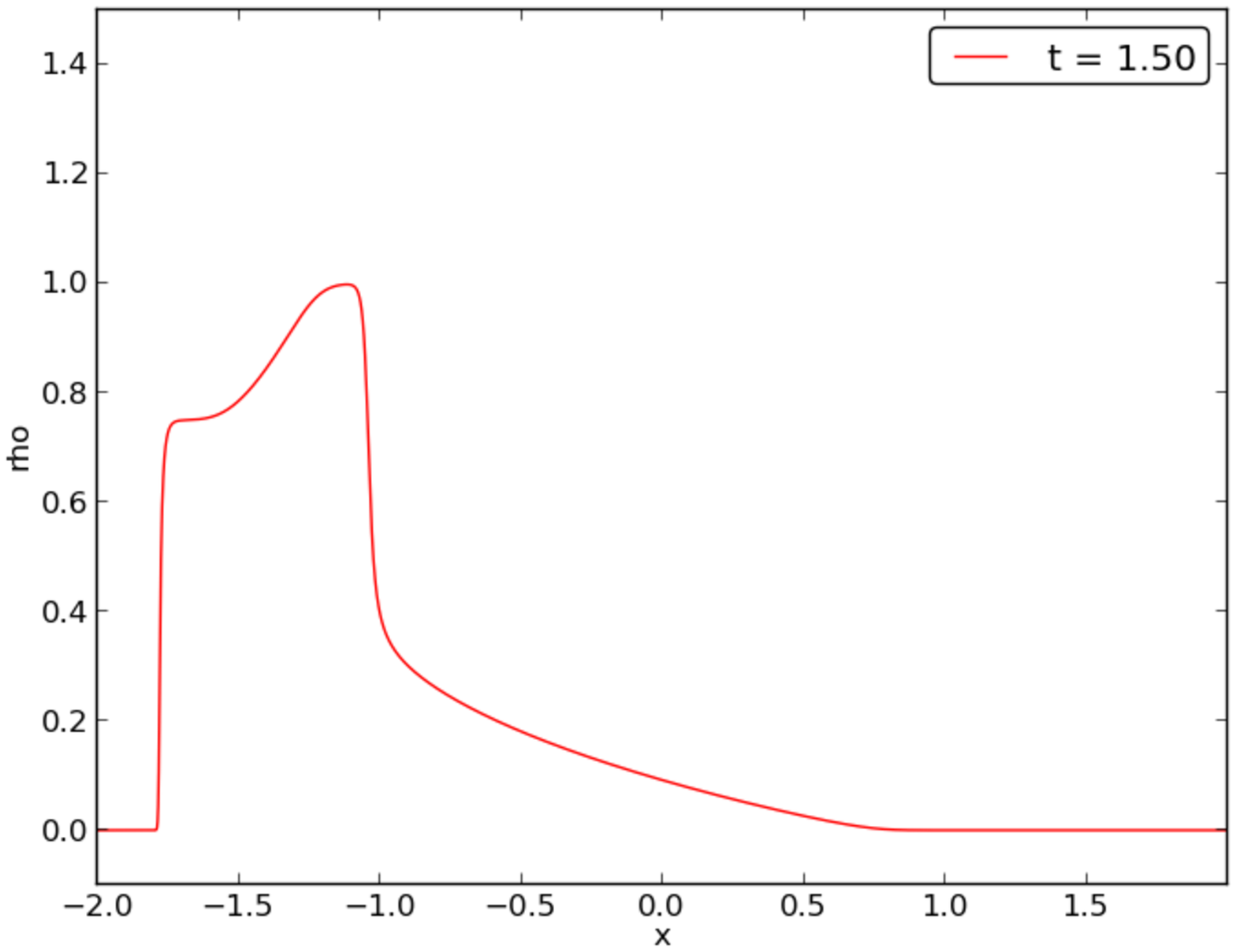}%
  \includegraphics[width=0.25\textwidth, trim=20 0 20 0]{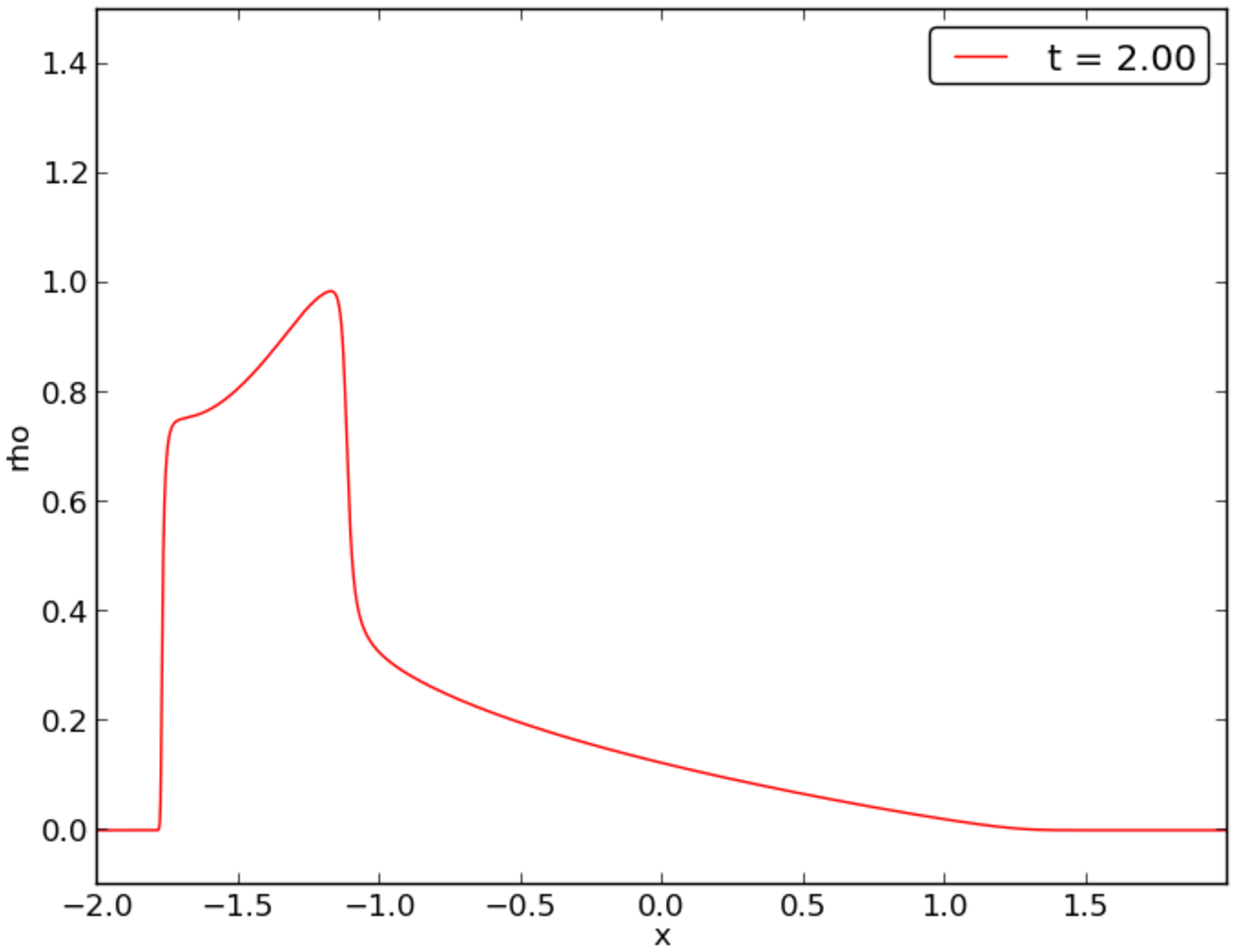}\\
  \caption{Integration of~\eqref{eq:1}--\eqref{eq:2}, first row with
    $a = 0.25$, second with $a = 0.1$, third with $a = 0.05$. On
    the last row, integration of~\eqref{eq:1loc}--\eqref{eq:2}. The
    four columns display the times $t=0.5,\, 1.0,\,1.5$ and $2.0$. The
    mixed waves are due to the non-convex flow~\eqref{eq:2}, see
    Figure~\ref{fig:IDandFlow}.}
  \label{fig:all}
\end{figure}
In the integration below, the solution $\rho$ attains positive values,
so that after an easy modification of $v$ on $\reali^-$ we can assume
that~\eqref{eq:hyp2} holds.

The resulting numerical integrations, carried out satisfying the CFL
condition~\eqref{eq:CFL}, give the diagrams in Figure~\ref{fig:all}.
In the limit case of~\eqref{eq:1loc}, the chosen initial datum leads
to the formation of a rarefaction wave, a shock and a mixed wave, due
to the change of convexity of the flow, see the lowest line in
Figure~\ref{fig:all}. The numerical integrations shown in
Figure~\ref{fig:all} qualitatively suggest that in the limit $a
\to 0$ the solution to~\eqref{eq:1}--\eqref{eq:2} converges to that
of~\eqref{eq:1loc}--\eqref{eq:2}.
\begin{figure}[!htpb]
  \centering
  \includegraphics[width=0.5\textwidth]{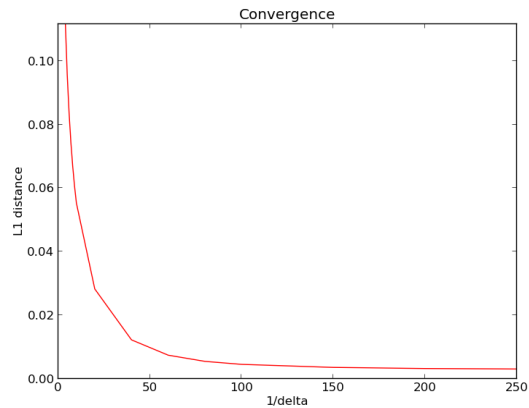}
  \caption{$\L1$-distance between the solution $\rho_a$
    to~\eqref{eq:1}--\eqref{eq:2} for the values of $a$
    in~\eqref{eq:table} and the solution $\rho$
    to~\eqref{eq:1loc}--\eqref{eq:2} at time $t=0.500$.}
  \label{fig:Convergence}
\end{figure}
A more quantitative hint in this direction is in
Figure~\ref{fig:Convergence}. Using the algorithm above, we computed
the solution $\rho_a$ to~\eqref{eq:1}--\eqref{eq:2} for different
values of $a$ and the solution $\rho$
to~\eqref{eq:1loc}--\eqref{eq:2}, all at time
$t=0.500$. Figure~\ref{fig:Convergence} presents the plot of the
$\L1$--distance $\norma{\rho_a - \rho}_{\L1}$ versus $1/a$, see also
table~\eqref{eq:table}.
\begin{equation}
  \label{eq:table}
  \begin{array}{@{}clllll@{}}
    1/a &
    4 &
    5 &
    6 &
    7 &
    8
    \\
    \norma{\rho_a - \rho}_{\L1} &
    0.11166027&
    0.09569174&
    0.08373053&
    0.0743367 &
    0.06674645
    \\[5pt]
    1/a &
    9 &
    10 &
    20 &
    40 &
    60
    \\
    \norma{\rho_a - \rho}_{\L1} &
    0.06049835&
    0.05526474&
    0.0282456 &
    0.0122137 &
    0.0073903
    \\[5pt]
    1/a&
    80 &
    100 &
    150 &
    200 &
    250
    \\
    \norma{\rho_a - \rho}_{\L1} &
    0.0054428 &
    0.00449935&
    0.00354348&
    0.00319927&
    0.0030382
  \end{array}
\end{equation}

\section{Technical Details}
\label{sec:TD}

For any $a,b \in \reali$, we denote $I (a,b) = \left]a,b \right[ \cup
\left]b, a \right[$. We use below the following classical notations:
\begin{displaymath}
  D^+ a_j = a_{j+1} - a_j,
  \qquad
  D^- a_j = a_{j} - a_{j-1},
  \qquad
  D^2 a_j = a_{j+1} -2a_j + a_{j-1} = (D^+ - D^-) a_j
\end{displaymath}
and recall the trivial identities
\begin{displaymath}
  \begin{array}{c}
    D^+ (a_j b_j) = (D^+ a_{j}) b_{j+1} + (D^+ b_{j}) a_{j},
    \qquad\qquad
    D^- (a_j b_j) = (D^- a_{j}) b_{j} + (D^- b_{j}) a_{j-1},
    \\[6pt]
    D^2 (a_j b_j) = (D^2 a_{j}) b_j + (D^2 b_j )a_j + D^+a_j
    D^+b_j + D^-a_j D^-b_j .
  \end{array}
\end{displaymath}
For later use, we note that the algorithm~(\ref{LF}) can then be
rewritten as
\begin{displaymath}
  \rho^{n+1}_j
  =
  \rho^n_j
  -
  \lambda \,
  \frac{D^+ \! \left(f (t^n, x_{j-1/2}, \rho^n_j) \, v (c_{j-1/2})\right)
    +
    D^- \! \left(f (t^n, x_{j+1/2}, \rho^n_j) \, v (c_{j+1/2})\right)}{2}
  +
  \frac{1}{6} \, D^2 \rho^n_j \,.
\end{displaymath}

\begin{proofof}{Lemma~\ref{lem:pos}}
  Note that, by~\eqref{LF}, standard computations lead to
  \begin{equation}
    \label{LINF:10}
    \rho^{n+1}_j
    =
    (1 - \alpha^n_j - \beta^n_j) \rho^n_{j}
    +
    \alpha^n_j \rho^n_{j-1}
    +
    \beta^n_j \rho^n_{j+1}
    -
    \lambda
    \left(
      \mathbf{f}^n_{j+1/2}(\rho^n_{j}, \rho^n_{j})
      -
      \mathbf{f}^n_{j-1/2}(\rho^n_{j}, \rho^n_{j})
    \right).
  \end{equation}
  where
  \begin{equation}
    \label{eq:alpha}
    \alpha^n_j
    =
    \lambda \,
    \frac{
      \mathbf{f}^n_{j+1/2}(\rho^n_{j}, \rho^n_{j+1})
      -
      \mathbf{f}^n_{j+1/2} (\rho^n_{j}, \rho^n_{j})}
    {\rho^n_{j+1} - \rho^n_{j}}
    \quad \mbox{ and } \quad
    \beta^n_j
    =
    \lambda \,
    \frac{
      \mathbf{f}^n_{j-1/2} (\rho^n_{j-1}, \rho^n_{j})
      -
      \mathbf{f}^n_{j-1/2} (\rho^n_{j}, \rho^n_{j})}
    {\rho^n_{j-1} - \rho^n_{j}} \,.
  \end{equation}
  We now show that under condition~\eqref{eq:CFL}, the following
  inequalities hold:
  \begin{equation}
    \label{eq:Cond0}
    \begin{array}{rcl}
      \alpha^n_j & \in & [0, 1/3]
      \\
      \beta^n_j & \in & [0, 1/3]
      \\
      1 - \alpha^n_j - \beta^n_j & \in & [1/3, 1]
    \end{array}
    \qquad \mbox{ and } \qquad
    \lambda \,
    \left(
      \mathbf{f}^n_{j+1/2}(\rho^n_{j},\rho^n_{j})
      -
      \mathbf{f}^n_{j-1/2}(\rho^n_{j}, \rho^n_{j})
    \right)
    \leq
    \frac{1}{3}\, \rho^n_j \,.
  \end{equation}
  Indeed,
  \begin{eqnarray*}
    \alpha^n_j
    & = &
    -
    \lambda \frac{
      \mathbf{f}^n_{j+1/2} (\rho^n_{j}, \rho^n_{j+1})
      -
      \mathbf{f}^n_{j+1/2} (\rho^n_{j}, \rho^n_{j})}{\rho^n_{j+1} - \rho^n_{j}}
    \\
    & = &
    -
    \frac{\lambda}{2}
    \left(
      \frac{
        f(t^n, x_{j+1/2},\rho^n_{j+1})
        -
        f(t^n, x_{j+1/2},\rho^n_{j})}{\rho^n_{j+1}-\rho^n_{j}} \, v(c^n_{j+1/2})
      -
      \frac{1}{3\, \lambda}
    \right)
    \\
    & =&
    -
    \frac{\lambda}{2} \; \partial_\rho f(t^n, x_{j+1/2},\zeta_{j+1/2}) \;
    v(c^n_{j+1/2})
    +
    \frac{1}{6},
  \end{eqnarray*}
  So that
  \begin{eqnarray}
    \nonumber
    \alpha^n_j
    & \geq &
    \displaystyle
    \frac{1}{2}
    \left(
      \frac{1}{3}
      -
      \lambda \, \norma{\partial_\rho f}_{\L\infty} \, \norma{v}_{\L\infty}
    \right)
    \; \geq \; 0
    \\
    \label{eq:bound6}
    \alpha^n_j
    & \leq &
    \displaystyle
    \frac{1}{2}
    \left(
      \frac{1}{3}
      +
      \lambda \, \norma{\partial_\rho f}_{\L\infty} \, \norma{v}_{\L\infty}
    \right)
    \; \leq \;
    \displaystyle
    \frac{1}{6}
    % & \leq &
    % \displaystyle
    % \frac{1}{3}
    \,.
  \end{eqnarray}
  Entirely similar computations lead to analogous estimates for
  $\beta^n_j$. The bounds on $1-\alpha^n_j-\beta^n_j$ follow. The
  last term in~\eqref{eq:Cond0}, using~\eqref{LFa}
  and~\eqref{eq:lambda}, is estimated as follows
  \begin{eqnarray*}
    \!\!\!
& &
    \mathbf{f}^n_{j+1/2}(\rho^n_{j},\rho^n_{j})
    -
    \mathbf{f}^n_{j-1/2}(\rho^n_{j}, \rho^n_{j})
    \\
    \!\!\!
    & \leq &
    \modulo{
      f(t^n,x_{j+1/2},\rho^n_{j}) \, v (c^n_{j+1/2})
      -
      f(t^n,x_{j-1/2},\rho^n_{j}) \, v (c^n_{j-1/2})
    }
    \\
    \!\!\!
    & \leq &
    \modulo{
      f(t^n,x_{j+1/2},\rho^n_{j})
      -
      f(t^n,x_{j-1/2},\rho^n_{j})} v (c^n_{j-1/2})
    +
    \modulo{f(t^n,x_{j-1/2},\rho^n_{j})}
    \modulo{v (c^n_{j+1/2}) - v (c^n_{j-1/2})}
    \\
    \!\!\!
    & \leq &
    h \, \modulo{\partial_xf (t^n,\zeta_j,\rho^n_j)} \, \norma{v}_{\L\infty}
    +
    2 \, \norma{\partial_\rho f}_{\L\infty} \,
    \norma{v}_{\L\infty} \, \modulo{\rho^n_j}
    \\
    \!\!\!
    & \leq &
    \left(
      C \, \norma{v}_{\L\infty} \, h
      +
      2 \, \norma{\partial_\rho f}_{\L\infty} \, \norma{v}_{\L\infty}
    \right)
    \modulo{\rho^n_j}
    \\
    \!\!\!
    & \leq &
    \left(1 + 2 \, \norma{\partial_\rho f}_{\L\infty} \right)
    \norma{v}_{\L\infty} \, \modulo{\rho^n_j}
    \\
    \!\!\!
    & \leq &
    \frac{1}{3\, \lambda} \, \rho^n_j \,.
  \end{eqnarray*}
  Using the bounds~\eqref{eq:Cond0} in~\eqref{LINF:10}, we obtain
  \begin{displaymath}
    \rho^{n+1}_j
    \geq
    (1 - \alpha^n_j - \beta^n_j) \, \rho^n_{j}
    +
    \alpha^n_j \, \rho^n_{j-1}
    +
    \beta^n_j \, \rho^n_{j+1}
    -
    \frac{1}{3} \, \rho^n_j
    \geq
    0 \,,
  \end{displaymath}
  proving the positivity of the discrete solution.
\end{proofof}

\begin{proofof}{Lemma~\ref{lem:L1}}
  Thanks to the positivity of the discrete solution, it is sufficient
  to compute
  \begin{eqnarray*}
    \norma{\rho^{n+1}}_{\L1}
    & = &
    \sum_j h \, \rho^{n+1}_j
    \\
    & = &
    \sum_j
    h
    \left(
      \rho^n_j
      -
      \lambda
      \left(
        \mathbf{f}^n_{j+1/2}(\rho^n_j, \rho^n_{j+1})
        -
        \mathbf{f}^n_{j-1/2}(\rho^n_{j-1}, \rho^n_{j})
      \right)
    \right)
    \\
    & = &
    \sum_j h \, \rho^n_j
    -
    h \, \lambda \, \left(
      \lim_{i\to -\infty} \mathbf{f}^n_{i+1/2}(\rho^n_i, \rho^n_{i+1})
      -
      \lim_{i\to+\infty} \mathbf{f}^n_{i-1/2}(\rho^n_{i-1}, \rho^n_{i})
    \right)
    \\
    & = &
    \norma{\rho^n}_{\L1}
  \end{eqnarray*}
  completing the proof.
\end{proofof}

\begin{proofof}{Lemma~\ref{lem:Linfty}}
  For later use, estimate the quantity
  \begin{eqnarray*}
    \modulo{c^n_{j+1/2} - c^n_{j-1/2}}
    & \leq &
    \sum_{k\in\interi} h \,
    \modulo{\rho^n_{k+1/2} (\eta_{k - (j+1/2)} - \eta_{k - (j-1/2)})}
    \\
    & \leq &
    \sum_{k\in\interi} h \, \rho^n_{k+1/2}
    \int_{x_{k-j-1/2}}^{x_{k-j+1/2}} \modulo{\eta'(s)} \, ds
    \\
    & \leq &
    h \, \norma{\rho^n}_{\L1} \, \norma{\eta'}_{\L\infty}
    \\
    & = &
    h \, \norma{\rho^o}_{\L1} \, \norma{\eta'}_{\L\infty} \, ,
  \end{eqnarray*}
  where Lemma~\ref{lem:L1} was used.  Using the same estimates as in
  the proof of Lemma~\ref{lem:pos}, equality~\eqref{LINF:10} yields
  \begin{eqnarray*}
    \nonumber
    \rho^{n+1}_{j}
    & \leq &
    % (1 - \alpha^n_j - \beta^n_j) \, \rho^{n}_{j}
    % +
    % \alpha^n_j \, \rho^{n}_{j-1}
    % +
    % \beta^n_j \, \rho^{n}_{j+1}
    % +
    % \lambda \, \modulo{
    % \mathbf{f}^n_{j+1/2}(\rho^n_{j}, \rho^n_{j})
    % -
    % \mathbf{f}^n_{j-1/2}(\rho^n_{j}, \rho^n_{j})
    % }
    %   \\
    %   \nonumber
    %   & = &
    %   (1 - \alpha^n_j - \beta^n_j) \, \rho^{n}_{j}
    %   +
    %   \alpha^n_j \, \rho^{n}_{j-1}
    %   +
    %   \beta^n_j \, \rho^{n}_{j+1}
    %   \\
    %   \nonumber
    %   & &
    %   +
    %   \lambda
    %   \modulo{
    %   f(t^n,x_{j+1/2},\rho^n_{j}) \, v (c^n_{j+1/2})
    %   -
    %   f(t^n,x_{j-1/2},\rho^n_{j}) \, v (c^n_{j-1/2})
    % }
    %   \\
    %   \nonumber
    %   & = &
    (1 - \alpha^n_j - \beta^n_j) \, \rho^{n}_{j}
    +
    \alpha^n_j \, \rho^{n}_{j-1}
    +
    \beta^n_j \, \rho^{n}_{j+1}
    \\
    \nonumber
    & &
    +
    \lambda
    \modulo{
      f(t^n,x_{j+1/2},\rho^n_{j})
      -
      f(t^n,x_{j-1/2},\rho^n_{j})} \, v (c^n_{j-1/2})
    \\
    \nonumber
    & &
    +
    \lambda
    \modulo{f(t^n,x_{j-1/2},\rho^n_{j})} \,
    \modulo{D^-v (c^n_{j+1/2})}
    \\
    \nonumber
    & \leq &
    (1 - \alpha^n_j - \beta^n_j) \, \norma{\rho^n}_{\L\infty}
    +
    \alpha^n_j \, \norma{\rho^n}_{\L\infty}
    +
    \beta^n_j \, \norma{\rho^n}_{\L\infty}
    \\
    \nonumber
    & &
    +
    \lambda
    \left(
      h \, \modulo{\partial_xf (t^n,\zeta_j,\rho^n_j)} \, \norma{v}_{\L\infty}
      +
      \rho^n_j \, \norma{\partial_\rho f}_{\L\infty} \,
      \norma{v'}_{\L\infty} \, \modulo{c^n_{j+1/2} - c^n_{j-1/2}}
    \right)
    \\
    \nonumber
    & \leq &
    \norma{\rho^n}_{\L\infty}
    +
    \tau
    \left(
      C \, \norma{v}_{\L\infty}
      +
      \norma{\partial_\rho f}_{\L\infty} \,
      \norma{v'}_{\L\infty} \,
      \norma{\rho^o}_{\L1} \, \norma{\eta'}_{\L\infty}
    \right)
    \norma{\rho^n}_{\L\infty}
    \\
    & \leq &
    e^{\mathcal{L}\tau} \,     \norma{\rho^n}_{\L\infty}
  \end{eqnarray*}
  provided
  \begin{equation}
    \label{eq:L}
    \mathcal{L}
    =
    C \, \norma{v}_{\L\infty}
    +
    \norma{\partial_\rho f}_{\L\infty} \, \norma{v'}_{\L\infty} \,
    \norma{\rho^o}_{\L1} \, \norma{\eta'}_{\L\infty} \,.
  \end{equation}
  A standard iterative argument completes the proof.
\end{proofof}

\begin{proofof}{Proposition~\ref{prop:TV}}
  First, we write~\eqref{LF} for $j$ and for $j+1$, subtract and get
  \begin{eqnarray*}
    & &
    \rho^{n+1}_{j+1} - \rho^{n+1}_j
    =
    \rho^n_{j+1} - \rho^n_j
    \\
    & &
    -
    \lambda
    \left(
      \mathbf{f}^n_{j+3/2} (\rho^{n}_{j+1}, \rho^{n}_{j+2})
      -
      \mathbf{f}^n_{j+1/2} (\rho^{n}_{j}, \rho^{n}_{j+1})
      +
      \mathbf{f}^n_{j+1/2} (\rho^{n}_{j}, \rho^{n}_{j+1})
      -
      \mathbf{f}^n_{j-1/2} (\rho^{n}_{j-1}, \rho^{n}_{j})
    \right).
  \end{eqnarray*}
  Now add and subtract $\mathbf{f}^n_{j+3/2} (\rho^{n}_{j},
  \rho^{n}_{j+1}) + \mathbf{f}^n_{j+1/2} (\rho^{n}_{j-1},
  \rho^{n}_{j})$, then rearrange to obtain
  \begin{equation}
    \label{eq:AB}
    \rho^{n+1}_{j+1} - \rho^{n+1}_j
    =
    \mathcal{A}^n_j - \lambda \, \mathcal{B}^n_j
  \end{equation}
  where
  \begin{eqnarray}
    \nonumber
    \mathcal{A}^n_j
    & = &
    \rho^n_{j+1} - \rho^n_j
    \\
    \nonumber
    & &
    -\lambda
    \left(
      \mathbf{f}^n_{j+3/2} (\rho^{n}_{j+1}, \rho^{n}_{j+2})
      -
      \mathbf{f}^n_{j+1/2} (\rho^{n}_{j}, \rho^{n}_{j+1})
      -
      \mathbf{f}^n_{j+3/2} (\rho^{n}_{j}, \rho^{n}_{j+1})
      +
      \mathbf{f}^n_{j+1/2} (\rho^{n}_{j-1}, \rho^{n}_{j})
    \right)
    % \\
    % A
    % & = &
    % \mathbf{f}^n_{j+3/2} (\rho^{n}_{j+1}, \rho^{n}_{j+2})
    % -
    % \mathbf{f}^n_{j+1/2} (\rho^{n}_{j}, \rho^{n}_{j+1})
    % -
    % \mathbf{f}^n_{j+3/2} (\rho^{n}_{j}, \rho^{n}_{j+1})
    % +
    % \mathbf{f}^n_{j+1/2} (\rho^{n}_{j-1}, \rho^{n}_{j})
    \\
    \label{eq:B}
    \mathcal{B}^n_j
    & = &
    \mathbf{f}^n_{j+3/2} (\rho^{n}_{j}, \rho^{n}_{j+1})
    -
    \mathbf{f}^n_{j+1/2} (\rho^{n}_{j}, \rho^{n}_{j+1})
    +
    \mathbf{f}^n_{j-1/2} (\rho^{n}_{j-1}, \rho^{n}_{j})
    -
    \mathbf{f}^n_{j+1/2} (\rho^{n}_{j-1}, \rho^{n}_{j})
    \,.
  \end{eqnarray}

  Consider first the term $\mathcal{A}^n_j$.  % We
  % claim that
  % \begin{equation}
  %   \label{eq:90}
  %   \sum_{j\in\interi} \modulo{\mathcal{A}_j}
  %   \leq
  %   (1 + \mathcal{L} \, \tau)
  %   \sum_{j\in\interi} \modulo{\rho^{n}_{j+1} - \rho^{n}_{j}} \,.
  % \end{equation}
  Recall~\eqref{eq:alpha} and observe that, after suitable
  rearrangements,
  \begin{eqnarray}
    \nonumber
    & &
    \mathcal{A}^n_j
    \\
    \nonumber
    & = &
    \frac{2}{3} \, (\rho^n_{j+1} - \rho^n_j)
    \\
    \nonumber
    & + &
    (\rho^n_{j+2} - \rho^n_{j+1})
    \left(
      \frac{1}{6}
      -
      \frac{\lambda}{2} \,
      \frac{
        f (t^n,x_{j+3/2},\rho^n_{j+2}) \, v (c^n_{j+3/2})
        -
        f (t^n,x_{j+1/2},\rho^n_{j+1}) \, v (c^n_{j+1/2})
      }{\rho^n_{j+2} - \rho^n_{j+1}}
    \right)
    \\
    \nonumber
    & + &
    (\rho^n_{j} - \rho^n_{j-1})
    \left(
      \frac{1}{6}
      +
      \frac{\lambda}{2} \,
      \frac{
        f (t^n,x_{j+3/2},\rho^n_{j}) \, v (c^n_{j+3/2})
        -
        f (t^n,x_{j+1/2},\rho^n_{j-1}) \, v (c^n_{j+1/2})
      }{\rho^n_{j} - \rho^n_{j-1}}
    \right)
    \\
    \nonumber
    & = &
    \frac{2}{3} \, (\rho^n_{j+1} - \rho^n_j)
    \\
    \label{eq:u1}
    & + &
    \!\!
    (\rho^n_{j+2} - \rho^n_{j+1})
    \left[
      \frac{1}{6}
      -
      \frac{\lambda}{2}
      \frac{
        f (t^n,x_{j+3/2},\rho^n_{j+2}) v (c^n_{j+3/2})
        -
        f (t^n,x_{j+3/2},\rho^n_{j+1}) v (c^n_{j+3/2})
      }{\rho^n_{j+2} - \rho^n_{j+1}}
    \right]
    \\
    \label{eq:u2}
    & + &
    (\rho^n_{j} - \rho^n_{j-1})
    \left[
      \frac{1}{6}
      +
      \frac{\lambda}{2} \,
      \frac{
        f (t^n,x_{j+3/2},\rho^n_{j}) \, v (c^n_{j+3/2})
        -
        f (t^n,x_{j+3/2},\rho^n_{j-1}) \, v (c^n_{j+3/2})
      }{\rho^n_{j} - \rho^n_{j-1}}
    \right]
    \\
    \label{eq:u3}
    & + &
    -
    \frac{\lambda}{2}
    \left(
      f (t^n,x_{j+3/2},\rho^n_{j+1})
      -
      f (t^n,x_{j+3/2},\rho^n_{j-1})
    \right)
    v (c^n_{j+3/2})
    \\
    \label{eq:u4}
    & + &
    \frac{\lambda}{2}
    \left(
      f (t^n,x_{j+1/2},\rho^n_{j+1})
      -
      f (t^n,x_{j+1/2},\rho^n_{j-1})
    \right)
    v (c^n_{j+1/2}).
  \end{eqnarray}
  Remark that the term~\eqref{eq:u1} equals $\alpha^n_{j+1}$, as
  defined in~\eqref{eq:alpha}. Hence, it can be bounded
  using~\eqref{eq:Cond0} as follows:
  \begin{displaymath}
    \modulo{\eqref{eq:u1}} \leq \frac{1}{6} \, \modulo{\rho^n_{j+2}-\rho^n_{j+1}}.
  \end{displaymath}
  The estimate for the term~\eqref{eq:u2} is exactly as that of
  $\alpha^n_j$ in Lemma~\ref{lem:pos}, so that
  \begin{displaymath}
    \modulo{\eqref{eq:u2}} \leq \frac{1}{6} \, \modulo{\rho^n_{j}-\rho^n_{j-1}}.
  \end{displaymath}
  Consider now the terms~\eqref{eq:u3} and~\eqref{eq:u4}:
  \begin{eqnarray*}
    \eqref{eq:u3} + \eqref{eq:u4}
    % & &
    % -
    % \left(
    %   f (t^n,x_{j+3/2},\rho^n_{j+1})
    %   -
    %   f (t^n,x_{j+3/2},\rho^n_{j-1})
    % \right)
    % v (c^n_{j+3/2})
    % \\
    % & &
    % \qquad
    % +
    % \left(
    %   f (t^n,x_{j+1/2},\rho^n_{j+1})
    %   -
    %   f (t^n,x_{j+1/2},\rho^n_{j-1})
    % \right)
    % v (c^n_{j+1/2})
    % \\
    & = &
    \frac{\lambda}{2}
    \int_{\rho^n_{j-1}}^{\rho^n_{j+1}}
    \left(
      - \partial_\rho f (t^n,x_{j+3/2},r)  \, v (c^n_{j+3/2})
      + \partial_\rho f (t^n,x_{j+1/2},r) \, v (c^n_{j+1/2})
    \right)
    \d r
    \\
    & = &
    -\frac{\lambda}{2}
    \int_{\rho^n_{j-1}}^{\rho^n_{j+1}}
    \left(
      \partial_\rho f (t^n,x_{j+3/2},r)
      -
      \partial_\rho f (t^n,x_{j+1/2},r)
    \right)
    \d r
    \; v (c^n_{j+3/2})
    \\
    & &
    \qquad
    - \frac{\lambda}{2}
    \int_{\rho^n_{j-1}}^{\rho^n_{j+1}}
    \partial_\rho f (t^n,x_{j+1/2},r)
    \d r
    \; \left(
      v (c^n_{j+3/2})
      -
      v (c^n_{j+1/2})
    \right)
    \\
    & = &
    - \frac{\lambda}{2}
    \int_{\rho^n_{j-1}}^{\rho^n_{j+1}}
    \int_{x_{j+1/2}}^{x_{j+3/2}}
    \partial^2_{x,\rho} f (t^n,\xi,r) \,
    \d\xi \, \d r
    \; v (c^n_{j+3/2})
    \\
    & &
    - \frac{\lambda}{2}
    \int_{\rho^n_{j-1}}^{\rho^n_{j+1}}
    \partial_\rho f (t^n,x_{j+1/2},r)
    \d r
    \; \left(
      v (c^n_{j+3/2})
      -
      v (c^n_{j+1/2})
    \right)
  \end{eqnarray*}
  So that, passing to the absolute value
  \begin{displaymath}
    \modulo{\eqref{eq:u3} + \eqref{eq:u4}}
    % & &
    % \Big|
    % -
    % \left(
    %   f (t^n,x_{j+3/2},\rho^n_{j+1})
    %   -
    %   f (t^n,x_{j+3/2},\rho^n_{j-1})
    % \right)
    % v (c^n_{j+3/2})
    % \\
    % & &
    % \qquad
    % +
    % \left(
    %   f (t^n,x_{j+1/2},\rho^n_{j+1})
    %   -
    %   f (t^n,x_{j+1/2},\rho^n_{j-1})
    % \right)
    % v (c^n_{j+1/2})
    % \Big|
    % \\
    \leq
    \frac12 \left[
      \norma{\partial^2_{x,\rho}f}_{\L\infty} \, \norma{v}_{\L\infty}
      +
      \norma{\partial_\rho f}_{\L\infty} \, \norma{v'}_{\L\infty}\,
      \norma{\rho^o}_{\L1} \, \norma{\eta'}_{\L\infty}
    \right]
    \tau \,
    \modulo{\rho^n_{j+1} - \rho^n_{j-1}}
  \end{displaymath}
  Grouping the estimates obtained we get
  \begin{equation}
    \label{eq:90}
    \begin{array}{@{}rcl@{}}
      \sum_{j \in \interi} \modulo{\mathcal{A}^n_j}
      & \leq &
      \left[
        1
        +
        \frac12 \left(
          \norma{\partial^2_{x,\rho}f}_{\L\infty} \,
          \norma{v}_{\L\infty} + \norma{\partial_\rho f}_{\L\infty} \,
          \norma{v'}_{\L\infty} \, \norma{\rho^o}_{\L1} \, \norma{\eta'}_{\L\infty}
        \right)\, \tau
      \right]
      \\
      & &
      \qquad\qquad
      \times
      \sum_{j \in \interi} \modulo{\rho^n_{j+1} - \rho^n_j} \,.
    \end{array}
  \end{equation}

  We now turn to the term $\mathcal{B}^n_j$ in~\eqref{eq:B}. Since
  {\small
    \begin{eqnarray}
      \label{eq:uu1}
      \!\!\!\!\!\!\!\!\!
      \mathcal{B}^n_j
      & = &
      \frac{
        f (t^n,x_{j+3/2},\rho^n_j) v (c_{j+3/2})
        -
        2 f (t^n,x_{j+1/2},\rho^n_j) v (c_{j+1/2})
        +
        f (t^n,x_{j-1/2},\rho^n_j) v (c_{j-1/2})
      }{2}
      \\
      \label{eq:uu2}
      & &
      +
      \frac{
        f (t^n,x_{j+3/2},\rho^n_{j+1}) \, v (c_{j+3/2})
        -
        f (t^n,x_{j+1/2},\rho^n_{j+1}) \, v (c_{j+1/2})
      }{2}
      \\
      \label{eq:uu3}
      & &
      -
      \frac{
        f (t^n,x_{j+1/2},\rho^n_{j-1}) \, v (c_{j+1/2})
        -
        f (t^n,x_{j-1/2},\rho^n_{j-1}) \, v (c_{j-1/2})
      }{2}\,,
    \end{eqnarray}
  }we consider the various terms separately.
  \begin{eqnarray*}
    \eqref{eq:uu1}
    & = &
    v (c_{j+1/2}) \;
    \frac{
      f (t^n, x_{j+3/2}, \rho^n_j)
      -
      2 f (t^n, x_{j+1/2}, \rho^n_j)
      +
      f (t^n, x_{j-1/2}, \rho^n_j)
    }{2}
    \\
    & &
    +
    f (t^n, x_{j+1/2}, \rho^n_j) \;
    \frac{v (c_{j+3/2}) - 2 v (c_{j+1/2}) + v (c_{j-1/2})}{2}
    \\
    & &
    +
    \frac{f (t^n, x_{j+3/2}, \rho^n_j) - f (t^n, x_{j+1/2}, \rho^n_j)}{2} \;
    \frac{v (c_{j+3/2}) - v (c_{j+1/2})}{2}
    \\
    & &
    +
    \frac{f (t^n, x_{j+1/2}, \rho^n_j) - f (t^n, x_{j-1/2}, \rho^n_j)}{2} \;
    \frac{v (c_{j+1/2}) - v (c_{j-1/2})}{2}
    \\
  \end{eqnarray*}
  where
  \begin{eqnarray*}
    & &
    \modulo{\frac{
        f (t^n, x_{j+3/2}, \rho^n_j)
        -
        2 f (t^n, x_{j+1/2}, \rho^n_j)
        +
        f (t^n, x_{j-1/2}, \rho^n_j)
      }{2}
    }
    \\
    & = &
    \modulo{\frac{
        f (t^n, x_{j+3/2}, \rho^n_j)
        -
        f (t^n, x_{j+1/2}, \rho^n_j)
      }{2}
      -
      \frac{
        f (t^n, x_{j+1/2}, \rho^n_j)
        +
        f (t^n, x_{j-1/2}, \rho^n_j)
      }{2}
    }
    \\
    & \leq &
    \frac{h}{2} \,
    \modulo{
      \partial_x f (t^n, \zeta_{j+1}, \rho^n_j)
      -
      \partial_x f (t^n, \zeta_{j-1}, \rho^n_j)}
    \\
    & \leq &
    \frac{h}{2} \,
    \int_{\zeta_{j-1}}^{\zeta_{j+1}}
    \modulo{ \partial^2_{xx} f (t_x, x, \rho^n_j)} \d x
    \\
    & \leq &
    C \, h^2 \, \modulo{\rho^n_j}
  \end{eqnarray*}
  where~\eqref{eq:hyp1} was used to get to the last line.  Moreover,
  \begin{eqnarray*}
    & &
    \frac{v (c_{j+3/2}) - 2 v (c_{j+1/2}) + v (c_{j-1/2})}{2}
    \\
    & = &
    \frac{v (c_{j+3/2}) - v (c_{j+1/2})}{2}
    -
    \frac{v (c_{j+1/2}) + v (c_{j-1/2})}{2}
    \\
    & = &
    \frac{1}{2}
    \left(
      v'(\zeta_j) ( c^n_{j+1/2} - c^n_{j-1/2})
      -
      v'(\zeta_{j+1})( c^n_{j+3/2} - c^n_{j+1/2})
    \right)
    \\
    & = &
    \frac{1}{2}
    \left(v' (\zeta_j) - v' (\zeta_{j+1})\right) (c^n_{j+1/2} - c^n_{j-1/2})
    -
    \frac{1}{2}
    v' (\zeta_{j+1}) (c^n_{j+3/2} + 2c^n_{j+1/2} - c^n_{j-1/2})
    \\
    & = &
    \frac{1}{2} \,
    v'' (\xi_j) \, (\zeta_{j} - \zeta_{j+1}) \, (c^n_{j+1/2} - c^n_{j-1/2})
    -
    \frac{1}{2} \,
    v' (\zeta_{j+1}) \, (c^n_{j+3/2} + 2c^n_{j+1/2} - c^n_{j-1/2}) \,.
  \end{eqnarray*}
  Note that we have $\modulo{\zeta_{j} - \zeta_{j+1}} \leq
  \modulo{c^n_{j+3/2} - c^n_{j+1/2}} + \modulo{c^n_{j+1/2} -
    c^n_{j-1/2}}$, and so using Young's inequality,
  \begin{eqnarray*}
    \!
    \modulo{\frac{v (c_{j+3/2}) - 2 v (c_{j+1/2}) + v (c_{j-1/2})}{2}}
    \!\!\!& \leq &\!\!\!
    \frac{1}{2}
    \norma{v''}_{\L\infty}
    \left[
      \frac32 \modulo{c^n_{j+1/2} - c^n_{j-1/2}}^2
      +
      \frac12 \modulo{c^n_{j+3/2} - c^n_{j+1/2}}^2
    \right]
    \!
    \\
    & &
    +
    \frac{1}{2}
    \norma{v'}_{\L\infty} \, \modulo{c^n_{j+3/2} + 2c^n_{j+1/2} - c^n_{j-1/2}} \,.
  \end{eqnarray*}
  We now estimate the terms involving the discrete derivatives of
  $c^n_j$ in the expression above, exploiting the
  regularity~\eqref{eq:hyp2} of $\eta$. By~\eqref{eq:cnj}, we have
  \begin{eqnarray}
    \nonumber
    \modulo{c^n_{j+1/2} - c^n_{j-1/2}}
    & = &
    \modulo{\sum_{k\in\interi} h \, \rho^n_{k+1/2} \,
      (\eta_{k - (j+1/2)} - \eta_{k - (j-1/2)})}
    \\
    \nonumber
    & \leq &
    \sum_{k\in\interi}
    h \, \modulo{\rho^{n}_{k-j-1/2}} \,
    \modulo{\eta_{k+1/2} - \eta_{k -1/2}}
    \\
    \nonumber
    & \leq &
    \sum_{k\in\interi} h \,
    \modulo{\rho^{n}_{k-j-1/2}}
    \int_{x_{k-1/2}}^{x_{k+1/2}} \modulo{\eta'(s)}  \d s
    \\
    \label{eq:c-c}
    & \leq &
    h \, \norma{\rho^n}_{\L1} \, \norma{\eta'}_{\L\infty} \,.
  \end{eqnarray}
  Similarly,
  \begin{eqnarray}
    \nonumber
    \modulo{c^n_{j+3/2} + 2c^n_{j+1/2} - c^n_{j-1/2}}
    & \leq &
    \sum_{k\in\interi} h \,
    \modulo{\rho^{n}_{k-j-1/2}} \, \modulo{\eta_{k-1/2} - 2\eta_{k +1/2} + \eta_{k +3/2}}
    \\
    \nonumber
    & \leq &
    h \, \sum_{k\in\interi} h \, \modulo{\rho^{n}_{k-j-1/2}}
    \, \modulo{\eta'(\zeta_{k+1}) - \eta'(\zeta_k)}
    \\
    \nonumber
    & = &
    h \sum_{k\in\interi} h \, \modulo{\rho^{n}_{k-j-1/2}}
    \int_{\zeta_{k}}^{\zeta_{k+1}} \modulo{\eta''(s)} \d s
    \\
    \nonumber
    & \leq &
    h \sum_{k\in\interi} h \, \modulo{\rho^{n}_{k-j-1/2}}
    \int_{x_{k-1/2}}^{x_{k+3/2}} \modulo{\eta''(s)} \d s
    \\
    \label{eq:c-c2}
    & = &
    2 \, h^2 \, \norma{\rho^n}_{\L1} \, \norma{\eta''}_{\L\infty} \,,
  \end{eqnarray}
  to complete the estimate of~\eqref{eq:uu1} we use the results above
  to bound the remaining terms
  \begin{eqnarray*}
    \modulo{\frac{f (t^n, x_{j+3/2}, \rho^n_j) - f (t^n, x_{j+1/2}, \rho^n_j)}{2}}
    & \leq &
    \frac{1}{2} \, h \, C \, \modulo{\rho^n_j}
    \\
    \modulo{\frac{v (c_{j+3/2}) - v (c_{j+1/2})}{2}}
    & \leq &
    \frac{1}{2} \, h \, \norma{v'}_{\L\infty} \, \norma{\eta'}_{\L\infty} \,
    \norma{\rho^o}_{\L1}
    \\
    \modulo{\frac{f (t^n, x_{j+1/2}, \rho^n_j) - f (t^n, x_{j-1/2}, \rho^n_j)}{2}}
    & \leq &
    \frac{1}{2} \, h \, C \, \modulo{\rho^n_j}
    \\
    \modulo{\frac{v (c_{j+1/2}) - v (c_{j-1/2})}{2}}
    & \leq &
    \frac{1}{2} \, h \, \norma{v'}_{\L\infty} \, \norma{\eta'}_{\L\infty} \,
    \norma{\rho^o}_{\L1}
  \end{eqnarray*}
  and we are now able to complete the estimate of~\eqref{eq:uu1}:
  \begin{eqnarray*}
    \eqref{eq:uu1}
    & \leq &
    h^2 \, C \, \norma{v}_\L\infty \, \modulo{\rho^n_j}
    \\
    & &
    +
    \norma{\partial_\rho f}_{\L\infty} \, \modulo{\rho^n_j}
    \left(
      h^2 \, \norma{v''}_{\L\infty} \norma{\rho^o}_{\L1}^2 \, \norma{\eta'}_{\L\infty}
      +
      h^2 \, \norma{v'}_{\L\infty} \, \norma{\rho^o}_{\L1}^2 \,
      \norma{\eta''}_{\L\infty}
    \right)
    \\
    & &
    +
    \frac{1}{2} \, h^2 \, C \, \norma{v'}_{\L\infty} \, \norma{\eta'}_{\L\infty} \,
    \norma{\rho^o}_{\L1} \, \modulo{\rho^n_j}
    \\
    & = &
    h^2
    \Big[
    C  \norma{v}_\L\infty
    +
    \Big(
    \norma{\partial_\rho f}_{\L\infty}
    \left(
      \norma{v''}_{\L\infty}   \norma{\eta'}_{\L\infty}
      +
      \norma{v'}_{\L\infty}   \norma{\eta''}_{\L\infty}
    \right)
    \\
    & &
    \qquad\qquad\qquad
    +
    \frac{C}{2} \norma{v'}_{\L\infty}  \norma{\eta'}_{\L\infty}
    \Big)
    \norma{\rho^o}_{\L1}
    \Big]
    \modulo{\rho^n_j}
    \\
    & \leq &
    h^2 \, \norma{v}_{\W2\infty}
    \left(
      C
      +
      \left(
        \norma{\partial_\rho f}_{\L\infty}
        +
        \frac{C}{2}
      \right)
      \norma{\eta'}_{\W1\infty} \, \norma{\rho^o}_{\L1}
    \right)
    \modulo{\rho^n_j}
  \end{eqnarray*}

  We now pass to estimate~\eqref{eq:uu2} and~\eqref{eq:uu3}:
  \begin{eqnarray*}
    \eqref{eq:uu2} + \eqref{eq:uu3}
    & = &
    \frac{1}{2} \,
    \left(f (t^n, x_{j+3/2}, \rho^n_{j+1}) - f (t^n, x_{j+1/2}, \rho^n_{j+1})\right)
    v (c_{j+3/2})
    \\
    & &
    +
    \frac{1}{2} \,
    f (t^n, x_{j+1/2}, \rho^n_{j+1}) \left(v (c_{j+3/2}) - v (c_{j+1/2})\right)
    \\
    & &
    -
    \frac{1}{2} \,
    \left(f (t^n, x_{j+1/2}, \rho^n_{j-1}) - f (t^n, x_{j-1/2}, \rho^n_{j-1})\right)
    v (c_{j+1/2})
    \\
    & &
    -
    \frac{1}{2} \,
    f (t^n, x_{j-1/2}, \rho^n_{j-1}) \left(v (c_{j+1/2}) - v (c_{j-1/2})\right)
    \\
    & = &
    \frac{1}{2} \, h \,
    \left(
      \partial_x f (t^n, \xi_{j+1}, \rho^n_{j+1}) \,   v (c_{j+3/2})
      -
      \partial_x f (t^n, \xi_{j}, \rho^n_{j-1}) \,   v (c_{j+1/2})
    \right)
    \\
    & &
    +
    \frac{1}{2}
    \Big(
    f (t^n, x_{j+1/2}, \rho^n_{j+1}) \, v' (\gamma_{j+1}) \,(c_{j+3/2} - c_{j+1/2})
    \\
    & &
    \qquad\qquad\qquad
    -
    f (t^n, x_{j-1/2}, \rho^n_{j-1}) \, v' (\gamma_{j}) \,(c_{j+1/2} - c_{j-1/2})
    \Big)
  \end{eqnarray*}
  for suitable $\xi_j \in \left]x_{j-1/2}, x_{j+1/2}\right[$ and
  $\gamma_j \in I (c_{j-1/2}, c_{j+1/2})$.  Introducing $\hat
  \xi_j \in \left]\xi_j, \xi_{j+1} \right[$, $\hat \zeta_j \in
  I(\rho^n_{j-1}, \rho^n_{j+1})$, $\hat \gamma_j \in I (\gamma_j,
  \gamma_{j+1})$, $\check \xi_j \in \left]x_{j-1/2}, x_{j+1/2}
  \right[$, $\check \zeta_j \in I(\rho^n_{j-1}, \rho^n_{j+1})$,
  $\check \gamma_j \in I (\gamma_j, \gamma_{j+1})$, $\delta_j \in I
  (c_{j+3/2}-c_{j+1/2},c_{j+1/2}-c_{j-1/2})$ and
  using~\eqref{eq:hyp1}, \eqref{eq:c-c}, \eqref{eq:c-c2}
  \begin{eqnarray*}
    \!\!\!
    & &
    \modulo{\eqref{eq:uu2} + \eqref{eq:uu3}}
    \\
    \!\!\!
    & \leq &
    \frac{1}{2} \, h \,
    \Big(
    \modulo{\partial^2_{x x} f (t^n, \hat \xi_j, \hat \zeta_j)}
    \, v (\hat \gamma_{j+1})
    \, h
    +
    \modulo{\partial^2_{\rho x}f (t^n, \hat \xi_j, \hat \zeta_j)}
    \, v (\hat \gamma_{j+1})
    \, \modulo{\rho^n_{j+1} - \rho^n_{j-1}}
    \\
    \!\!\!
    & &
    \qquad\qquad\qquad
    +
    \modulo{\partial_x f (t^n, \hat \xi_j, \hat \zeta_j)}
    \, \modulo{v' (\hat \gamma_{j+1})} \, \modulo{c_{j+3/2} - c_{j+1/2}}
    \Big)
    \\
    \!\!\!
    & &
    +
    \frac{1}{2} \,
    \Big(
    \modulo{\partial_x f (t^n, \check \xi_j, \check \zeta_j)}
    \, \modulo{v' (\check \gamma_j)} \, \modulo{\delta_{j}} \, h
    +
    \modulo{f (t^n, \check \xi_j, \check \zeta_j)}
    \, \modulo{v'' (\check \gamma_j)} \, \modulo{\delta_{j}}
    \, \modulo{\gamma_{j+1} - \gamma_j}
    \\
    \!\!\!
    & &
    \qquad\qquad\qquad
    +
    \modulo{f (t^n, \check \xi_j, \check \zeta_j)}
    \, \modulo{v'(\check \gamma_j)} \, \modulo{c_{j+3/2} - 2c_{j+1/2} - c_{j-1/2}}
    \Big)
    \\
    \!\!\!
    & \leq &
    \frac{1}{2} \, h
    \Big(
    C \, \norma{v}_{\L\infty} \, \modulo{\hat \zeta_j} \, h
    +
    \norma{\partial^2_{\rho x}f}_{\L\infty} \, \norma{v}_{\L\infty} \,
    \modulo{\rho^n_{j+1} - \rho^n_{j-1}}
    +
    C \,
    \norma{v'}_{\L\infty} \,
    \norma{\rho^o}_{\L1} \,
    \norma{\eta'}_{\L\infty} \,
    \modulo{\hat \zeta_j} \, h
    \Big)
    \\
    \!\!\!
    & &
    +
    \frac{1}{2} \,
    \Big(
    C \,
    \norma{v'}_{\L\infty} \,
    \norma{\rho^o}_{\L1} \,
    \norma{\eta'}_{\L\infty} \,
    \modulo{\check \zeta_j} \, h^2
    +
    \norma{\partial_\rho f}_{\L\infty} \,
    \norma{v''}_{\L\infty} \,
    \norma{\rho^o}_{\L1}^2 \, \norma{\eta'}_{\L\infty}^2 \,
    \modulo{\check \zeta_j} \,h^2
    \\
    \!\!\!
    & &
    \qquad\qquad\qquad
    +
    2 \, \norma{\partial_\rho f}_{\L\infty}  \,
    \norma{v'}_{\L\infty} \,
    \norma{\rho^o}_{\L1} \, \norma{\eta''}_{\L\infty} \,
    \modulo{\check \zeta_j} \, h^2
    \Big)
    \\
    \!\!\!
    & = &
    \frac{1}{2} \,
    C
    \left(
      1
      +
      \norma{\eta'}_{\L\infty} \, \norma{\rho^o}_{\L1}
    \right)
    \norma{v}_{\W1\infty} \,
    \modulo{\hat \zeta_j}
    \, h^2
    \\
    \!\!\!
    & &
    +
    \frac{1}{2}
    \left(
      C \,
      +
      \norma{\partial_\rho f}_{\L\infty} \,
      \norma{\rho^o}_{\L1} \, \norma{\eta'}_{\L\infty}
      +
      2 \, \norma{\partial_\rho f}_{\L\infty}
    \right)
    \norma{v'}_{\W1\infty} \,
    \norma{\eta'}_{\W1\infty} \,
    \norma{\rho^o}_{\L1} \,
    \modulo{\check \zeta_j} \, h^2
    \\
    \!\!\!
    & &
    +
    \frac{1}{2} \,
    \norma{\partial^2_{\rho x}f}_{\L\infty} \, \norma{v}_{\L\infty}
    \modulo{\rho^n_{j+1} - \rho^n_{j-1}} \, h \,.
  \end{eqnarray*}
  The above bound allows to obtain the estimate for $\mathcal{B}_j$:
  \begin{eqnarray*}
    \mathcal{B}_j
    & \leq &
    \norma{v}_{\W2\infty}
    \left(
      C
      +
      \left(
        \norma{\partial_\rho f}_{\L\infty}
        +
        \frac{C}{2}
      \right)
      \norma{\eta'}_{\W1\infty} \, \norma{\rho^o}_{\L1}
    \right)
    \modulo{\rho^n_j} \, h^2
    \\
    & &
    +
    \frac{1}{2} \,
    C
    \left(
      1
      +
      \norma{\eta'}_{\L\infty} \, \norma{\rho^o}_{\L1}
    \right)
    \norma{v}_{\W1\infty} \,
    \modulo{\hat \zeta_j}
    \, h^2
    \\
    & &
    +
    \frac{1}{2}
    \left(
      C \,
      +
      \norma{\partial_\rho f}_{\L\infty} \,
      \norma{\rho^o}_{\L1} \, \norma{\eta'}_{\L\infty}
      +
      2 \, \norma{\partial_\rho f}_{\L\infty}
    \right)
    \norma{v'}_{\W1\infty} \,
    \norma{\eta'}_{\W1\infty} \,
    \norma{\rho^o}_{\L1} \,
    \modulo{\check \zeta_j} \, h^2
    \\
    & &
    +
    \frac{1}{2} \,
    \norma{\partial^2_{\rho x}f}_{\L\infty} \, \norma{v}_{\L\infty}
    \modulo{\rho^n_{j+1} - \rho^n_{j-1}} \, h \,.
  \end{eqnarray*}
  so that
  \begin{eqnarray*}
    \sum_{j \in \interi} \mathcal{B}_j
    & \leq &
    \norma{v}_{\W2\infty}
    \left(
      C
      +
      \left(
        \norma{\partial_\rho f}_{\L\infty}
        +
        \frac{C}{2}
      \right)
      \norma{\eta'}_{\W1\infty} \, \norma{\rho^o}_{\L1}
    \right)
    \norma{\rho^o}_{\L1} \, h
    \\
    & &
    +
    \frac{1}{2} \,
    C
    \left(
      1
      +
      \norma{\eta'}_{\L\infty} \, \norma{\rho^o}_{\L1}
    \right)
    \norma{v}_{\W1\infty} \,
    \norma{\rho^o}_{\L1} \, h
    \\
    & &
    +
    \frac{1}{2}
    \left(
      C \,
      +
      \norma{\partial_\rho f}_{\L\infty}
      \left(2 + \norma{\rho^o}_{\L1} \, \norma{\eta'}_{\L\infty}\right)
    \right)
    \norma{v'}_{\W1\infty} \,
    \norma{\eta'}_{\W1\infty} \,
    \norma{\rho^o}_{\L1} \, h
    \\
    & &
    +
    \frac{1}{2} \,
    \norma{\partial^2_{\rho x}f}_{\L\infty} \, \norma{v}_{\L\infty}
    \left(\sum_{j \in \interi} \modulo{\rho^n_{j+1} - \rho^n_{j-1}}\right) h \,.
    % \\
    % & \leq &
    % \left(\mbox{horrible constant}\right) \, \norma{\rho^o}_{\L1} \,
    % h
    % +
    % \left(\mbox{less horrible constant}\right)
    % \left(\sum_{j \in \interi} \modulo{\rho^n_{j+1} -
    %   \rho^n_{j-1}}\right) h \,.
  \end{eqnarray*}
  Recall now~\eqref{eq:AB} and~\eqref{eq:90} to obtain
  \begin{displaymath}
    \sum_{j \in \interi} \modulo{\rho^{n+1}_{j+1} - \rho^{n+1}_j}
    \leq
    \sum_{j \in \interi}
    \modulo{\mathcal{A}_j}
    +
    \lambda \sum_{j \in \interi} \modulo{\mathcal{B}_j}
    \leq
    (1+\mathcal{K}_1 \, \tau) \sum_{j \in \interi} \modulo{\rho^n_{j+1} - \rho^n_j}
    +
    \mathcal{K}_2 \, \tau
  \end{displaymath}
  where
  \begin{equation}
    \label{eq:K1K2}
    \begin{array}{rcl}
      \mathcal{K}_1
      & = &
      \frac12\, \norma{\partial_\rho f}_{\L\infty} \, \norma{v'}_{\L\infty}  \, \norma{\rho^o}_{\L1} \, \norma{\eta'}_{\L\infty}
      +
      \norma{\partial^2_{\rho x}f}_{\L\infty} \, \norma{v}_{\L\infty}
      \\
      \mathcal{K}_2
      & = &
      \Bigg[
      \frac{3}{2} \, C
      +
      \left(
        \norma{\partial_\rho f}_{\L\infty}
        +
        C
      \right)
      \norma{\eta'}_{\W1\infty} \, \norma{\rho^o}_{\L1}
      \\
      & &
      +
      \frac{1}{2}
      \left(
        C \,
        +
        \norma{\partial_\rho f}_{\L\infty}
        \left(2 + \norma{\rho^o}_{\L1} \, \norma{\eta'}_{\L\infty}\right)
      \right)
      \norma{\eta'}_{\W1\infty}
      \Bigg]
      \norma{v}_{\W2\infty} \, \norma{\rho^o}_{\L1}
    \end{array}
  \end{equation}
  The estimate~\eqref{TV1} now follows from standard iterative
  procedure.  The proof of Proposition~\ref{prop:TV} follows
  immediately.
\end{proofof}

\begin{proofof}{Lemma~\ref{lem:DepOnT}}
  We follow the same line as
  in~\cite[Section~3]{KarlsenRisebro2001}. Using~(\ref{eq:c-c}),
  Lemma~\ref{lem:pos}, Lemma~\ref{lem:L1} and
  Proposition~\ref{prop:TV} compute preliminarily
  \begin{eqnarray*}
    & &
    \sum_{j \in \interi}
    \modulo{D^+ \! \left(f (t^n, x_{j-1/2}, \rho^n_j) \, v (c^n_{j-1/2})\right)}
    \\
    & \leq &
    \sum_{j \in \interi}
    \Big[
    \modulo{\partial_x f (t^n, \xi_j, \zeta_j) \, v (\gamma_j)} \, h
    +
    \modulo{\partial_\rho f (t^n, \xi_j, \zeta_j) \, v (\gamma_j)}
    \, \modulo{\rho^n_{j+1}- \rho^n_j}
    \\
    & &
    \qquad\qquad
    +
    \modulo{f (t^n, \xi_j, \zeta_j) \, v' (\gamma_j)}
    \, \modulo{c^n_{j+1/2} - c^n_{j-1/2}}
    \Big]
    \\
    & \leq &
    \sum_{j \in \interi}
    \Big[
    C \, \norma{v}_{\L\infty} \, \rho^n_j \, h
    +
    \norma{\partial_\rho f}_{\L\infty}\, \norma{v}_{\L\infty} \,
    \modulo{\rho^n_{j+1} - \rho^n_j}
    \\
    & &
    \qquad\qquad
    +
    \norma{\partial_\rho f}_{\L\infty} \,
    \norma{v'}_{\L\infty} \, \norma{\eta'}_{\L\infty} \,
    \norma{\rho^o}_{\L1} \, \max\left\{\rho^n_j, \rho^n_{j+1}\right\} \,h
    \Big]
    \\
    & \leq &
    C
    \norma{v}_{\L\infty} \,
    \norma{\rho^o}_{\L1}
    +
    2 \norma{\partial_\rho f}_{\L\infty} \, \norma{\rho^o}_{\L1}^2 \,
    \norma{\eta'}_{\L\infty}
    \norma{v'}_{\L\infty} \,
    +
    \norma{\partial_\rho f}_{\L\infty} \norma{v}_{\L\infty}
    \sum_{j \in \interi} \modulo{\rho^n_{j+1} - \rho^n_j}
    \\
    & \leq &
    C
    \norma{v}_{\L\infty} \,
    \norma{\rho^o}_{\L1}
    +
    2 \norma{\partial_\rho f}_{\L\infty} \, \norma{\rho^o}_{\L1}^2 \,
    \norma{\eta'}_{\L\infty}
    \norma{v'}_{\L\infty}
    \\
    & &
    \qquad\qquad
    +
    \norma{\partial_\rho f}_{\L\infty} \norma{v}_{\L\infty}
    \left(
      \mathcal{K}_2 \, t
      +
      \sum_{j\in\interi} \modulo{\rho^o_{j+1} - \rho^o_j}
    \right)
    e^{\mathcal{K}_1 t}\,.
  \end{eqnarray*}
  The term $\sum_{j\in \interi} \modulo{D^- \! \left(f (t^n,
      x_{j+1/2}, \rho^n_j) \, v (c^n_{j+1/2})\right)}$ admits an
  analogous estimate.  Moreover,
  \begin{displaymath}
    \sum_{j \in \interi} \modulo{D^2 \rho^n_j}
    \leq
    2 \sum_{j \in \interi} \modulo{\rho^n_{j+1} - \rho^n_j}
    \leq
    2
    \left(
      \mathcal{K}_2 \, t
      +
      \sum_{j\in\interi} \modulo{\rho^o_{j+1} - \rho^o_j}
    \right)
    e^{\mathcal{K}_1 t} \,.
  \end{displaymath}
  Using the above estimates and~(\ref{LF}) we get
  \begin{eqnarray*}
    \norma{\rho^{n+1} - \rho^n}_{\L1}
    & = &
    \sum_{j\in \interi} h \, \modulo{\rho^{n+1}_j - \rho^n_j}
    \\
    & \leq&
    \frac{\tau}{2}
    \sum_{j\in \interi}
    \modulo{D^+ \! \left(f (t^n, x_{j-1/2}, \rho^n_j) \, v (c^n_{j-1/2})\right)}
    \\
    & &
    +
    \frac{\tau}{2}
    \sum_{j\in \interi}
    \modulo{D^- \! \left(f (t^n, x_{j+1/2}, \rho^n_j) \, v (c^n_{j+1/2})\right)}
    +
    \frac{\lambda \, \tau}{6} \sum_{j \in \interi} \modulo{D^2 \rho^n_j}
    \\
    & \leq & {\cal C} (t) \, \tau
  \end{eqnarray*}
  where
  \begin{equation}
    \label{eq:C}
    \begin{array}{rcl}
      {\cal C} (t)
      & = &
	 C
      \norma{v}_{\L\infty} \,
      \norma{\rho^o}_{\L1}
      +
      2 \norma{\partial_\rho f}_{\L\infty} \, \norma{\rho^o}_{\L1}^2 \,
      \norma{\eta'}_{\L\infty}
    \norma{v'}_{\L\infty} \,
      \\
      & &
      +
      \left(\norma{\partial_\rho f}_{\L\infty} \, \norma{v}_{\L\infty}
        + \frac{\lambda}{3}\right)
      \left(
        \mathcal{K}_2 \, t
        +
        \sum_{j\in\interi} \modulo{\rho^o_{j+1} - \rho^o_j}
      \right)
      e^{\mathcal{K}_1 t} \,,
    \end{array}
  \end{equation}
  completing the proof.
\end{proofof}

\begin{proofof}{Proposition~\ref{prop:Entropy}}
  Fix $n \in \naturali$ and for any sequence $(\rho)_{j\in\interi}$
  define the transformation $\rho \mapsto H (\rho)$ given by
  \begin{equation}
    \label{HH}
    H^n_j (\rho) = \rho_j -   \left(
      \mathbf{f}^n_{j+1/2}(\rho_j, \rho_{j+1})
      -
      \mathbf{f}^n_{j-1/2}(\rho_{j-1}, \rho_{j})
    \right),
  \end{equation}
  where the functions $\mathbf{f}^n_{j+1/2}$ are given by~\eqref{LFa},
  but where, instead of~\eqref{eq:cnj}, the sequence $(c^n_{j+1/2})_{j
    \in \interi}$ is now an arbitrary fixed sequence. Thus,
  $H^n_j(\rho)$ depends only on $\rho_{j-1}$, $\rho_j$ and
  $\rho_{j+1}$.  Then, $H^n$ is monotone, in the sense that
  \begin{equation}
    \label{Hmonot}
    \frac{\partial H^n_j}{\partial \rho_i} \geq 0 \,,
    \quad i = j-1, j, j+1.
  \end{equation}
  The cases $i = j\pm1$ are easily verified. If $i = j$,
  using~\eqref{LFa} we find
  \begin{eqnarray*}
    \frac{\partial H^n_j}{\partial \rho_j}
    & = &
    \frac13
    -
    \frac{\lambda}{2}
    \left(
      \partial_\rho f (t^n, x_{j+1/2}, \rho_{j} ) \, v(c_{j+1/2})
      -
      \partial_\rho f (t^n, x_{j-1/2}, \rho_{j} ) \, v(c_{j-1/2})
    \right)
    \\
    & \geq &
    \frac13
    -
    \lambda \, \norma{\partial_\rho f}_{\L\infty} \, \norma{v}_{\L\infty}
    \\
    & \geq &
    0
  \end{eqnarray*}
  by the CFL condition~\eqref{eq:CFL}. The definition~\eqref{HH} of
  $H^n$ and~\eqref{eq:Kruzhkov} imply that for any $k \in \reali$
  \begin{equation}
    \label{H10}
    \modulo{\rho_{j} - k}
    -
    \lambda
    \left(
      F^k_{j+1/2}(\rho_{j}, \rho_{j+1})
      -
      F^k_{j-1/2} (\rho_{j-1}, \rho_{j})
    \right)
    =
    H^n_j (\rho \wedge k) - H^n_j (\rho \vee k) \,,
  \end{equation}
  where $k$ in the right-hand side above is understood as the sequence
  identically equal to $k$. The monotonicity condition~\eqref{Hmonot}
  and the scheme~\eqref{LF}--\eqref{LFa} ensure that
  \begin{eqnarray}
    \nonumber
    & &
    H^n_j (\rho \wedge k) - H^n_j (\rho \vee k)
    \\
    \nonumber
    & \geq &
    H^n_j(\rho) \wedge H^n_j ( k ) - H^n_j (\rho) \vee H^n_j ( k )
    \\
    \nonumber
    & = &
    \sgn \left[
      H^n_j (\rho)
      -
      k
      +
      \lambda \left( f(t^n, x_{j+1/2}, k) - f(t^n, x_{j-1/2}, k) \right)
    \right]
    \times
    \\
    \nonumber
    & &
    \qquad
    \times
    \left[
      H^n_j (\rho)
      -
      k
      +
      \lambda \left( f(t^n, x_{j+1/2}, k) - f(t^n, x_{j-1/2}, k) \right)
    \right]
    \\
    \nonumber
    & \geq &
    \sgn \left( H^n_j (\rho) - k \right)
    \left[
      H^n_j (\rho)
      -
      k
      +
      \lambda \left( f(t^n, x_{j+1/2}, k) - f(t^n, x_{j-1/2}, k) \right)
    \right]
    \\
    \label{H20}
    & = &
    \modulo{H^n_j (\rho) - k}
    +
    \lambda \,
    \sgn \left( H^n_j(\rho) - k\right)
    \left( f(t^n, x_{j+1/2}, k) - f(t^n, x_{j-1/2}, k) \right) .
  \end{eqnarray}
  In the last inequality we used also the non-negativity of the
  function $(a,b) \mapsto (\sgn(a+b) - \sgn(a) )
  (a+b)$. From~\eqref{H10} and~\eqref{H20} we conclude that
  \begin{equation}
    \label{H30}
    \begin{array}{rcl}
      \modulo{H^n_j(\rho) - k}
      -
      \modulo{\rho_{j} - k}
      +
      \lambda
      \left(
        F^k_{j+1/2} (\rho_{j}, \rho_{j+1})
        -
        F^k_{j-1/2} (\rho_{j-1},\rho_{j})
      \right)
      \\
      \qquad +
      \lambda \sgn ( H^n_j(\rho) - k)
      \left( f(t^n, x_{j+1/2}, k)- f(t^n, x_{j-1/2}, k) \right)
      & \leq &
      0 \,.
    \end{array}
  \end{equation}

  Consider now the numerical approximation $\rho^{n}_{j}$ given by the
  algorithm~\eqref{LF}. Then, we apply~\eqref{H30} to $\rho^{n}$, with
  the sequence $c_{j+1/2}$ appearing in~\eqref{HH} as given by the
  convolution~\eqref{eq:cnj}. Observing that $H^n_j(\rho^n) =
  \rho^{n+1}_{j}$, we conclude that~\eqref{eq:Entropy} holds.
\end{proofof}

\small{

  \bibliography{ACT}

  \bibliographystyle{abbrv}

}

\end{document}